\documentclass[times,review]{elsarticle}
\usepackage{jcomp}
\usepackage{framed,multirow}

\usepackage{latexsym}

\usepackage{url}
\usepackage{xcolor}
\definecolor{newcolor}{rgb}{.8,.349,.1}

\usepackage{ae}
\usepackage{amsfonts}
\usepackage{amsmath}
\usepackage{amsopn}
\usepackage{amssymb}
\usepackage{array}
\usepackage{booktabs}
\usepackage[hang,small,bf]{caption}
\usepackage{color}
\usepackage{colortbl}
\usepackage{enumitem}
\usepackage{epsfig} % For figure environment
\usepackage{epstopdf}
\usepackage{etoolbox}
\usepackage{float}
\usepackage[T1]{fontenc}
\usepackage{graphics}
\usepackage{graphicx}
\usepackage[utf8x]{inputenc}
\usepackage{lipsum}
\usepackage{multirow}
\usepackage{parskip}
\usepackage{pstricks}
\usepackage{pstricks-add} % PSTricks
\usepackage{pst-grad, pst-plot, pst-eps}
\usepackage{setspace}
\usepackage[euler]{textgreek}
\usepackage{xcolor}
\usepackage[caption=false,labelformat=simple]{subfig}

\definecolor{backgr}{rgb}{0.1,0.1,0.3}
\definecolor{darkblue}{rgb}{0,0,0.6}
\definecolor{darkgreen}{rgb}{0,0.6,0}
\definecolor{darkred}{rgb}{0.6,0,0}
\definecolor{orange}{rgb}{1,0.34,0}
\definecolor{darkorange}{rgb}{1,0.4,0.3}
\definecolor{MatBlue}{rgb}{0,0.4431,0.7373}
\definecolor{MatOrange}{rgb}{0.9255,0.6902,0.1216}
\definecolor{MatRed}{rgb}{0.8471,0.3216,0.0941}

\newcommand{\todo}[1]{ {\color{red}{! #1 !}}}

\newcommand{\dontshow}[1]{}

\DeclareMathOperator{\R}{\mathbb{R}}

\newcommand{\Ft}[1]{\left(#1\right)^\wedge}

\newcommand{\Pat}{\mathcal P}
\def\kwave/{\textbf{k}-\texttt{Wave}}
\def\restrict#1{\raise-.5ex\hbox{\ensuremath|}_{#1}}
\newcommand{\D}{\textrm{d}}
\newcommand{\Marta}[1]{{\textcolor{blue}{#1}}}

\newcommand{\Kiko}[1]{{\textcolor{darkgreen}{#1}}}

\definecolor{myGreen}{rgb}{0, 0, 0}
\definecolor{myBlue}{rgb}{0, 0. 0}

\newcommand{\rA}[1]{\textcolor{myGreen}{#1}}
\newcommand{\rB}[1]{\textcolor{myBlue}{#1}}

\ifpdf
  \DeclareGraphicsExtensions{.eps,.pdf,.png,.jpg}
\else
  \DeclareGraphicsExtensions{.eps}
\fi

\verso{Francesc Rul$\cdot$lan \textit{etal}}
\title{Hamilton-Green Solver for the Forward and Adjoint Problems in Photoacoustic Tomography}%
\author[1]{Francesc \snm{Rul$\cdot$lan}}

\author[1]{Marta M.~\snm{Betcke}\corref{cor1}}
\cortext[cor1]{Corresponding author: 
  \ead{m.betcke@ucl.ac.uk}
  }

\address[1]{Department of Computer Science, University College London, London WC1E 6BT, UK}
\begin{document}

\begin{frontmatter}

%==============================
% Abstract
%==============================
% REQUIRED
\begin{abstract}                                               
The majority of the solvers for the acoustic problem in Photoacoustic Tomography (PAT) rely on  full solution of the wave equation, which makes them less suitable for \textit{real-time} and \textit{dynamic} applications where only partial data is available.
This is in contrast to other tomographic modalities, e.g. X-ray tomography, where partial data implies partial cost for the application of the forward and adjoint operators.
In this work we present a novel solver for the forward and adjoint wave equations for the
acoustic problem in PAT. We term the proposed solver Hamilton-Green as it approximates the fundamental solution to the respective wave equation along the trajectories of the Hamiltonian system resulting from the high frequency \rB{asymptotic approximate solution} for the wave equation.
This approach is flexible and scalable in the sense that it allows computing the solution for each sensor independently at a fraction of the cost of the full wave solution. 
The theoretical foundations of our approach are rooted in results available in seismics and ocean acoustics. To demonstrate the feasibility of our approach we present results for 2D domains with homogeneous and heterogeneous sound speeds and evaluate them against a full wave solution obtained with a pseudospectral finite difference method implemented in the \texttt{k-Wave} toolbox \cite{treeby2010k}.
%We also include a 2D phantom simulation using heterogeneous sound speed that demonstrates the feasibility of our approach.
\end{abstract}


\end{frontmatter}

%==============================
% 1. Introduction
%==============================
%%  % REQUIRED
%%  \begin{AMS}
%%    35L05, %\KikoCom{Partial differential equations - Wave equation\\}
%%    65M80, %\KikoCom{Numerical analysis - Partial differential equations, initial value and time-dependent initial-boundary value problems -  Fundamental solutions, Green's function methods, etc.\\}
%%    65M32  %\KikoCom{Numerical analysis - Inverse problems}
%%  \end{AMS}

%==============================
% Introduction
%==============================
\section{Photoacoustic tomography}
\label{sec:intro}
%================================================================================
%====================   PHOTOACOUSTIC TOMOGRAPHY   ==============================
%================================================================================

Photoacoustic tomography (PAT) is a hybrid imaging technique based on the photoacoustic effect. 
In PAT, the biological tissue is irradiated with a laser pulse \cite{xu2006photoacoustic, kuchment2011mathematics}. 
Part of the energy of this pulse is absorbed in the form of heat, which is then released through thermoelastic expansion and emission of an ultrasonic wave. % associated to this expansion. 
This wave propagates outwards, where a set of detectors placed at the boundary of the tissue records the ultrasonic pressure over time.
This data is then used to reconstruct the initial pressure (PAT image). 
The main advantage of this modality is its ability to simultaneously obtain high spatial resolution of the ultrasonic wave and the high contrast of the optical absorption of the laser pulse. 
%In optical imaging techniques, the spatial resolution is greatly degraded as we increase the depth, mainly due to optical scattering.
%Additionally, the poor contrast obtained via ultrasound techniques, as a consequence of the acoustic properties of the biological tissue, is a limiting factor in terms of its use
%as diagnosis tools. 
%In PAT the objective is to image the absorption of the tissue in the near infrared range (high contrast and relatively high penetration)
%via the emitted associated ultrasound (high spatial resolution).
%
%As a consequence of this hybrid nature of PAT, in order to reconstruct the images we need two physical models. 
%The first one addresses the electromagnetic propagation of the laser from the source towards and inside the tissue. 
%The mathematical tools used to solve this problem are based on the radiative transfer equation \cite{cox2009challenges, ren2013hybrid}. 
%The second model is related to the acoustic propagation induced by the generation of the ultrasonic wave after the thermoelastic expansion.
To make PAT quantitative requires a solution of a coupled acoustic and optical problem, where the initial acoustic pressure constitutes the ``interior data'' for the optical problem \cite{cox2009challenges}. In this work we focus on the acoustic problem. %, which is modeled by a linear scalar wave equation \cite{morse1968theoretical}.
%The PAT image is obtained by concatenating the solutions to the acoustic problem and the optical problem. The acoustic problem is the subject of this work. This is modeled by the linear scalar wave equation, derived from the conservation laws of mass and momentum \cite{morse1968theoretical}.

\subsection{Acoustic propagation in PAT} 
The forward problem in PAT is modeled by the initial value problem for the wave equation in free space $\mathbb{R}^d$ see e.g.~
\cite{kuchment2011mathematics, arridge2016adjoint}
%The propagation of the pressure wave inside a domain $\mathbb{R}^d$ is modeled by
\begin{subequations}\label{eq:PATfwd}%\label{problem}
\begin{align}
    \square^2 u(t,x) := \frac{1}{c^2(x)} \frac{\partial^2 u(t,x)}{\partial t^2} - \Delta u(t,x) & = 0, & (t, x)\in (0, \infty)\times\mathbb{R}^d, \label{eq:PATfwd-1}\\[5pt]
    u_t(0, x)              & = 0,  & x\in\mathbb{R}^d, \label{eq:PATfwd-2}\\[5pt]  
    u(0, x)                & = u_0(x), & x\in\mathbb{R}^d, \label{eq:PATfwd-3}
\end{align}
\end{subequations}
where $c(x)\in\mathcal{C}^\infty(\mathbb{R}^d)$ denotes the sound speed in $\mathbb{R}^d$.
Reconstruction of the PAT image amounts to recovery of the initial pressure $u(0,x) = u_0(x)$
from pressure over time measurements at a set of ideal (isotropic and point-like) detectors $\Gamma = \{ x_0^m \in \R^d, \; m=1,\dots,M\}$.
\rB{The detectors are typically placed on the boundary of the object, covering the boundary fully or, more frequently, in part, the latter scenario potentially resulting in a limited view problem. Other factors like limited propagation time or trapping sound speed may also result in incomplete data.}
%the tissue, hence limiting the region of interest to the domain circumscribed by them. In case of using partially surrounding sensors, the domain is then limited by other factors such as the propagation time or the depth penetrated by the laser beam.}
We denote the data obtained at the sensors over a finite measuring time $T > 0$ with 
\begin{equation}\label{sensorData}
g(t, y) = u(t, y), \quad (t, y) \in (0, T)\times \Gamma.
\end{equation}
%and $\omega: C^\infty(0,T) \times \Gamma \rightarrow C_0^\infty(0,T) \times \Gamma$ %is an \textit{aperture function} which restricts the data to the sensor support and applies smooth cut off to the measurement time.
%$\omega(t, y) \in C_0^{\infty}(0,T) \times \Gamma$ the \textit{aperture} function that restricts the field of view to the set of detectors and applies smooth cut off to the measurement time.
%\begin{definition}
We define the \textit{forward operator $\mathcal P$} for PAT with sensors placed on $\Gamma$ and the measurement time $T$ as
\begin{align*} 
\mathcal P: \mathcal{C}_0^\infty(\R^d) & \longrightarrow \mathcal{C}_0^\infty(0, T)\times \Gamma\\
                         u_0(x)& \longmapsto   \omega(t)  g(t, y), 
\end{align*}
with $g(t, y)$ given by \eqref{sensorData}
and \rB{an aperture function $\omega(t)\in C^\infty_0(0,T)$} %is an \textit{aperture function} which restricts the data to the sensor support and applies smooth cut off to the measurement time.
which applies a smooth cut off to the measurement time.
%\end{definition}
The corresponding adjoint \textit{operator $\mathcal{P}^*$} is defined as
\begin{align*} 
\mathcal{P}^*: \mathcal{C}_0^\infty(0, T)\times \Gamma  & \longrightarrow \mathcal{C}_0^\infty(\R^d)\\
                                    g(t, y) & \longmapsto      v(x, T),
\end{align*}
where $v$ is the solution to the wave equation with a time varying source supported on $(0,T) \times \Gamma$ 
\begin{subequations}\label{eq:PATadj}
\begin{align}%\label{timeReversalProblem}
\label{eq:PATadj-1} %\frac{1}{c^2(x)} \frac{\partial^2 v(t,x)}{\partial t^2} - \nabla^2 v(t,x)   =  
\square^2 v(t,x) &=
%& \left\{\begin{array}{l} \sum_{m=1}^M \frac{\partial}{\partial t} \left( g(T-t, x_0^m)) \omega(T-t) \right) \delta (x-x_0^m) \\[5pt]
%                                                        0  \end{array}\right. & 
%                                       \begin{array}{l} (t, x)\in(0,T) \times \R^d, \\[5pt]
%                                                      \textnormal{everywhere else,}\end{array} \\[5pt]
\sum_{m=1}^M \frac{\partial}{\partial t} \Big( g(T-t, x_0^m)\, \omega(T-t) \Big) \delta (x-x_0^m), \quad  & (t, x) \in(0,T) \times \R^d, \\[5pt]
\label{eq:PATadj-2} v_t(0, x)                &=   0, & x\in\mathbb{R}^d,\\[5pt]
\label{eq:PATadj-3} v(0, x)                  &=   0, & x\in\mathbb{R}^d,
\end{align}
\end{subequations}
\noindent evaluated at the time $t=T$, $v(T,x)$ \cite{arridge2016adjoint, belhachmi2016direct}.
\rA{Thus, in the adjoint problem the time reversed sensor measurements $g(T-t, x_0^m)\,\omega(T-t),\, m=1,\dots,M$ act as time dependent point mass sources. %with magnitudes $\frac{\partial}{\partial t} \left( g(T-t, x_0^m)\,\omega(T-t) \right), \; m=1,\dots,M$.
%The inverse problem in PAT is to recover the initial pressure $u_0(x)$ from a pressure time series $g(t, x)\,\omega(t, x)$ measured over time time $T$ on the set of sensors $\Gamma$. 
In contrast, a popular technique for PAT inversion known as \textit{time reversal} consists of enforcing the same time reversed measurements $g(T-t, x_0^m)\,\omega(T-t),\, m=1,\dots,M$ as constraints, see e.g.~\cite{kuchment2011mathematics,stefanov2009thermoacoustic}, and yields an approximate-inverse operator to $\mathcal P$ \cite{givoli2014time, xu2004application}.
The iterated variant of time reversal proposed in \cite{qian2011efficient} effectively restores partially visible singularities which arise in limited view problems (e.g.~trapping sound speed, limited time, limited angle). However, when the measured data is sparse, i.e.~relatively few sensors are used, we increasingly encounter invisible singularities. Variational methods provide a simple way of incorporating prior knowledge to fill in the missing information in the data.  
%However, when the measured data is sparse i.e.~relatively few sensors are used, the approximation quality decays which may negatively impact the convergence of the Neumann series (e.g.~as in iterative time reversal). \KikoCom{Need to include citation of \cite{qian2011efficient} with explanation}.
Moreover, the adjoint operator to a sparse data problem can be directly obtained and used in the framework of variational methods with readily available convergence results. 
}
Finally, we would like to remark that we altered the functional framework introduced in \cite{arridge2016adjoint} to suit the proposed ray based approach. 
The adjoint operator from \cite{arridge2016adjoint} follows analogously where all the spatial integrals are evaluated in the sense of the definition of the Green's function
$$\square^2 G(t,x \mid t_0, x_0) = \delta (x-x_0)\, \delta(t-t_0),$$
which is the framework we adapt in this paper. 

%================================================================================
%====================   MOTIVATION FOR HF IN PAT   ==============================
%================================================================================

\rB{\subsection{Motivation}\label{sec:intro:motiv}}

The complexity of the full wave solvers such as
%\textcolor{green}{include k-wave, not sure about BEM (smigaj2015solving) - this needs homogeneous speed of sound and has lower complexity because it forms boundary integral - I would limit to time stepping schemes if not spectral they will have $n^4 k$ complexity}
\cite{holm2001ultrasim, gerstoft2004cabrillo} is independent from the number of sensors. For instance, for the pseudospectral finite difference solver on a cube-shaped domain of size $n^3$ we have a complexity of $\mathcal{O}(n^4 \log n)$.
The analogous scenario for the here proposed Hamilton-Green (HG) solver corresponds to the complete data case (with some abuse of notation as the sensor is finite) i.e.~collecting $n^2$ measurements on one side of the cube, which results in the complexity of $\mathcal{O}(n^5)$. Thus, at the first glance, our new method has a worse theoretical complexity than the existing full wave solvers. However, the Hamilton-Green solver also offers a number of advantages with respect to the full wave solvers which could make it competitive in scenarios with e.g.~incomplete / sparse data, various region of interest problems and in conjunction with ultrasound tomography: 
\begin{itemize}\itemsep -0.2cm
    \item Solution for each sensor is computed independently at a much lower cost, i.e.~$\mathcal{O}(n^3)$, than the full wave solution; 
    \item In contrast to the finite difference or pseudospectral solvers, it can take advantage of irregular or multiscale grids;
    \item The solution can be computed on a part of the domain (a region of interest) at a largely reduced computational cost;
    \item The coupling of the homogeneous and inhomogeneous parts of the domain is inherent and efficiently accomplished via the step size choice;
    \item The Hamiltonian system underlying the trajectories %used in the approximation in 
    computed by the Hamilton-Green solver is the one of the transmission ultrasound problem, rendering the PAT solution essentially free when coupled with ultrasound tomography.
\end{itemize}

%\Kiko{As mentioned above, the complexity of the HG solution scales proportionally to the number of sensors and depends on the position and the size of the region of interest. Below we discuss a number of scenarios in which a flexible ray based solution could offer a computational advantage warranting deployment of the Hamilton-Green solver.
%This implementation could be conducted either on its own or in conjunction with other techniques.} \Marta{Meaning of last sentence?}

\rB{Below we discuss a few scenarios which could benefit from the aforementioned advantages of the flexible ray based Hamilton-Green solver.}

\subsubsection{Joint Photoacoustic Tomography and Ultrasound Tomography Reconstruction}

The quality of the PAT reconstruction hinges upon the accurate knowledge of the speed of sound in the domain. 
In the absence of sound speed  measurements, there are two commonly taken approaches: i) assume the sound speed homogeneous and potentially fit its value via optimization; 
ii) deduce an inhomogeneous sound speed model from a priori information obtained from another non-acoustic modality (e.g.~MRI as in \cite{lou2017generation}). 
The first approach hits its limitations when the inhomogeneity becomes pronounced, which occurs even in soft tissue, 
with differences up to 10\% between different tissue types \cite{mast2000empirical}. 
Assuming sound speed difference as in \cite{mast2000empirical}, imaging to 7 [cm] depth in breast disregarding the inhomogeneity between the fat and muscle,  will result in an aberration on the order of mm and cause artefacts in the reconstruction. The second needs an extra scan which is potentially expensive, time consuming, requires registration of the images and the conversion of the other contrast to speed of sound is potentially difficult to justify. Therefore, recently there has been an interest in developing integrated \rB{ultrasound} (US) / PAT scanners, capable of simultaneous acquisition of the US and PAT data and hence providing a registered set of measurements from which both the speed of sound and PAT image can be reconstructed e.g.~\cite{matthews2017joint, matthews2017recon}.
As the forward solve for the complete data transmission US problem requires solution of $\mathcal O(n^2)$ full wave problems, the standard approach is to settle for the time of travel tomography i.e.~first time of arrival e.g.~\cite{li2009vivo, duric2005development}.
The Hamiltonian system underlying the bent ray tomography approaches such as \cite{li2009vivo, duric2005development} is the same as the one used in Hamilton-Green PAT solver thus the PAT solution can be obtained essentially for ``free'' when solving the transmission US problem. 
%\textcolor{green}{maybe the next sentence is irrelevant}
%While it is known that the speed of sound cannot be obtained from the PAT measurements, it is an open question what is the minimal complete set of US and PAT measurements for joint reconstruction under the usual smoothness and non-trapping assumptions.     

\subsubsection{Region of Interest}

The advances in the scanner hardware led to miniature endoscopic devices which can be deployed to guide surgical  intervention such as e.g.~the separation of the twin fetuses involving separation of veins in the placenta. In this case, the probe can be at a substantial distance from the target. This region is usually filled with a clear fluid. Furthermore, the  measurements have to be acquired almost instantaneously. 
The region of interest tomography, to the best of our knowledge, has not been widely considered in the context of PAT. 
The existing full wave solvers cannot take advantage of this scenario as the wave has to be propagated through the entire domain containing the region of interest and all the sensors.
On the other hand, the Hamilton-Green solver is inherently flexible and it is sufficient to restrict the propagation domain for each sensor to this sensor and the region of interest. Furthermore, as the clear fluid region can be assumed to have a homogeneous sound speed, and hence in the fluid the rays are straight lines, such a hybrid problem can be efficiently handled by simply increasing the time step in this region without deteriorating the quality of the solution. %\textcolor{green}{not sure if to say} On the other hand the full-wave solver would have to be explicitly coupled with a ray solver.

\subsubsection{Stochastic Methods and Subsampled Data}
The data acquisition time in PAT presents a bottleneck in pre- and clinical environment. Hence the recent interest in studying subsampling schemes and compressed sensing to reduce the acquisition time. The computational cost of  the full wave solvers is nearly independent of the number of sensors, while  the cost of Hamilton-Green solver scales linearly with the number of sensors. Furthermore, stochastic methods have been shown to outperform their deterministic counterparts in applications such as PET \cite{ehrhardt2017faster, chambolle2017stochastic}. Stochastic methods
use partial operators i.e.~a subset of the sensors in each iteration and so are particularly well suited for the Hamilton-Green solver.

\rB{\subsection{Contribution}}

%The main contribution of this paper is a new way to construct an approximation to the solution of the time-domain wave equation which arises in forward and adjoint problems in PAT. The bulk of the paper is devoted to establishing mathematical foundations of the proposed Hamilton-Green solver. The main building blocks of the solver are derived from methods used in seismic and underwater acoustics, including the computation of phase and amplitude along the rays.

%Based on this approximation we develop a novel acoustic solver for the forward and adjoint PAT problems for homogeneous and heterogeneous sound speeds. The proposed method is based on ray tracing which enables evaluation of partial forward/adjoint operators at a fractional cost which is not possible with full-wave solvers and is a common and desirable feature of  tomographic problems exploited by many algebraic reconstruction methods.

%In this paper, we propose a novel acoustic solver for the forward and adjoint problems in PAT. Our ray based solver effectively approximates the Green's function of the underlying acoustic problem along the trajectories of the Hamiltonian system resulting from the high frequency asymptotic \Kiko{approximate solution} to the wave equation. Hence, we term the solver Hamilton-Green (HG). 

\rB{
The main theoretical contribution of this paper is a new way to construct an approximation to the solution of the time-domain wave equation which arises in forward and adjoint problems in PAT. We term our solver Hamilton-Green (HG), as it effectively approximates the time derivative of the unknown Green's function of the underlying acoustic problem along the trajectories of the Hamiltonian system resulting from the high frequency asymptotic approximation to the solution of the wave equation. 
The bulk of the paper is devoted to establishing the mathematical foundations of the proposed Hamilton-Green solver.
The main building blocks of the solver are derived from methods used in seismic and underwater acoustics, including the computation of phase and amplitude along the rays.}

\rB{
Based on this approximation we develop a novel acoustic solver for the forward and adjoint PAT problems for homogeneous and heterogeneous sound speeds. 
% A proof of concept Matlab implementation for 2D problems is available from a Github repository \texttt{green-ray}\footnote{\texttt{green-ray} on Github: \href{https://github.com/kikorulan/green-ray}{https://github.com/kikorulan/green-ray}}. 
The proposed method is inherently ray based which enables evaluation of partial forward/adjoint operators (for each sensor or even each ray) at a fractional cost which is not possible with full-wave solvers and is a common and desirable feature of tomographic problems exploited by many algebraic reconstruction methods.} 

%Finally, the solver is evaluated on a 2D domain with homogeneous sound speed to numerically prove the numerical convergence of our solution.

\rB{
We discuss the approximation error of the Hamilton-Green solver for the forward and adjoint problems for PAT and support our conclusions with 2D examples in domains with homogeneous and heterogeneous sounds speeds. We compare our results to full-wave solution obtained with a pseudo-spectral finite difference method implemented in \texttt{k-Wave}\footnote{\href{http://www.k-wave.org/}{http://www.k-wave.org/}}, a well established and widely used toolbox for biomedical PAT and ultrasonic simulations.}
%
%vessel phantom  with complete data and compare our results to the full-wave solution using a pseudo-spectral method.
%
%We include forward and adjoint problems with complete data and compare our results to the full-wave solution using a pseudo-spectral method. 

\rB{
We consider full-data problems which, as alluded to in section \ref{sec:intro:motiv}, while native to the full-wave solvers are not a realistic deployment scenario for our HG solver. However, we stress that our purpose here is not to compete with the full-wave solver in this context, nor to produce PAT image reconstruction, but to evaluate the impact of the approximation underlying the HG solver on both the forward and the adjoint problems. This is also the reason for restricting ourselves to 2D domains which allow for a better visualization.} 

In future work, the HG solver will be deployed in a number of realistic 3D scenarios which can benefit from the additional flexibility it offers.

\subsection{Outline}
The remainder of this paper is organized as follows. In section \ref{sec:rayequations} we recall the high frequency \rB{approximation to the solution of the wave equation} including the equations for the phase, the amplitude and two methods for computation of the Jacobian of the change of coordinates between the Cartesian and ray based coordinates (which corresponds to the inverse ray density). Furthermore, we derive expressions for the reversed quantities needed by the adjoint Hamilton-Green solver. 
%An alternative approach, the discrete interface method, based on the eikonal equation and Snell's law is proposed in section \ref{sec:dim}.   
Section \ref{sec:solver} constitutes the core of the paper where we introduce the novel Hamilton-Green acoustic solver for the forward and adjoint problems in PAT. 
%In \cref{sec:rayequations} we derive the ray equations and discuss their limitations in terms of accuracy of the solution to the wave equation. %some of the mathematical challenges that arise from their use.
%In \cref{sec:solver} we introduce a novel Hamilton-Green acoustic solver for the forward and adjoint problems in PAT. 
%addressing both numerical and implementation aspects of the adaptation of ray tracing to the PAT context.
\rB{In Section \ref{sec:errana} we analyse the error of the forward and adjoint Hamilton-Green solver. 
In section \ref{sec:simulations} we report simulations for both the forward and the adjoint problems. We illustrate HG convergence in homogeneous case, and evaluate HG error in inhomogeneous case against solutions obtained with the first order finite difference pseudo-spectral method implemented in \texttt{k-Wave} toolbox.} 
Section \ref{sec:conclusions} summarizes the contributions and conclusions of the paper and provides outlook on future research.
Technical details of the solver pertaining mapping between Cartesian and ray based grids and Green's function derivative approximation are postponed to  \ref{app}. \rA{A case study of the effect of caustics on HG solver on an acoustic lens example is provided in \ref{app:lens}}.

%==============================
% 2. Ray Equations
%==============================
%================================================================================
%====================   HIGH FREQ APPROX FOR THE WAVE EQUATION ==================
%================================================================================
\section{High frequency \rB{approximate solution} to the wave equation}
\label{sec:rayequations}
We consider the linear scalar wave equation
\begin{equation}\label{wave}
u_{tt} - c^2(x)\Delta u = 0, \quad(t, x) \in (0, \infty)\times\mathbb{R}^d,
\end{equation}
with $d = 2, 3$ and appropriate initial / boundary conditions to be specified later,
where $c(x)$ is the speed of the wave in the medium. 
%The initial data and boundary conditions define the frequency spectrum of the solution.
In some scenarios, the initial/boundary/source conditions induce relatively high (\rB{with respect to the size of the domain}) \textit{essential frequencies} \cite{engquist2003computational} in the propagating waves. 
Accurate resolution of such high frequencies would require correspondingly fine discretization, rendering the problem computationally infeasible.
% 
%Since the precision of the direct numerical simulation depends heavily on the number of points taken per wavelength, 
%solving the wave equation as in (\ref{wave}) becomes unfeasible in terms of computational resources.
%In a direct numerical simulation, in order to capture those essential frequencies the spatial discretization would have to be very fine, 
%which would make this simulation computationally infeasible.

In the high frequency limit $\omega\rightarrow \infty$, one can consider the following approximation to (\ref{wave}).
%The \textit{geometrical optics} equations \review{equations for geometrical optics}
%provide a proper mathematical structure for this purpose \cite{born2013principles}.
Let $u$ admit a series expansion of the form
\begin{equation}\label{highFreq}
u(t, x) = e^{i\omega\phi(t, x)} \sum_{k = 0}^\infty A_k(t, x)(i\omega)^{-k},
\end{equation}
where $\phi(t, x)$ is the \textit{phase}, $A_k$ the \textit{coefficients of the amplitude} and $\omega$ the \textit{frequency} of the oscillating wave. 
%The approximation of the wave solution as the series given in (\ref{highFreq}), also known as WKB expansion, holds valid as $\omega\rightarrow\infty$.
In this representation we expect the phase $\phi$ and the amplitude coefficients $A_k$ to vary at a much lower \rB{temporal and spatial rate than the wave field $u$ since they are independent of $\omega$}. 
%Therefore, we can use a coarser discretization to compute them. 

The \textit{geometrical optics equations} are obtained by substituting the WKB expansion \eqref{highFreq} into the wave equation \eqref{wave} and equating terms of the same order to ensure that \eqref{wave} holds down to $\mathcal{O}(\omega)$:
terms $\mathcal{O}(\omega^2)$ result in the \textit{eikonal equation},
%\begin{equation}\label{req:eikonal}
%\phi_t \pm c(x) \|\nabla\phi\| = 0,
%\end{equation}
while terms $\mathcal{O}(\omega^1)$ result in the \textit{transport equation}.
%\begin{equation}\label{req:transportTime}
%(A_0)_t + c(x)\frac{\nabla\phi\cdot \nabla A_0}{\|\nabla \phi\|} + \frac{c^2(x)\Delta \phi -\phi_{tt}}{2c(x) \|\nabla \phi\|} = 0.
%\end{equation}
%In the limit $\omega\rightarrow\infty$ the terms of order $\mathcal{O}(\omega)^{-n}$ for $n\geq 0$ can be neglected.

\rA{The literature on solution of the geometrical optics equations is rich, thus the following list is necessarily incomplete.} \rB{We particularly mention Beylkin's seminal work on Fourier integral operators \cite{beylkin1985imaging},}
\rA{phase-space methods for solution of eikonal equation with multi-valued phases \cite{benamou2000eulerian, cheng2002level, fomel2002fast, engquist2002high} 
%Hamilton-Jacobi methods \cite{cao1994finite, gremaud2006computational, qian2006local}, 
and methods that solve for both the multi-valued phases and amplitudes \cite{qian2002adaptive, qian2006local, qian2011efficient} and finally ray tra\-cing methods which solve the Hamiltonian system and the transport equation \cite{vcerveny1979ray, virovlyansky2003ray}.} 
Our solver is based on \textit{ray tracing} \rB{approach, which has immediate parallels to other tomographic modalities, most prominently X-ray computed tomography, and} has been studied extensively
both theoretically and numerically, in particular in the context of seismic imaging \cite{richards1980quantitative, cerveny2005seismic, um1987fast}
and ocean acoustics \cite{jensen2000computational, porter1987gaussian}.

\rA{The immediate role of the ray equations in our approach is as means to \textit{reparametrisation} of the domain $\Omega \subset \R^d$. Therefore we
assume a fully implicit time dependence of the amplitude and phase, i.e.~solely through the time dependence of the trajectory $x(t)$
\begin{equation*}
u(t, x(t)) = A(x(t))\exp^{i\omega \phi(x(t))}
\end{equation*}
and solve the frequency domain version of the eikonal and transport equations.
}

%================================================================================
%====================             PHASE                        ==================
%================================================================================
\subsection{Phase}

The frequency domain version of the eikonal equation reads
\begin{equation}\label{eq:eikonalFreq}
\|\nabla \phi\| = 1/c = \eta,
\end{equation}
where $\eta(x) = 1/c(x)$ is the \textit{slowness} of the medium.
Following the construction in \cite{runborg2007mathematical}, we introduce the Hamiltonian \hbox{$H(x, p) = c(x)\|p\|$} defined in the phase space $\mathbb{R}^d\times\mathbb{R}^d$ %, with $d$ the dimensionality of the domain.
of double the dimensionality of the ambient space, $d$.
Let $(x(t), p(t))$ be a bicharacteristic pair associated with this Hamiltonian. 
The Hamiltonian $H$ is constant along these bicharacteristics and is set to the initial value $H(x_0, p_0) = 1$,
which corresponds to $\|p\| = \eta$.
Therefore we have
\begin{subequations}\label{req:hamiltonianEta}
\begin{align}
    &\dfrac{\D x}{\D t} = \nabla_p H(x, p) = \dfrac{p}{\eta^2}, && x(0) = x_0, \label{req:1stHamil}\\[15pt]
    &\dfrac{\D p}{\D t} = -\nabla_x H(x, p) = \dfrac{\nabla\eta}{\eta}, && p(0) = p_0, \qquad \|p_0\| = \eta(x_0), \label{req:2ndHamil}
\end{align}
\end{subequations}
where $(x(0), p(0)) = (x_0, p_0)$ are the initial conditions for the Hamiltonian system \eqref{req:hamiltonianEta}.
\rB{The point $x_0 \in\mathbb{R}^d$ is the \textit{origin} of the ray and should be chosen on a non-characteristic hypersurface $S \subset \mathbb{R}^d$}.

\paragraph{\rB{Method of characteristics}} 
Consider $\phi(x)$ smooth such that
\begin{equation} H(x, \nabla\phi(x)) = 1. \end{equation}
Due to the uniqueness of the solution of \eqref{req:hamiltonianEta}, \rB{applying the method of characteristics we find that} $(x(t), \nabla\phi(x(t)))$ 
is a bicharacteristic pair of $H$ with $p(t) = \nabla\phi(x(t))$.
This pair can be interpreted as follows: $x(t)$ represents the \textit{trajectory} in the domain, while the slowness vector $p(t)$ is the direction of propagation at each point along that trajectory.
%Additionally, because of the ansatz \eqref{eq:EquivFreqT}, diffe\-ren\-tiating $\phi$ with respect to $t$ yields constant, implying a linear relation between the time and phase 
Furthermore, for $\phi(t,x) = \phi(x(t))$ the eikonal equation \eqref{eq:eikonalFreq} implies a linear relation between the time and phase
\begin{equation}\label{eq:linphase} \phi(x(t)) = \phi(x_0) + t. \end{equation}
%The function $\phi(x(t))$ solves the eikonal equation along the ray. 
The independence of bicharacteristic pairs (corresponding to different initial conditions for the bicharacteristic system \eqref{req:hamiltonianEta}) naturally accommodates multiple phase solutions, in contrast to the inherently single phased viscosity solution (first time of arrival). 
The system of ODEs (\ref{req:hamiltonianEta}) is numerically solved using a 2nd order Runge-Kutta method.
The details of the interpolation between the Cartesian and ray based grids are given in \ref{app}.

%================================================================================
%====================             AMPLITUDE                    ==================
%================================================================================

\subsection{Amplitude}
To compute the amplitude we solve the frequency domain version of the transport equation %\eqref{req:transportTime}
\begin{equation}\label{req:transport}
2\nabla\phi\cdot\nabla A + \Delta \phi A = 0.
\end{equation}
%where we have dropped the subindex $0$ for the amplitude term.
The solution of the first order equation \eqref{req:transport} for the bicharacteristic pair $(x(t; x_0), p(t; p_0))$ at a point $x(t; x_0)$ can be explicitly written as
\begin{equation} \label{req:amplitude}
A(x(t; x_0)) = A(x_0)\frac{\eta(x_0)}{\eta(x(t; x_0))}\sqrt{\frac{q(0; x_0)}{q(t; x_0)}},
\end{equation}
where $q$ is the determinant of the Jacobian $J$ of $x$ with respect to the initial data,
\rB{
\begin{equation}\label{req:jacobian}
q(t; x_0) = \det J := \det D_{x_0} x(t; x_0) = \det \left(\cfrac{\partial x(t; x_0)}{\partial x_{0, i}}\right).
\end{equation}
}
%================================================================================
%====================        DETERMINANT OF THE JACOBIAN       ==================
%================================================================================

\subsection{Jacobian determinant}
\rA{For completeness we describe two methods for numerical evaluation of the determinant of the Jacobian, $q$, along the ray trajectory \cite{engquist2003computational,cerveny2005seismic}}. We restrict our presentation to 2D as 3D follows analogously.
%In the following we present two different methods for numerical evaluation of the determinant of the Jacobian, $q$, along the ray trajectory.We restrict the presentation to 2D as 3D follows analogously.
%\Kiko{An alternative computation of the Jacobian $q$ can be found in \cite{lu2016babich}.}

\begin{figure}
\begin{center}
 \includegraphics[width=0.5\textwidth]{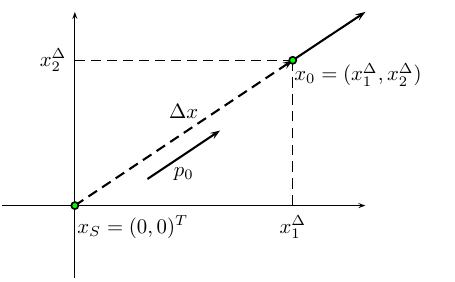}
\caption{The relation of the origin (shooting point) of the ray $x_S$ to its initial point $x_0$ as used to calculate the initial conditions for the ODE system \eqref{eq:odeInit}.}
\label{ODEsource}
\end{center}
\end{figure}

%========================================
% ODE Method
%========================================
\paragraph{ODE method}
The time evolution of the determinant $q$ can be computed solving the ODE system obtained by differentiating the Hamiltonian with respect to the initial conditions \cite{engquist2003computational} and changing the order of the differentiation
\rB{
\begin{equation}\label{eq:ODEsystem}
\frac{\D}{\D t}
\begin{pmatrix}D_{x_0} x\\ D_{x_0} p\end{pmatrix} =
\begin{pmatrix}D^2_{px} H & D^2_{pp} H\\
    -D^2_{xx} H&  -(D^2_{px} H)^T\end{pmatrix}
\begin{pmatrix}D_{x_0} x\\ D_{x_0} p\end{pmatrix}
\end{equation}
}
with the initial conditions
\rB{
\begin{equation}
\label{eq:odeInit}
D_{x_0} x(0; x_0) = I, \quad D_{x_0} p(0; x_0) = D^2\phi(x_0).
\end{equation}}
For this $t$ parametrization of the Hamiltonian, $H = c(x)\|p\|$,
the entries of the system matrix (\ref{eq:ODEsystem}) are
\rB{\begin{equation*}
 D^2_{px} H = c\left(p\nabla^T c\right), \qquad 
 D^2_{xx} H = cD^2c, \qquad
 D^2_{pp} H = c^2 \left[ I - c^2 (pp^T)\right].
\end{equation*}}
It remains to specify the second initial condition, $D^2\phi(x_0)$.
Assuming $c$ to be constant in a neighborhood of $x_S$, the \textit{origin of the ray} or \textit{shooting point}, %resulting
results in an out-propagating spherical wave (phase)\rB{. Setting $\phi(x_S) = 0$ we have}
\begin{equation*} \phi(x) = \frac{\|x-x_S\|}{c} . \end{equation*}
Furthermore, without loss of generality %we consider 
choosing the wave centered at the origin
$x_S = (0, 0)^T$ we obtain
\begin{equation}\label{eq:phiSph}
\phi(x) = \frac{\|x\|}c = \frac{\sqrt{x_1^2 + x_2^2}}c,\end{equation}
with $x = (x_1, x_2)^T$ the Cartesian coordinates of $x$. 
The Hessian of \eqref{eq:phiSph} 
\rB{
\begin{equation*}
D^2(\phi(x)) = \frac{1}{c(x_1^2 + x_2^2)^{3/2}}\begin{pmatrix} x_2^2 & -x_1x_2 \\ -x_1x_2 & x_1^2 \end{pmatrix}
\end{equation*}}
has a singularity at the shooting point $x_S$. 
To circumvent this issue we assume that the initial conditions are evaluated at a nearby point $x_0$  further down the ray. To be precise, the \textit{initial point} $x_0$ at a distance $\Delta x$ in the shooting direction $p_0$ (see figure \ref{ODEsource})
% Figure
\begin{equation}
 x_0 = \Delta x\dfrac{p_0}{\|p_0\|} = \begin{pmatrix} x_{1}^\Delta \\ x_{2}^\Delta \end{pmatrix}. 
\end{equation}
The Hessian at $x_0$
\rB{
\begin{equation*}
D^2\phi(x_0) = \frac{1}{c((x^\Delta_1)^2 + (x^\Delta_2)^2)^{3/2}}\begin{pmatrix} (x^\Delta_2)^2 & -x^\Delta_1x^\Delta_2 \\ -x^\Delta_1x^\Delta_2 & (x^\Delta_1)^2 \end{pmatrix}.
\end{equation*}
}
is not singular but it depends on the initial point $x_0$ and hence implicitly on $\Delta x$, the distance from $x_0$ to the origin of the ray, $x_S$.

%========================================
% Proximal Ray Method
%========================================
\paragraph{Proximal ray method}
%The \Marta{here presented} numerical approximation of the Jacobian $J$ of $x$ with respect to the initial conditions $x_0$, which we refer to as \textit{proximal ray method}, 
Inspired by an analogous construction in \cite{cerveny2005seismic}, we use a nearby ray to construct a finite difference approximation to the Jacobian $J$.

Let $\{\theta, \tau\}$ be a local coordinate system anchored at the initial point $x(0) = x_0$ of the ray $x(\tau)$, with $\theta$ the shooting angle and $\tau$ the parameter along the ray. We can write the Jacobian with respect to the initial conditions $x_0$ as %is given by the matrix
\begin{equation*}
J_\tau(x) = \begin{pmatrix} \dfrac{\partial x}{\partial \theta} & \dfrac{\partial x}{\partial \tau}\end{pmatrix}
\end{equation*}
and approximate the derivatives with forward finite differences
%\begin{subequations}
\begin{align*}
 &\dfrac{\partial x(\tau, p_0; x_0)}{\partial \theta} \approx \dfrac{x(\tau, p_0^\theta; x_0) - x(\tau, p_0; x_0)}{\Delta\theta}, \\[15pt]
 &\dfrac{\partial x(\tau, p_0; x_0)}{\partial \tau} \approx \dfrac{x(\tau+\Delta\tau, p_0; x_0) - x(\tau, p_0; x_0)}{\Delta\tau}. 
\end{align*}
%\end{subequations}

%%  \begin{figure}
%%  \begin{center}
%%  \hspace*{1cm}\includegraphics[width=0.7\linewidth]{Figures/2_JacobianDet.pdf}
%%    \vspace{-0.5cm}
%%  \caption{The interpretation of the determinant of the Jacobian via proximal ray method. 
%%  %\todo{update notation to ";" to correspond to text.}
%%  }
%%  \end{center}
%%  \label{jacobian}
%%  \end{figure}

Here $x(\tau, p_0; x_0)$ denotes the \textit{main ray},
$x(\tau, p_0^\theta; x_0)$ an \textit{auxiliary ray} shot at the angle $p_0^\theta = p_0 + \Delta\theta$ and $\Delta\tau$ is the step length along the ray (chosen equal to the step length used while solving for the trajectories). 
With this notation, we have the determinant of the Jacobian $q$ \cite{cerveny2005seismic}
\begin{equation}\label{qSing}
q_\tau(\tau; x_0) = 
\begin{vmatrix} \dfrac{x(\tau, p_0^\theta; x_0) - x(\tau, p_0; x_0)}{\Delta\theta} & \dfrac{x(\tau+\Delta\tau, p_0; x_0) - x(\tau, p_0; x_0)}{\Delta\tau} \end{vmatrix}.
\end{equation}

\begin{figure}
\begin{center}
 \includegraphics[width=0.8\linewidth]{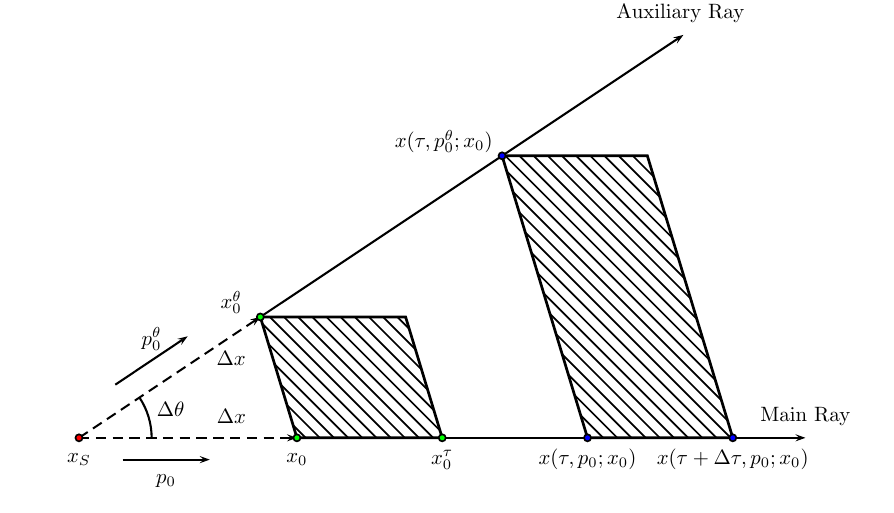}
 \vspace{-0.5cm}
 \caption{The relation of the origin (shooting point) of the ray $x_S$ to its initial point $x_0$ in the proximal ray method for the determinant of the Jacobian. 
 %\todo{";" same as in previous figure}
 }
\label{jacobianProx} 
\end{center}
\end{figure}

The determinant admits an interpretation as the area of the parallelogram spanned by the vectors $\left\{\dfrac{\partial x}{\partial \theta}, \dfrac{\partial x}{\partial \tau}\right\}$ (dashed area in figure \ref{jacobianProx}). 
%divided by the differential elements of the local coordinates $\Delta \theta$ and $\Delta \tau$.

The amplitude \eqref{req:amplitude} requires the evaluation of the Jacobian \eqref{qSing} at the ray initial point $x_0$, at which the determinant is singular. We circumvent the problem analogously as in the ODE method by shifting the origin of the ray to a nearby point $x_S = x_0 - \Delta x \cdot p_0 /\|p_0\|$ (see figure \ref{jacobianProx}), which yields $x_S = x(-\Delta\tau, p_0; x_0), \,x_0^\tau = x(\Delta\tau, p_0; x_0) \,$ and 
\rB{$\, x_0^\theta = x(\Delta\tau, p_0^\theta; x_S)$}
%$\, x_0^\theta = x(0, p_0^\theta; x_0)$ \todo{the last is $x_0$ which is not equal $x_0^\theta$ how about $\, x_0^\theta = x(\Delta\tau, p_0^\theta; x_S)$?} 
and the non-singular approximation to the determinant
%With our numerical approximation as it is, this determinant has a singularity at $x_0$. To remedy this, we assume that the rays are shot from a proximal point $x_S$, as in the ODE method. The determinant is evaluated along the ray (including their shooting point) as follows, shown in figure \ref{jacobianProx}. The ray shooting point is $x_S = x_0 - \Delta x \cdot p_0 /\|p_0\|$.
%Now we need to obtain the Jacobian at $x_0$, for which we should consider the points $x_0 = x(0, p_0, x_S), \,x_0^\tau = x(\Delta\tau, p_0, x_S) \,$ and $\, x_0^\theta = x(0, p_0^\theta, x_S).$
%Substituting these in formula (\ref{qSing}), the determinant of the Jacobian is given by
\begin{equation*}
q_\tau(0; x_S) \approx \begin{vmatrix} \cfrac{x_0^\theta - x_0}{\Delta\theta} & \cfrac{x_0^\tau - x_0}{\Delta\tau}\end{vmatrix}.
\end{equation*}
Finally, the point $x_0^\theta$ corresponds to a rotation of $x_0$ by an angle $\Delta\theta$ around the origin of the ray $x_S$,
\begin{equation*}
 x_0^\theta = x_S + \Delta x \dfrac{p_0^\theta}{\|p_0^\theta\|}, \quad
 p_0^\theta = \begin{pmatrix} \cos \Delta\theta & -\sin\Delta\theta \\ \sin\Delta\theta & \cos\Delta\theta\end{pmatrix} p_0.
\end{equation*}
We observe that $q$ (and hence the amplitude $A$) depends on %the chosen proximal point 
the choice of the ray origin $x_S$, which is determined by the increment $\Delta x$.
In a homogeneous 2D medium, the amplitude decays as $A(x) \sim \sqrt{\Delta x}/\sqrt{\Delta x + \|x\|}$. Assuming $\Delta x \ll \|x\|$, the amplitude behaves as $A(x) \approx \sqrt{\Delta x}/\sqrt{\|x\|}$ indicating that the shift by $\Delta x$ essentially corresponds to a scaling factor. 
%Furthermore, the amplitude is independent of the angle increment $\Delta\theta$.

%Assuming that the wave propagates in a homogeneous 2D medium, the amplitude has a decay of the form $A(x) \sim \sqrt{\Delta x}/\sqrt{\Delta x + \|x\|}$. 
%Then, for any $\|x\| \ll \Delta x$, the amplitude is approximately $A(x) \approx \sqrt{\Delta x}/\sqrt{\|x\|}$, that shows that the term $\Delta x$ does not affect the shape of the amplitude and can be simplified as a known multiplicative factor, given that $\Delta x$ is chosen to be sufficiently small compared to the spatial resolution in the domain $\Omega$. Additionally, the amplitude is independent of the incremental angle $\Delta\theta$.

%================================================================================
%====================       REVERSE AMPLITUDE                  ==================
%================================================================================

\subsection{Reversing rays}\label{sec:RevRay}
The computation of the ray trajectories is at the core of our Hamilton-Green (HG) solver. 
For the solver to handle both the forward and the adjoint problems, it has to be capable of propagating the pressure wave both from the domain towards the sensors and from the sensors back into the domain. A naive approach would entail treating these two problems separately i.e.~computing the ray trajectories for each independently. However, here the natural choice is to shoot the rays from the sensors into the domain. \rB{In the following we show that such trajectories are reversible. Hence, they can be used to propagate the wave in both directions. We derive an expression for the reversed phase $\phi_R$ and the reversed amplitude $A_R$ over the reversed ray $x_R(t)$ in terms of the original phase $\phi$ and amplitude $A$ over the original ray $x(t)$.}
%As we explain in the next section, for the forward problem we need to propagate the pressure from the domain to the sensors. A possible way to do this would be to shoot a ray from each pixel towards the corresponding sensor. However, this solution is computationally infeasible, as this naive approach would require shooting an extremely large number of rays.  In contrast, we propose to shoot rays from the sensor into the domain, which allows us to efficiently  cover the entire domain.  The  pressure amplitudes need to be computed along the reversed rays (i.e. from the domain to the sensors). We remark that the such computed ray trajectories can be used for both the forward and adjoint problem solution.  
%
%To verify that this approach is correct we prove that the trajectories that solve the ray equations are reversible.
%This means that given a ray $x(t)$ shot from $x(0) = x_A,\, p(0) = p_A$ to $x(T) = x_B, p(T) = p_B$ we can obtain  its reverse ray 
\rB{\subsubsection{Reversed trajectory and phase}} Given a ray $x(t)$ shot from $x(0) = x_A,\, p(0) = p_A$ to $x(T) = x_B, p(T) = p_B$, we obtain its reversed ray as
\begin{equation*}
x_R(t) = x(T-t),\qquad
 p_R(t) = -p(T-t).
\end{equation*}
%It is easy to see that the reversed ray also satisfies the Hamiltonian system
%\begin{subequations}\label{eq:revRayHamil}
%\begin{align}
%&&&&\dfrac{\D x_R(t)}{\D t} &= \dfrac{p_R(t)}{\eta^2}   & \Longleftrightarrow &&  \dfrac{\D x(T-t)}{\D t} &= -\dfrac{p(T-t)}{\eta^2}&&&& \\
%&&&&\dfrac{\D p_R(t)}{\D t} &= \dfrac{\nabla\eta}{\eta} & \Longleftrightarrow && -\dfrac{\D p(T-t)}{\D t} &= \dfrac{\nabla\eta}{\eta}&&&&
%\end{align}
%\end{subequations}
%with the initial conditions $x_R(0) = x_B, p_R(0) = -p_B$.
%Undoing the change of variables on the right in \eqref{eq:revRayHamil}, $t' = T-t, \D t' = -\D t$, we recover the original Hamiltonian for $x(t)$ with the initial conditions $x(0) = x_A, p(0) = p_A$
%These equivalences are proven by undoing the change of variables $t' = T-t$ and applying the chain rule. For the initial and end points we have
%\begin{subequations}
\rB{Following the principle of reciprocity, we observe that the reversed ray also satisfies the Hamiltonian system with the initial conditions $x_R(0) = x_B, p_R(0) = -p_B$.
Summarising, for $x$ and $x_R$ it holds}
\begin{align*}
&x_R(0) = x(T) = x_B, && x_R(T) = x(0) = x_A,\\
&p_R(0) = -p(T) = -p_B, && p_R(T) = -p(0) = -p_A.
\end{align*}
%\end{subequations}
%so the rays are indeed reversible. The reversibility of the rays opens the question of how to obtain the phase $\phi_R(t)$ and amplitude $A_R(t)$ for the reverse ray $x_R(t)$ as in (\ref{reverseRay}) given $\phi(t)$ and $A(t)$ associated to $x(t)$. 
%The remaining question is how to obtain the phase $\phi_R(t)$ and the amplitude $A_R(t)$ along the reversed ray  $x_R(t)$ in terms of the phase $\phi(t)$ and amplitude $A(t)$ of $x(t)$.

Recalling that the phase is essentially the propagation time (c.f.~change of variables $t' = T-t$), we immediately obtain
%\todo{I think this only holds true for initial phase 0, otherwise one would need + $\phi(0)$ (which we are using to avoid singularity, there also is inconsistency with time $t_0$ not being 0). Also actually we used the notation $\phi(x(t))$ before, did we say at any point we are going to simplify like that?}
\begin{equation}\label{eq:phaseRev}
%\phi_R(t) = \phi(x(T)) - \phi(x(T-t)) + \phi(x_0).
\phi_R(t) = \phi(x(T-t)) =  T-t + \phi(x_0).
\end{equation}
%\todo{changed this formula, check if you agree}
%This holds because the propagation time is also reversible. The formula for the amplitude $A_R(t)$ requires a little more effort. To derive it we use the arguments presented in the discrete interface approach method. Let us recall the total attenuation due to spherical spreading and interface losses along the ray, as defined in (\ref{req:totalAt}), 
%To derive the formula for the amplitude $A_R(t)$, we recall
%the total attenuation due to spherical spreading and interface losses along the ray, \eqref{req:totalAt}
%\begin{equation}
%\alpha_T(t) = \alpha_S(t)\cdot\alpha_I(t) \equiv \frac{A_0}{A(t)}.
%\end{equation}
%\rev{I think you want this the other way around? The numbers should be smaller than 1? This has impact above e.g. (42), (43) and below see comments.}

\subsubsection{Reversed amplitude}
In the same spirit we would like to express the reversed amplitude $A_R(x_R(t))$ along the reversed ray $x_R(t)$ in terms of the original trajectory $x$, phase $\phi$ and the  inverse ray density (Jacobian determinant) $q$.
From the first principle the reversed inverse ray density, $q_R$, is the determinant of the Jacobian of $x_R$ w.r.t.~the initial point of $x_R$ (which is the end point of $x$) i.e.~$x_T := x(T;x_0) = x_R(0;x_T)$ 
%with $A_0(x_0)$ the amplitude at $x_0$, $\xi(x_0) = \sqrt{q(0, x_0) \eta(x_0)}$ a constant that depends on $x_0$ and $\xi(x(t, x_0)) = (\sqrt{q(x(t, x_0))\eta(x(t, x_0))})^{-1}$ the term that varies along the trajectory.
%The determinant for the reversed ray is 
\begin{equation}\label{eq:qR}
q_R(t; x_T) = \left|D_{x_T} x_R(t; x_T)\right|, 
\end{equation}
with 
\begin{equation}\label{req:reverse} 
x_R(t; x_T) := x(T-t; x_0)
\end{equation} 
Substituting \eqref{req:reverse} into \eqref{eq:qR}, applying the chain rule and the determinant calculus yields
\begin{equation*}
q_R(t; x_T) = |D_{x_T} x(T-t; x_0)| = |D_{x_0} x(T-t; x_0)| \left|\frac{\partial x_0}{\partial x_T}\right|,
\end{equation*}
where $|D_{x_0} x(T-t; x_0)| = q(T-t; x_0)$ by definition and 
$$|\partial x_0/\partial x_T| = | \left. \partial x(t; x_0)/\partial x_0 \right|_{t=T} |^{-1} = | \left. D_{x_0} x(t;x_0) \right|_{t=T} |^{-1} = 1/q(T; x_0).$$ 
Thus we expressed $q_R$ in terms of $q$
\begin{equation}\label{req:reverseQ}
q_R(t; x_T) = \frac{q(T-t; x_0)}{q(T; x_0)}.
\end{equation}

The reversed amplitude $A_R(x_R(t; x_T))$ along the reversed ray $x_R(t)$ in terms of reversed quantities follows immediately from \eqref{req:amplitude}
\begin{equation}\label{req:revAmplitude}
A_R(x_R(t; x_T)) = A_0(x_T) \sqrt{\frac{q_R(0; x_T)}{q_R(t; x_T)}}\frac{\eta(x_T)}{\eta(x_R(t; x_T))}.
\end{equation}
Substituting the reversed quantities with \eqref{req:reverse} and \eqref{req:reverseQ} we obtain an expression in terms of the original quantities $x$ and $q$
\begin{equation*}
A_R(x_R(t; x_T)) = A_0(x_T) \sqrt{\frac{q(T; x_0)}{q(T-t; x_0)}}\frac{\eta(x_T)}{\eta(x(T-t; x_0))}.
\end{equation*}

\dontshow{
It is illuminating to decompose the amplitude \eqref{req:amplitude} into quantities dependent on the initial values only: $A(x_0)$ - the amplitude at $x_0$ and $\xi_0(x_0)$, and varying along the trajectory: $\xi(x(t, x_0))$
%\begin{equation}\label{req:amplitudeXi}
%A(x(t, x_0)) = A_0(x_0)\,\xi_1(x_0)\, \xi(x(t, x_0)),
%\end{equation}
\begin{equation}\label{req:amplitudeXi}
A(x(t, x_0)) = A_0(x_0) \, \underbrace{  q(0, x_0)^{\frac12}  \eta(x_0) }_{=:\xi_0(x_0)} \,  \underbrace{  q(x(t, x_0))^{-\frac12}  \eta(x(t, x_0))^{-1}}_{=:\xi(x(t, x_0))}.
\end{equation}
%We seek an analogous decomposition for the reversed amplitude $A_R(x_R(t'))$.
Rearranging terms and following the notation in \eqref{req:amplitudeXi} we obtain an analogous decomposition for the reversed amplitude 
\begin{equation}
A_R(x_R(t', x_T)) = A_0(x_T)\, \underbrace{q(T, x_0)^{\frac12}\eta(x_T)}_{=:\xi_{0R} (x_T)}\, \underbrace{q(T-t', x_0)^{-\frac12}\eta(x(T-t', \Marta{x_0}))^{-1}}_{=:\xi_R(x_R(t', x_T)) = \xi_R(x(T-t', x_0))}.
\end{equation}
%with $\xi_{0R}(x_T) = q(T, x_0)^{\frac12}\eta(x_T)$ and $\xi_R(x_R(t', x_T)) = q(T-t', x_0)^{-\frac12}\eta(x(T-t', x_T))^{-1}$.
\todo{The factors a naturally defined for the amplitude expressed in terms of the reversed quantities, not after the substitution. I am wondering if we should get rid of them }
}
%================================================================================
%====================        CAUSTICS AND SHADOW REGIONS       ==================
%================================================================================

\subsection{Caustics and shadow regions}
\label{re:sec:caustics}
We finish this section with the discussion of the two limitations of the ray tracing approach, the \textit{caustics} and the \textit{shadow regions}.
%In the last part of this section we introduce the two most relevant mathematical challenges we have  encountered so far in the development of this PAT solver, namely \textit{caustics} and \textit{shadow regions}. 

%========================================
% Caustics
%========================================
\subsubsection*{Caustics}
\rB{
A \textit{caustic} is a set of points at which the determinant of the Jacobian of the coordinate change, \eqref{req:jacobian}, becomes singular.}
%with respect to the initial conditions using the Cartesian coordinates \eqref{req:jacobian}, becomes singular.}

%\rev{is this the common definition? Note, that this is only for the Jacobian wrt Cartesian coordinates, so I am not really liking it ... If it is not the formal definition, either use the formal definition and make the connection below, or do not formalize this as definition, just formulate a normal sentence.}
As the determinant of the Jacobian \eqref{req:jacobian} is in the denominator of the amplitude formula \eqref{req:amplitude}, the latter becomes unbounded at caustics. Thus \eqref{req:amplitude} cannot be used to compute the amplitude directly at and beyond the caustic. The following modification
%Our interest in these sets is due to the fact that the formula that we gave for the amplitude along the rays becomes unbounded on them. This means that the amplitude is not defined on caustics for ray tracing. As we show in the next section, in order to compute the time signals measured at the sensors, we need the amplitude attenuations on the rays, which cannot be fully achieved in the presence of caustics. The amplitude formula (\ref{req:amplitude}) given for the ODE and Proximal Ray methods breaks down at caustics: when the determinant $q$ of the Jacobian with respect to the initial conditions is zero, the amplitude is not defined. Therefore, additional mathematical tools are needed.
has been proposed beyond the caustic (but not at the caustic itself)
%An improvement over the formula for the amplitude given in \eqref{req:amplitude} consists of introducing 
\begin{equation}\label{eq:AMaslov}
A(x(t, x_0)) = A(x_0)\frac{\eta(x_0)}{\eta(x(t; x_0))} \sqrt{\left|\frac{q(0; x_0)}{q(t; x_0)}\right|} e^{-i m(t) \frac{\pi}{2}},
\end{equation}
where $m(t)$ is the \textit{Keller-Maslov} index which counts the number of times the ray crosses a caustic (see e.g.~\cite{runborg2007mathematical}). 
\rA{To control the numerical integration error of the amplitude in the neighborhood of the caustic point, we set the ray amplitudes $A(x(t,x_0))$ to 0 for all times for which $|q(t; x_0)| < \epsilon$ with a small parameter $\epsilon = q(0; x_0)$. 
%Unfortunately, \eqref{eq:AMaslov} still contains the inverse ray density $1/q(t; x_0)$, which blows up at the caustic and could cause a significant error in numerical integration of the amplitude in the neighborhood of the caustic point. 
%\rA{To prevent this, we set to 0 the ray amplitudes $A(x(t,x_0))$ for all times for which $|q(t; x_0)| < \epsilon$ with a small parameter $\epsilon = q(0; x_0)$. 
Once the determinant magnitude becomes large enough $|q(t; x_0)| \geq \epsilon$, we compute the amplitudes $A(x(t,x_0))$ behind the caustic according to \eqref{eq:AMaslov} neglecting the complex factor $\exp(-i m(t) \pi/2)$ which is essentially a phase shift experienced by the rays that cross the caustic.
%suggesting that one could recover the amplitude after the caustic.
%However, this phase shift term lacks direct translation to our time domain formulation. 
Taking the absolute value of $q(t; x_0)$ in \eqref{eq:AMaslov} annihilates the reversing of the associated ray tube coordinate system due to the caustic.}
%However, this seems not to be the case in our 2d simulations. 
A solution which holds at and beyond the caustic could be obtained e.g.~using dynamic ray tracing or Gaussian beams, which are out of scope of this manuscript but will be explored in future work. 
%where $m(t)$ represents the number of times the ray crosses a caustic (see e.g.~ \cite{runborg2007mathematical}). It holds valid after the caustic but unfortunately not on the caustic itself, which might be a significant limitation for this improvement.

%========================================
% Shadow Regions
%========================================
\subsubsection*{Shadow regions}

%\Kiko{\begin{definition}
%We say that $c(x)$ satisfies the \textit{non-trapping condition} in the domain $\Omega \subset \mathbb{R}^d$ if all ray trajectories from any point $x\in\Omega$ leave the domain in finite time.
%\end{definition}}

%\begin{definition}
%We say that the domain $\Omega$ is \textit{simple} if there exists at least one ray that travels from $x$ to $y$ 
%for any given $x\neq y \in \Omega$.
%\end{definition}

In geometrical optics a shadow region is defined with respect to a point (here we choose $x_0 \not\in \Omega$). A shadow region is a subset of the domain $\Omega_{x_0}^S \subset \Omega$ which cannot be reached along any ray trajectory originating from $x_0$, 
in other words we cannot ``illuminate'' it with rays shot from $x_0$. 
Consequently, any pressure propagating from a shadow region will not affect the pressure at the sensor located at $x_0$ when using the Hamilton-Green solver.
%: the wavefront from a shadow region never reaches the sensor. Since the wavefront is the only energy captured by the HG solver, we would like to avoid this situation.

%Therefore, in the sequel we assume that the domain $\Omega$ is \textit{simple},  that is for any two points in the domain $x\neq y \in \Omega$ there exists at least one ray trajectory from $x$ to $y$. Due to relatively small variation of the sound speed in soft tissue, with a high probability a domain in a typical PAT experiment will be simple. 
\rB{Therefore, in the sequel we assume that the sound speed $c(x)$ in $\Omega$ is \textit{non-trapping} \cite{kuchment2011mathematics},  that is for any phase space point $(x_0, p_0)$ there exists a finite time $T$ such that the ray $x(t, x_0)$ with shooting angle $p_0$, satisfies $x(t, x_0) \notin \Omega\, \; \forall T > t$. 
Due to relatively small variation of the sound speed in soft tissue, with a high probability a domain in a typical PAT experiment will be non-trapping.}
However, in the numerical simulations the issue can become more subtle. 
Even in absence of a shadow region there can be a region  with only few rays.
Such low penetration region can in turn become an effective numerical shadow region, when we cannot illuminate it with rays given a fixed numerical precision. 
In such case, adaptive techniques should be deployed to ensure that all the domain is covered with rays sufficiently densely.
Such methods could employ inserting new rays into the cone if divergence is observed. For inspirational ideas we refer to the literature on wavefront tracking methods \cite{sethian1996fast,sethian1999fast,vinje1993traveltime}. It is not immediately clear how adaptive methods can be implemented maintaining the efficiency of the one shot methods discussed here and a detailed study will be subject of future research.

\section{Hamilton-Green solver for the forward and adjoint problems}\label{sec:solver}
%Ray tracing based solution of forward and adjoint problem in PAT

%============================================================
% Discretization
%============================================================

\dontshow{
In this section we explain how we use the ray tracing approach to solve the forward and adjoint problems in PAT. 
We term the solver Hamilton-Green as our approach can be viewed as a ray tracing approximation (i.e.~numerical solution of the Hamiltonian system \eqref{req:hamiltonianEta}) to the unknown heterogeneous Green's function in the Green's integral formula for the solution of the respective forward and adjoint wave equations evaluated using the ray induced discretization of the domain. 
We note that the PAT sensors provide the natural shooting points for the rays. As shown in section \ref{sec:RevRay} the same trajectories can be used in both the forward and the adjoint problem. In particular, in the proposed Hamilton-Green solver, the adjoint problem solution follows the ray trajectories while the forward problem solution the reversed ray trajectories.}

\rB{
In this section we derive the Hamilton-Green solver for the forward and adjoint problems in PAT. Our approach can be viewed as approximating the time derivative of the unknown heterogeneous Green's function in the Green's integral solution of the respective forward/adjoint wave equation along the rays via the ray equations.  
We note that the PAT sensors provide the natural shooting points for the rays. As shown in section \ref{sec:RevRay} the same trajectories can be used in both the forward and the adjoint problem. In particular, in the proposed Hamilton-Green solver, the adjoint problem solution follows the ray trajectories while the forward formula uses the reversed ray trajectories.}

\rB{
We see the following benefits of our formulation against direct approximation of the pressure along rays: i) the resulting integral formulations highlight the parallels to the forward and adjoint/inverse transforms in ray based tomography e.g.~X-ray CT, ii) the Green integral in ray coordinates is a convenient and natural formulation for PAT where we integrate along isotemporal surfaces, iii) our formulation allows to seamlessly plug in any numerical approximation to the homogeneous Green's function $G_0$ which makes the ray approach more directly comparable to the numerical method used to obtain that approximation to $G_0$ and makes explicit the point source normalisation used, iv) the asymptotic error $\mathcal O(1/\omega)$ is confined to the approximation of the time derivative of the Green's function $\tilde G':= \partial G/\partial t$, $\tilde G' = G' + \mathcal O(1/\omega)$ which leads to a higher order asymptotic error for the solution of both the forward and adjoint problems than when these solutions are directly approximated using ray equations.  
}

% \subsection{Discretization}
% Numerical solution of the wave equation requires spatial and temporal discretization. 
% Typically the spatial domain $\Omega$ is discretized using an equispaced rectangular  
% grid of points $\Omega_D$ of the form $\omega_{m,n}= (m\Delta x, n\Delta y)$ with $m, n\in\mathbb Z$
% and $\Delta x, \Delta y$ some fixed positive quantities associated to each coordinate named \textit{grid spacing}. 
% Additionally, time is also discretized taking points of the form $j\Delta t$ with $j\in\mathbb N$ and $\Delta t\in \mathbb R^+$, 
% known as the \textit{time step}.

% From this point and for the remainder of the paper we assume that the information given for solving the forward and inverse problems
% in PAT is of the form of an equi\-spaced rectangular grid for the spatial domain and evenly discretized for the time domain.
% Moreover, we consider only rectangular spatial domains $\Omega$ for the sake of simplicity. 
% This assumption is not very restrictive, since the acoustic propagation occurs in the free space and therefore we can arbitrarily establish the limits of the domain. 

% In the following we are going to adapt the ray tracing approach to the forward and adjoint problems in PAT. 

%============================================================
% Response Function      
%============================================================
\subsection{Green's solution to the wave equation}

%The fundamental solution of the wave equation plays a key role in our RT solver. In this subsection we review the so called \textit{magic rule} \cite{barton1989elements}, which expresses the solution to the equation in terms of the associated Green's function and the set of boundary conditions (BCs) and initial conditions (ICs).

The general form of the fundamental solution of the wave equation for bounded and unbounded domains can be found e.g.~in \cite{barton1989elements}. Here we restrict the presentation to the initial value and the source problem for the unbounded domain which are relevant to PAT, see section \ref{sec:intro}. 
%\rev{The statements about homogeneous boundary conditions were not right. We are using the free space solutions, BC are not defined in the way as you wrote, see Chapter 11 Barton. I changed a lot in this section 1) to restrict to what we need, 2) to be consistent with the rest of the paper, 3) there were some mistakes in the formula for f and the problem with free space function.}

The Green's solution to the initial value problem for the inhomogeneous wave equation with the density source term $\rho$ 
\begin{subequations}
\begin{align}\label{RT:inhomogeneous}
%\square^2 u:= \left(\frac{1}{c^2(x) }\frac{\partial^2}{\partial t^2} - \Delta \right) u &= \rho(t, x),\\
\square^2 u&= \rho(t, x),\\
u(0, x) &= \psi(x),\\
u_t(0, x) &= \psi_0(x),
\end{align}
\end{subequations}
can be written as%, on a volume $V\subset \R^d$ with boundary $S = \partial V$ is of the form 
\begin{equation}\label{RT:sum}
u(t, x) = f(t, x) + h_1(t, x) + h_2(t, x)
\end{equation}
with
\begin{subequations}
\begin{align}
\label{eq:GreenSource} f(t, x) &= \int_{0}^{t^+} \int_{\R^d} G(t, x \mid t', x') \rho(t',x') \,\D x'\, \D t', \\
h_1(t, x) & = \int_{\R^d} \frac{1}{c^2(x')} G(t, x\mid 0, x') \rB{\psi_{0}(x'})\, \D x',\\
\label{eq:GreenInitP}  h_2(t, x) & = \int_{\R^d} \frac{1}{c^2(x')} \frac{\partial G}{\partial t} (t, x\mid 0, x') \rB{\psi(x')}\, \D x'.
\end{align}
\end{subequations}
Here $G$ is the corresponding free space Green's function in $\mathbb{R}^d$ and $f$ and $h = h_1 + h_2$ independently solve
\begin{itemize}
\item[$f$]: %the inhomogeneous equation \eqref{RT:inhomogeneous} with homogeneous initial conditions 
$\square^2 f = \rho(t, x), \quad f(0, x) = 0,  \; f_t(0, x) = 0,$
\item[$h$]: %the homogeneous equation with initial conditions $u(t_0, x)$ and $u_t(t_0, x)$
$\square^2 h = 0, \quad h(0, x) = \psi(x),  \; h_t(0, x) = \psi_0(x).$
\end{itemize}
%\Marta{For the initial conditions $\psi, \psi_0$ compactly supported on a bounded domain $\Omega$} and the source $\rho$ compactly supported on $\Omega \times (0, T)$, both $f$ and $h$ and hence $u$ vanish at infinity as does the free space Green's function. \todo{I wonder we need this statement at all? I am also not sure, this is actually true?}

In acoustically homogeneous medium the free space Green's functions can be calculated analytically
\begin{equation}\label{RT:greenFunc2D}
\R^2:\quad G_0^2(t, x \mid t', x') = \frac{c}{2\pi}\frac{H(c(t-t') - \|x-x'\|)}{\sqrt{c^2(t-t')^2 - \|x-x'\|^2}}, \qquad 
%H(t) = \left\{\begin{array}{ll} 0 & t < 0,\\
%                                1 & t \geq0.
%        \end{array}\right.
\end{equation}
where $H(t)$ is the Heaviside step function 
\begin{equation*}
H(t) = \left\{\begin{array}{ll} 0 & t < 0,\\
                                1 & t \geq0,
        \end{array}\right.
\end{equation*}
and
\begin{equation}\label{RT:greenFunc3D}
\R^3:\quad  G_0^3(t, x\mid t', x') = \frac{\delta(c(t-t') - \|x-x'\|)}{4\pi\|x-x'\|},
\end{equation}
where $\delta(t)$ denotes the Dirac delta distribution with the usual \textit{weak} definition
\begin{equation}\label{eq:delta}
\int_{-\infty}^{\infty}\delta(t)\zeta(t) dt = \zeta(0)
\end{equation}
holds for any Schwarz test function $\zeta$, see e.g.~\cite{barton1989elements}.

In subsections \ref{sec:solver:fwd}, \ref{sec:solver:adj} we will make use of the following invariances of (in general heterogeneous) free space Green's functions 
\begin{subequations}\label{eq:GreenInv}
\begin{align}
G(t, x \, | \, t', x') &= G(-t', x' \, | \, -t, x)      \quad\quad \mbox{(time reversal)}\\
G(t, x \, | \, t', x') &= G(t + T, x \, | \, t'+T, x'), \quad T > 0 \quad\quad \mbox{(time shift)}.
\end{align}
\end{subequations}
\subsection{Strong approximation to the Green's function}
In order to implement the Green's formula within the ray tracing algorithm we need a strong approximation of the Green's function and its derivative. For strict comparability, in the numerical experiments in section \ref{sec:simulations} we are using the Green's function obtained numerically using \texttt{k-Wave} \cite{treeby2010k}. We note that such obtained Green's function is not strictly shift invariant which will contribute to the discrepancies between the results obtained with \texttt{k-Wave} and the proposed Hamilton-Green solver.

%============================================================
% Forward Problem       
%============================================================
\subsection{Forward problem}\label{sec:solver:fwd}
%The propagation of the initial pressure wave $u_0(x)$ inside a domain $\mathbb{R}^d$ is modeled by the wave equation \rev{we repeat this a lot...}
We recall the PAT forward problem introduced in \eqref{eq:PATfwd}
\begin{subequations} 
\begin{align} 
%    u_{tt} - c^2(x)\Delta u & = 0, & (t, x)\in (0, \infty)\times\mathbb{R}^d, \tag{\ref{eq:PATfwd-1}}\\[5pt]
     \square^2 u & = 0, & (t, x)\in (0, \infty)\times\mathbb{R}^d, \tag{\ref{eq:PATfwd-1}}\\[5pt]
    u_t(0, x)              & = 0,  & x\in\mathbb{R}^d, \tag{\ref{eq:PATfwd-2}}\\[5pt]  
    u(0, x)                & = u_0(x), & x\in\mathbb{R}^d, \tag{\ref{eq:PATfwd-3}}
\end{align}
\end{subequations}
which is an initial value problem for the wave equation with heterogeneous sound speed $c \in \mathcal C^{\infty}(\R^d)$.

In the Hamilton-Green forward solver the time dependent pressure is obtained at each of the sensors separately using the cone of rays originating from this particular sensor.
%The forward simulation entails computing the pressure time series induced by the initial pressure $u_0(x)$ at a sensor located at $x_0$. 
Using Green's formula \eqref{eq:GreenInitP} the pressure $u(t, x_0)$ at a given time $t$ at a sensor location point $x_0$ can be written as
\begin{equation}\label{rt:forward}
u(t, x_0) = \int_{\mathbb{R}^d} \frac{1}{c^2(x')} \frac{\partial}{\partial t}  G(t, x_0 \mid 0, x')u_0(x')\D x',
\end{equation}
where $G(t, x \mid t', x')$ is an, in general, heterogeneous free-space Green's function. Furthermore, for any $u_0$ compactly supported on $\Omega$ the integration can be restricted to the domain $\Omega$. 
%Assuming that the rays shot from $x_0 \in \Omega$ do not develop a caustic in the domain $\Omega$ i.e.~no two different rays from $x_0$ intersect in $\Omega$, and that they cover $\Omega$ sufficiently densely, these rays form a coordinate system for $\Omega$.
\rA{Let $(\ell, \theta) \in \R_{0+} \times \mathcal S^{d-1}$ be ray based coordinates for the rays originating from $x_0 \not\in \Omega$. Assuming that these rays do not intersect in $\Omega$, they form a coordinate system in $\Omega$ and we can perform a formal change of coordinates in \eqref{rt:forward}. However, the ray based coordinates parametrisation yields a useful approximation to \eqref{rt:forward} even if the rays intersect. 
In this parametrisation \eqref{rt:forward} becomes}
\begin{align}\label{req:sphericalIntegral}
u(t, x_0) =  \int_0^T \int_{\mathcal S^{d-1}} 
\frac{1}{c^2(x(\ell, \theta; x_0))}
 \frac{\partial}{\partial t}G(t, x_0 \mid 0, x(\ell, \theta; x_0))\, u_0(x(\ell, \theta; x_0)) \left|\frac{\D x'}{\D(\ell, \theta)}\right|\, \D(\ell, \theta).
\end{align}
%Here $(\ell, \theta) \in \R_{0+} \times \mathcal S^{d-1}$ are ray based coordinates for the rays originating from $x_0$.
\rA{The acquisition time $T \in \R_+$ is chosen sufficiently large so that all the rays in the cone have left the domain i.e.~$x(\ell, \theta; x_0) \notin \Omega \,\,\forall \ell > T, \theta \in \mathcal S^{d-1}$.} Furthermore, by the definition, the determinant of the Jacobian of the change of coordinates 
\begin{equation}\label{eq:dxdltheta}
\left|\frac{\D x}{\D (\ell, \theta)} \right|_{x(0)=x_0} = q(\ell, \theta; x_0).
\end{equation}
%immediately from the definition  
%\begin{equation} q(\ell, \theta; x_0) = \left|\frac{\partial x}{\partial \ell} \,\, \frac{\partial x}{ \partial \theta}\right|_{x(0) = x_0}.
%\end{equation}
Using the combined time reversal and shift invariances of Green's functions \eqref{eq:GreenInv} we can rewrite the Green's function with the source at the sensor 
\begin{equation}\label{eq:GreenRev}
G(t, x_0 \, | \, 0, x(\ell,\theta; x_0)) = G(t, x(\ell,\theta; x_0) \,|\, 0, x_0).
\end{equation} 

Substituting \eqref{eq:dxdltheta} and \eqref{eq:GreenRev} into \eqref{req:sphericalIntegral} we obtain
\begin{align}\label{eq:sphericalIntegralRev}
u(t, x_0) =  \int_0^T \int_{\mathcal S^{d-1}} &
\frac{1}{c^2(x(\ell, \theta; x_0))}
 \frac{\partial}{\partial t} G(t, x(\ell,\theta; x_0) \,|\, 0, x_0) \, u_0(x(\ell, \theta; x_0)) \, q(\ell,\theta; x_0) \, \D(\ell, \theta). 
\end{align}
As the Green's function for heterogeneous domain $G$ is generally unknown, we propose the following approximation in \eqref{eq:sphericalIntegralRev}. Assuming a constant sound speed in a small neighborhood $\mathcal N(x_0)$ of $x_0$, we can locally approximate the heterogeneous Green's function $G(t,x\, |\, 0, x_0)$ with the homogeneous Green's function $G_0(t,x\, |\, 0, x_0)$ with sound speed $c_0$ equal to sound speed at the emission point $x_0$ i.e.~$c_0 = c(x_0)$, 
$$ G(t,x\, |\, 0, x_0) \approx G_0(t, x\, |\, 0, x_0), \quad x\in \mathcal N(x_0).$$ 
%\todo{i changed to approx instead of equal, because $G$ has reflections and $G_0 $ does not. We could write something like "with equality up to reflections".}
Equation \eqref{eq:sphericalIntegralRev} however, requires $G$ everywhere and the homogeneous function approximation is not valid outside of $\mathcal N(x_0)$.
\rB{
To simplify notation we short hand the attenuation along the ray from $x_0$ to $x(\ell, \theta; x_0)$ in the amplitude formula \eqref{req:amplitude} as
\begin{equation}\label{req:mu}
\mu(\ell, \theta; x_0) = \frac{\eta(x_0)}{\eta(x(\ell, \theta; x_0))} 
\sqrt{\frac{q(0, \theta; x_0)}{q(\ell, \theta; x_0)}}.
\end{equation}}
Now, the ray based coordinate system allows us to extend the approximation propagating the homogeneous free space Green's function $G_0(t,x_0\, |\, 0, x_0)$ 
along the ray $(x_0, \theta)$ from $x_0$ to $x(\ell, \theta; x_0)$ and compute the corresponding amplitude according to \eqref{req:amplitude}
\begin{equation}\label{req:approxGreens}
%G(t, x(\ell, \theta; x_0) \mid 0, x_0) \approx G_0 (t+\ell , x_0 \mid 0, x_0)  \, \mu(\ell, \theta; x_0). 
G(t, x(\ell, \theta; x_0) \mid 0, x_0) \approx G_0 (t-\ell , x_0 \mid 0, x_0)  \, \mu(\ell, \theta; x_0). 
%= G_0 (t , x_0 \mid -\ell, x_0) \, \mu(\ell, \theta; x_0).
\end{equation}
Assuming the time parametrization along the ray, $\ell$ is the time of propagation from $x_0$ to $x(\ell,\theta; x_0)$ which directly  %is used to synchronize the emission time at $x_0$ to 0.
corresponds to the phase \eqref{eq:linphase}. The latter effects a time shift in the Green's function due to its explicit dependence on time (note that throughout our derivation the time dependence has been implicit through the parametrization of the trajectory).

We note that the choice of the approximation point to be the evaluation point, here the sensor $x_0$, is not coincidental but necessary for the appropriate speed of sound $c_0$ and hence time scaling of the homogeneous Green's function $G_0$. This indeed was the reason for reversing the Green's function in \eqref{eq:GreenRev}.
Substituting \eqref{req:approxGreens} and \eqref{req:mu} into \eqref{eq:sphericalIntegralRev} and rearranging the integration to reflect functional dependencies we arrive at
\begin{align}\label{eq:GHFwdNum}
%\nonumber u(t, x_0) = \eta(x_0) \int^T_0 &\frac{\partial}{\partial t} G_0(t+\ell, x_0 \mid 0, x_0)  \\
u(t, x_0) = \eta(x_0) \int^T_0 \frac{\partial}{\partial t} G_0(t-\ell, x_0 \mid 0, x_0) 
 \int_{\mathcal S^{d-1}}  \eta(x(\ell, \theta; x_0)) \sqrt{ q(0, \theta; x_0) q (\ell, \theta; x_0) } \, u_0(x(\ell, \theta; x_0)) \,\D\theta \, \D\ell. 
\end{align}

\begin{figure}
\begin{center}
\includegraphics[width=0.6\linewidth]{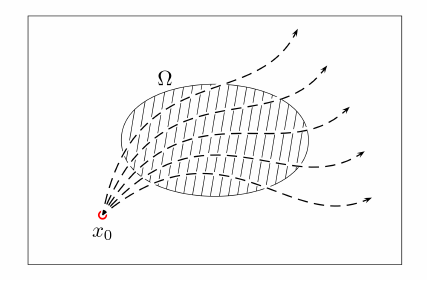}
\caption{Cone of rays from $x_0$ covering the domain $\Omega$ where the initial pressure $u_0$ is compactly supported.}
\label{initialPressure}
\end{center}
\end{figure}

Below we summarize the steps necessary to numerically evaluate \eqref{eq:GHFwdNum} which corresponds to solving
%In summary, the steps to solve 
the forward problem \eqref{eq:PATfwd} for a sensor located at $x_0$.
\begin{itemize}
\item Trace ray trajectories from $x_0$ such that they cover the domain $\Omega$ sufficiently densely (figure \ref{initialPressure}).
\item Interpolate $u_0$ and $\eta$  at equitemporal points along the ray trajectories ($q$ is already available at these points), see \ref{app:intC2R}. % using nearest-neighbor interpolation.
\item Numerically evaluate the spherical integral along the isotemporal surfaces, $\ell = \mbox{const}$. Using mid-point rule this corresponds to summing the appropriately weighted contributions from all rays from $x_0$ at $\ell$. 
%as a properly shifted signal in time using the propagation times and reverse amplitudes along them.
\item Convolve the result with the time derivative of the homogeneous Green's function at $x_0$, $\frac{\partial}{\partial t} G_0(t, x_0 \mid 0, x_0)$.
\end{itemize}

%============================================================
% Adjoint Problem
%============================================================
\subsection{Adjoint problem}\label{sec:solver:adj}

The adjoint problem in PAT \eqref{eq:PATadj} is a time varying source problem for the free space wave equation with homogeneous initial conditions. Due to the linearity of the wave equation \eqref{eq:PATadj}, the solution is a superposition of the corresponding solutions for each individual point mass source located at the sensor $x_0^m, \; m=1,\dots,M$. Similarly as for the forward problem, in the adjoint Hamilton-Green solver we solve the adjoint problem for each mass source separately using the ray cone originating from the sensor $x_0$
\begin{subequations}%\tag{\eqref{eq:PATadj-1}}
\begin{align}
%v_{tt}-c^2(x)\Delta v & = c^2(x)\, \delta(x-x_0)\,  \frac{\partial}{\partial t} \left( g(T-t, x) \omega(T-t, x) \right), & t>0,\, x\in\Omega,\\
\tag{\ref{eq:PATadj-1}'} \square^2 v &= \frac{\partial}{\partial t} \Big( g(T-t, x)\, \omega(T-t) \Big) \delta(x-x_0),  &(t, x) \in(0,T) &\times \R^d,\\
\tag{\ref{eq:PATadj-2}'} v_t(0, x) &= 0, &x&\in\R^d,\\
\tag{\ref{eq:PATadj-3}'} v(0, x) &= 0,  &x&\in\R^d.
\end{align}
\end{subequations}

Analogously as for the forward problem, we start from the corresponding Green's representation to the solution  of the source problem \eqref{eq:GreenSource} evaluated at $t = T$
\begin{equation}
v(T, x) = \int_{\R^d}\int_{0}^{T} G(T, x \mid t', x') \, \frac{\partial}{\partial t'} \Big( g(T-t', x')\,\omega(T-t') \Big) \delta(x'-x_0)\, \D t' \D x'. \label{eq:GreenAdjX}
\end{equation}
Upon the change to ray based coordinate system \eqref{eq:GreenAdjX} becomes
\begin{align}
\nonumber v(T, x(\ell,\theta; x_0)) = \int_0^T\int_{\mathcal S^{d-1}} & \underbrace{\int_{0}^{T} G(T, x \mid t', x(\ell',\theta'; x_0)) \, \frac{\partial}{\partial t'} \Big( g(T-t', x(\ell',\theta'; x_0))\,\omega(T-t') \Big) \D t'}_{(\star)}  \\ 
&\times \underbrace{\delta(\ell'-0)\delta(\theta' - \theta)}_{=\delta(x(\ell',\theta'; x_0)-x(0,\theta; x_0))} q(\ell',\theta'; x_0) \,\D \theta' \D \ell', \label{eq:GreenAdjRayD1}
\end{align}
where we chose to represent $x_0 = x(0,\theta; x_0)$ using the ray through $x = x(\ell, \theta; x_0)$.
Integrating $(\star)$ by parts, using the time shift invariance of the Green's function to set the emission time to 0 and the variable change $t = T-t'$ yields %we move the derivative to the Green's function
\begin{align}\label{eq:GreenSourceIntByParts}
(\nonumber \star) =& -\int_{0}^{T} \frac{\partial}{\partial t'} G(T, x \mid t', x(\ell',\theta'; x_0)) \, g(T-t', x(\ell',\theta'; x_0))\,\omega(T-t')  \D t' \\
\nonumber &+ G(T,x, \mid T, x(\ell',\theta'; x_0)) g(0, x(\ell',\theta'; x_0)) \underbrace{\omega(0)}_{=0}\\
\nonumber &- G(T,x, \mid 0, x(\ell',\theta'; x_0)) g(T, x(\ell',\theta'; x_0)) \underbrace{\omega(T)}_{=0}\\
=& \int_{0}^{T} \frac{\partial}{\partial t} G(t, x \mid 0, x(\ell',\theta'; x_0)) \, g(t, x(\ell',\theta'; x_0))\,\omega(t) \D t.
\end{align}
Substituting \eqref{eq:GreenSourceIntByParts} into \eqref{eq:GreenAdjRayD1} and using the distributional calculus to evaluate the ray cone integral against the $\delta(\ell'-0)\delta(\theta' - \theta)$, we obtain
\begin{equation}
v(T, x(\ell,\theta; x_0)) = q(0,\theta; x_0) \int_{0}^{T} \frac{\partial}{\partial t} G(t, x \mid 0, x_0) \, g(t,x_0)\,\omega(t) \D t.\label{eq:GreenAdjRayD2}
\end{equation}
To approximate the unknown heterogeneous Green's function we again make use of a ray based approximation setting 
$$G(t, x \mid 0, x(\ell,\theta; x_0)) \approx G_0(t, x \mid 0, x(\ell,\theta; x_0)), \quad x\in \mathcal N(x(\ell,\theta; x_0)),$$
and outside of $\mathcal N(x(\ell,\theta; x_0))$ (in particular at the sensor $x_0$)
\begin{align*}
%G(t, x_0 \mid 0, x(\ell,\theta; x_0)) \approx G_0(t + \ell, x(\ell, \theta; x_0) \mid 0, x(\ell, \theta; x_0)) \, \mu_R(\ell, \theta; x_0),
G(t, x_0 \mid 0, x(\ell,\theta; x_0)) \approx G_0(t - \ell, x(\ell, \theta; x_0) \mid 0, x(\ell, \theta; x_0)) \, \mu_R(\ell, \theta; x_0),
\end{align*}
where $\mu_R$ is the analogously defined reversed attenuation in \eqref{req:revAmplitude} i.e.~attenuation between $x(\ell, \theta; x_0)$ and $x_0$ along the reversed ray
\begin{align}\label{eq:muR}
\nonumber \mu_R(\ell, \theta; x_0) &= \frac{\eta(x(\ell, \theta; x_0))}{\eta(x_0)} 
\sqrt{\frac{q_R(0, \theta; x(\ell, \theta; x_0))}{q_R(\ell, \theta; x(\ell, \theta; x_0))}}\\
&=\frac{\eta(x(\ell, \theta; x_0))}{\eta(x_0)}\sqrt{\frac{q(\ell, \theta; x_0)}{q(0, \theta; x_0)}} \phantom{x(\ell,\theta; x_0)} = \frac{1}{\mu(\ell,\theta; x_0)}.
\end{align}
Using the combined time shift and reversal invariances of Green's function we rephrase the Green's function in \eqref{eq:GreenAdjRayD2} with a source at the evaluation point $x(\ell, \theta; x_0)$
\begin{equation}
v(T, x(\ell,\theta; x_0)) = q(0,\theta; x_0) \int_{0}^{T} \frac{\partial}{\partial t} G(t, x_0 \mid 0, x(\ell, \theta; x_0)) \, g(t,x_0)\,\omega(t) \D t  \label{eq:GreenAdjRay}
\end{equation}
which will guarantee the correct speed of sound $c_0 = c(x(\ell, \theta; x_0)))$ for the homogeneous Green's function $G_0$. Substituting the ray approximation and $\mu_R$ from \eqref{eq:muR} we arrive at the Hamilton-Green approximation to the adjoint solution
\begin{align}
v(T, x(\ell,\theta; x_0)) %\approx & q(0,\theta; x_0) \\
\approx \frac{\eta(x(\ell, \theta; x_0))}{\eta(x_0)}\sqrt{q(0,\theta; x_0) q(\ell,\theta; x_0)}
\label{eq:GHAdj}  \int_{0}^{T} \frac{\partial}{\partial t} G_0(t - \ell, x(\ell, \theta; x_0) \mid 0, x(\ell, \theta; x_0)) \, g(t,x_0)\,\omega(t) \D t.
\end{align}
Equation \eqref{eq:GHAdj} is naturally defined in ray based coordinates and needs to be interpolated on the Cartesian grid as described in the \ref{app:intR2C}.  
Summing over all sensors we obtain the adjoint of the original equation \eqref{eq:PATadj}
\begin{align}\label{eq:GHAdjFull}
v(T, x) \approx& \sum_{m=1}^M v(T, x(\ell^m,\theta^m; x_0^m)).
\end{align}
%\begin{align}\label{eq:GreenAdj}
%v(T, x) \approx& \sum_{m=1}^M q(0,\theta^m; x_0^m) \\
%&\times \int_{0}^{T} \frac{\partial}{\partial t} G_0(t + \ell^m, x(\ell^m, \theta^m; x_0^m) \mid 0, x(\ell^m, \theta^m; x_0^m)) \mu_R(\ell^m, \theta^m; x_0^m) \, g(t,x_0^m)\,\omega(t) \D t. 
%\end{align}
\noindent The steps necessary for numerical evaluation of \eqref{eq:GHAdj} which corresponds to the adjoint problem with one sensor located at $x_0$ are summarized below. We note that steps 1,2 are the same as for the forward problem and will be executed only once in an efficient implementation. 
\begin{itemize}
\item Trace ray trajectories from $x_0$ such that they cover the domain $\Omega$ sufficiently densely (figure \ref{initialPressure}).
\item Interpolate $\eta$ at equitemporal points in the domain ($q$ is already available at these points).
\item Numerically evaluate the time integral. Using mid-point rule this corresponds to summing the data time series weighted by a shifted derivative of the homogeneous Green's function $G_0$ with a sound speed $c_0 = c(x(\ell,\theta; x_0))$. The details of efficient evaluation of $G_0$ are given in the \ref{app:G0}.
\end{itemize}

\rB{\section{Error analysis of the Hamilton-Green solver}\label{sec:errana}}
{\color{myBlue}

In this section we identify and quantify various contributions to the error of the Hamilton-Green approximation to the solution of the forward and adjoint PAT problems. 

\paragraph{Asymptotic error in approximation of $G' := \partial G/\partial t$}
The HG solver uses the first term of the WKB expansion to approximate the time derivative of the unknown (heterogeneous) Green's function $G$, $G' := \partial G/\partial t$, on the domain along the rays. 
This implies that the asymptotic error of the approximation of $G'$ is $\tilde G' = G' + \mathcal O(1/\omega)$. 

The solution $u$ of the forward problem \eqref{eq:GHFwdNum} is a convolution of the derivative of the Green's function $G'$ with the initial pressure $u_0$. Using the approximate $\tilde G'$ instead, results in an approximation $\tilde u$. Switching to Fourier domain, we have the following asymptotic for this approximate solution $\Ft{\tilde u} = \Ft{\tilde G'} \Ft{u_0} = \left(\Ft{G'} +\mathcal O(1/\omega)\right) \Ft{u_0} = \Ft{G'} \Ft{u_0} + \mathcal O(1/\omega) \Ft{u_0} = \Ft{u} + \mathcal O(1/\omega) \Ft{u_0}$ as opposed to $\Ft{\breve u} = \Ft{u} + \mathcal O(1/\omega)$ when the solution is approximated directly by the ray equations. Making a conservative assumption that $u_0$ is compactly supported of bounded variation, implies the decay of its Fourier transform $| \Ft{u_0(\xi)}| = \mathcal O(1/|\xi|)$ as $|\xi| \rightarrow \infty$ yielding a higher order asymptotic error for the HG solution $\tilde u = u + \mathcal O(\omega^{-1}|\xi|^{-1}) = u + \mathcal O(\omega^{-2})$. 

When we assume $u_0 \in C_0^\infty(\mathbb R^d)$ (as we do here), its Fourier transform decays faster than any polynomial $|\Ft{u_0(\xi)}| =  o(|\xi|^{-n}), \; \forall n\in \mathbb N$ as $|\xi| \rightarrow \infty$ and consequently so does the asymptotic error of the HG solution $\tilde u = u + o(\omega^{-n})$. 

As the asymptotic error of $\tilde G'$ enters the adjoint problem \eqref{eq:GHAdj} in the same way i.e.~the solution of \eqref{eq:GHAdj} also contains the convolution of $G'$ but this time with the PAT data $g$, the same argument holds for the asymptotic error of the solution to the adjoint problem with analogous smoothness assumptions on data $g$. 

\paragraph{Accuracy of numerical integration of the ray trajectories}
The solution to the Hamiltonian system \eqref{req:hamiltonianEta} results in ray trajectories $x(t)$, which are an essential building block of our solver. 
These ray trajectories are computed using a second order Runge-Kutta method.
Assuming that the exact sound speed is known, the error commited at each step is $\mathcal O(h^3)$, where $h$ is the step size \cite{ascher1998computer}. 
Therefore, the accumulated overall error for the whole trajectory is $\mathcal O(h^2)$.

\paragraph{Grid interpolation error}
Another source of error is the interpolation between Cartesian and ray grids.
We use nearest-neighbour interpolation in both directions: from Cartesian to ray and from ray to Cartesian (see \ref{app}).
For a Cartesian grid with spacing $h$, the order of the error using nearest-neighbour interpolation is $\mathcal O(h)$. 
Conversely, the order of the error when converting from the ray grid to a Cartesian grid is $\mathcal O(h)$, where $h$ describes the maximum distance between neighbouring points along the trajectory or neighbouring points on the wavefront.

\paragraph{Accuracy of approximation of free space Green's function $G_0$ and its time derivative $G_0'$}
%The Green's integral formula 
The forward \eqref{eq:GHFwdNum} and adjoint \eqref{eq:GHAdj}, \eqref{eq:GHAdjFull} formulae require an approximation of the free space Green's function $G_0$ or its time derivative $G_0'$ . 
The free space Green's function is a distribution with a 
singularity at the origin (cf.~equations \eqref{RT:greenFunc2D} and \eqref{RT:greenFunc3D}) %a numerical approximation of such distribution with arbitrary precision is not possible. 
which we compute as an impulse response to a strong approximation to $\delta$. The so approximated $G_0$ is has no singularities and is finite (when e.g.~using Gaussian strong approximation with small tail truncation).   
Alternatively, $G_0$ can be obtained solving the wave equation with a numerical method and an approximation of $\delta$ consistent with the method, see also \ref{app:G0}.
%approximation of the Green's function, which is finite with no singularities.
As the $G_0$ approximation is obtained independently, we can tune it to be of the same or lower order as the interpolation errors.
The approximation of the derivative $G_0'$ follows analogously. 
%As the frequency range of the solution is already limited by the chosen discretisation, this $G_0$ approximation error can be tuned to be 

\paragraph{Error resulting from reflections, caustics and shadow regions}
When a wave traverses an interface between two media, part of its energy is reflected back and part of its energy is transmitted. The strength of the reflection depends on the smoothness of the sound speed respective to the wave length.
The Hamilton-Green solver does not account for the reflected part, since it only propagates the pressure at the wavefront. 
As discussed in Subsection \ref{re:sec:caustics} the presence of caustics also limits the accuracy of our solver, since the determinant becomes singular and amplitude blows up at caustics.
To control this error, we set to 0 the amplitudes along the rays in the neighbourhood of a caustic region.
As no rays pass through numerical shadow regions, no contribution to pressure is recorded from these regions. These regions however, naturally correspond to very low ray densities and their contribution should be minimal. 

In \ref{app:lens} we study the effect of the caustics in an acoustic lens which however highly exaggerates the issue. We performed no dedicated studies bounding the error introduced by reflections and shadow regions. 
%There are no studies bounding the error introduced by reflections, caustics and shadow regions. 
However, we know that caustics, reflections and shadow regions are absent in the homogeneous case thus there we have control over all other errors (ray trajectories, grid interpolation and Green's function), our solver converges in the homogeneous case at the order determined by the lowest order of the numerical scheme. 
Finally, we expect the errors introduced by reflections, caustics and shadow regions to be small when the changes of sound speed in the domain are subtle, which is in general the case in soft tissue \cite{mast2000empirical}.

%==============================
% 5. Simulations and comparison with k-Wave
%==============================
%==============================
% Simulations and comparison with k-Wave
%==============================
\dontshow{
\section{Simulations and Comparison with k-Wave}
\label{sec:simulations}
In this section we solve the forward and adjoint problems of a phantom that emulates veins in a soft tissue, an acoustically
heterogeneous rectangular domain $\Omega$.
This domain is discretized using $128\times256$ equispaced grid points along the $x$ and $y$ axes, with distance 
between points equal to $0.2$ [mm].
The sound speed $c$ in the domain (figure \ref{sim:domain} left) ranges from 1,350 to 1,650 [m/s], and therefore contains the typical values for sound speeds
in soft tissue, which are around $1,500 \pm 5\%$ [m/s] \cite{mast2000empirical}.
The phantom that we have chosen (figure \ref{sim:domain} right) is a simplified version of the type of image that we would like to obtain when using PAT.
It represents some veins in the tissue, which are the initial pressure wave $u_0(x)$ that propagates in the domain.
}

\section{Evaluation of the Hamilton-Green solver}\label{sec:simulations}
\rA{
In this section we evaluate the accuracy of the Hamilton-Green solver (HG) for the forward and adjoint PAT problems. We first discuss the numerical convergence of the solver for homogeneous sound speeds and then move to heterogeneous sound speed phantom emulating vessels embedded in a soft tissue. 
The HG solutions are compared against the first order pseudospectral method implemented in \texttt{k-Wave} toolbox \cite{treeby2010k}.
}
\rB{The acoustic solver in \texttt{k-Wave} was extensively validated and its performance compared against analytical solutions \cite{treeby2012modeling, martin2016simulating}, other numerical models \cite{demi2014comparison, robertson2017accurate} and experimental measurements \cite{Robertson_2017, wang2012modelling}.}
Our examples are restricted to 2D for ease of visualisation.
%The 2D Matlab implementation used in this paper can be downloaded from the Github repository  \href{https://github.com/kikorulan/green-ray}{\texttt{green-ray}}.
An efficient GPU implementation enabling handling large image volumes will be the subject of a sequel paper along with the integration of the solver
with iterative methods which can take full benefit of its flexibility.
Therefore, we defer the discussion of the computational efficiency to that later work.

%==============================
% SIMULATIONS - HETEROGENEOUS
%==============================

%====================================================================================================
% ERROR ANALYSIS FOR FORWARD PROBLEM
%====================================================================================================
\rA{
\subsection{Error Analysis of the Forward Problem}\label{sec:sim:fwd}}
{\color{myGreen}

% ball positions
% - 6.4 [mm] x1 x2 x3
% - 10.25 [mm] y1
% - 15.35 [mm] y2
% - 20.48 [mm] y3
% Sigma 400 [mum]
We consider a discretised rectangular domain $\Omega$ of dimensions 512$\times$1024 with grid spacing $\Delta x = 25$ [$\mu$m] and homogeneous sound speed $c_0 = 1,500$ [m/s].
We compute the solutions for three smooth initial pressure conditions $u_i$, $i\in\{1, 2, 3\}$ %, each of them corresponding to a Gaussian,
\begin{equation}
u_i(x, y) = e^{-\frac{(x-x_i)^2 + (y-y_i)^2}{\sigma^2}},
\end{equation}
where $(x_1, y_1) = (6.35, 10.25)$ [mm], $(x_2, y_2) = (6.35, 15.35)$ [mm] and $(x_3, y_3) = (6.35, 20.48)$ [mm] and $\sigma = 400$ [$\mu$m].
The Gaussian serves as a smooth alternative to a ball allowing to bypass the use of the smoothing filter in \texttt{k-Wave} simulations\footnote{http://www.k-wave.org/manual/k-wave\_user\_manual\_1.1.pdf} 
and reduce the errors at the boundary. For compactness of visualisation, a superposition of these individual initial pressure conditions, $\sum_i u_i$, is shown in figure \ref{sim:fig:InitialPressure}.

\begin{figure}
\begin{center}
\subfloat[][]{\includegraphics[height = 7.5cm]{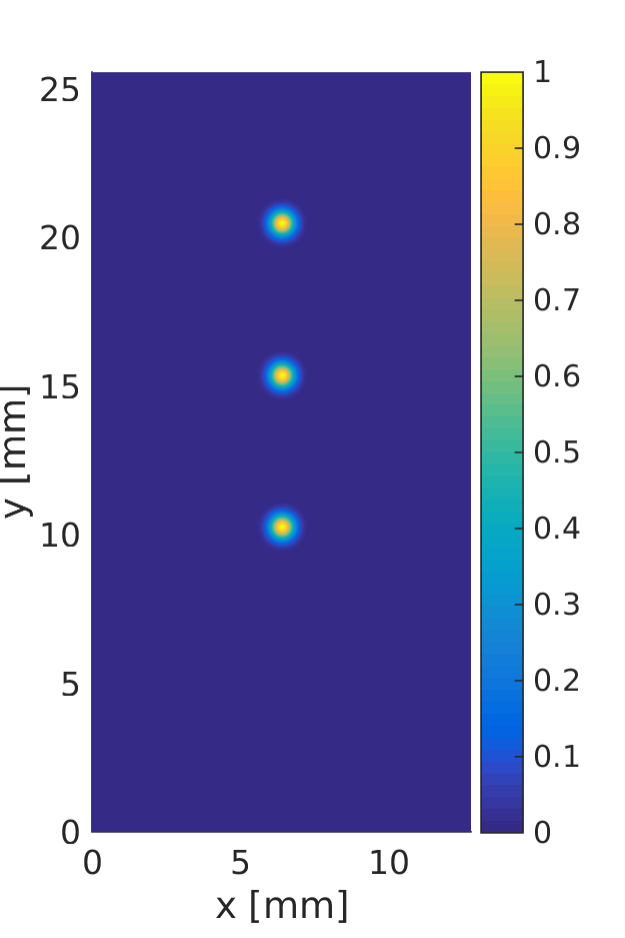}\label{sim:fig:InitialPressure}} \hspace*{0.2cm}
\subfloat[][]{\raisebox{0.2cm}{\includegraphics[height = 6.8cm]{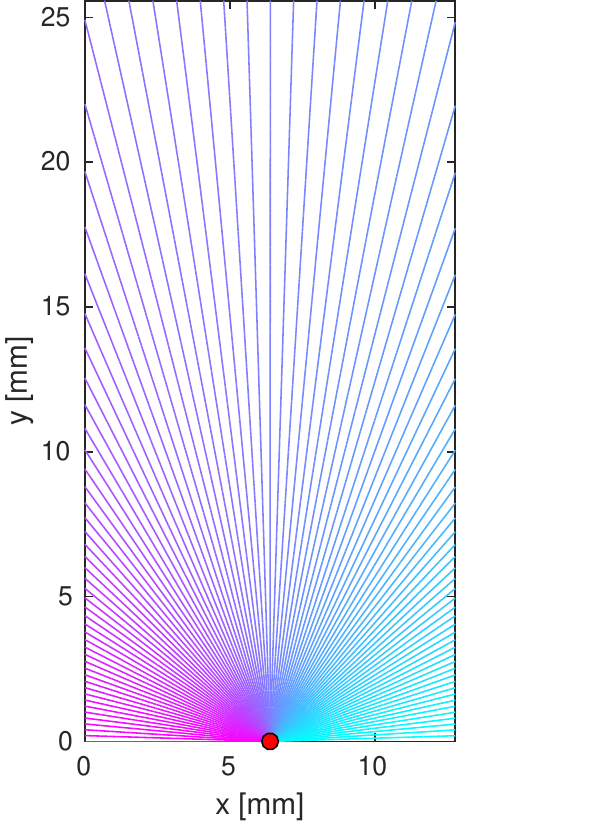}}\label{sim:fig:Rays_homo}}
\caption{Homogeneous sound speed. (a) Sum of initial pressure conditions $\sum u_i$. (b) Subset of ray trajectories from  $x_0 = (6.35, 0)$ [mm].}
\end{center}
\end{figure}

%==================== HOMOGENEOUS SOUND SPEED
\begin{figure}
\begin{center}
\subfloat[][{HG vs. \texttt{k-Wave} for $\Delta t$ = 8 [ns].}]{\includegraphics[width = 6.8cm]{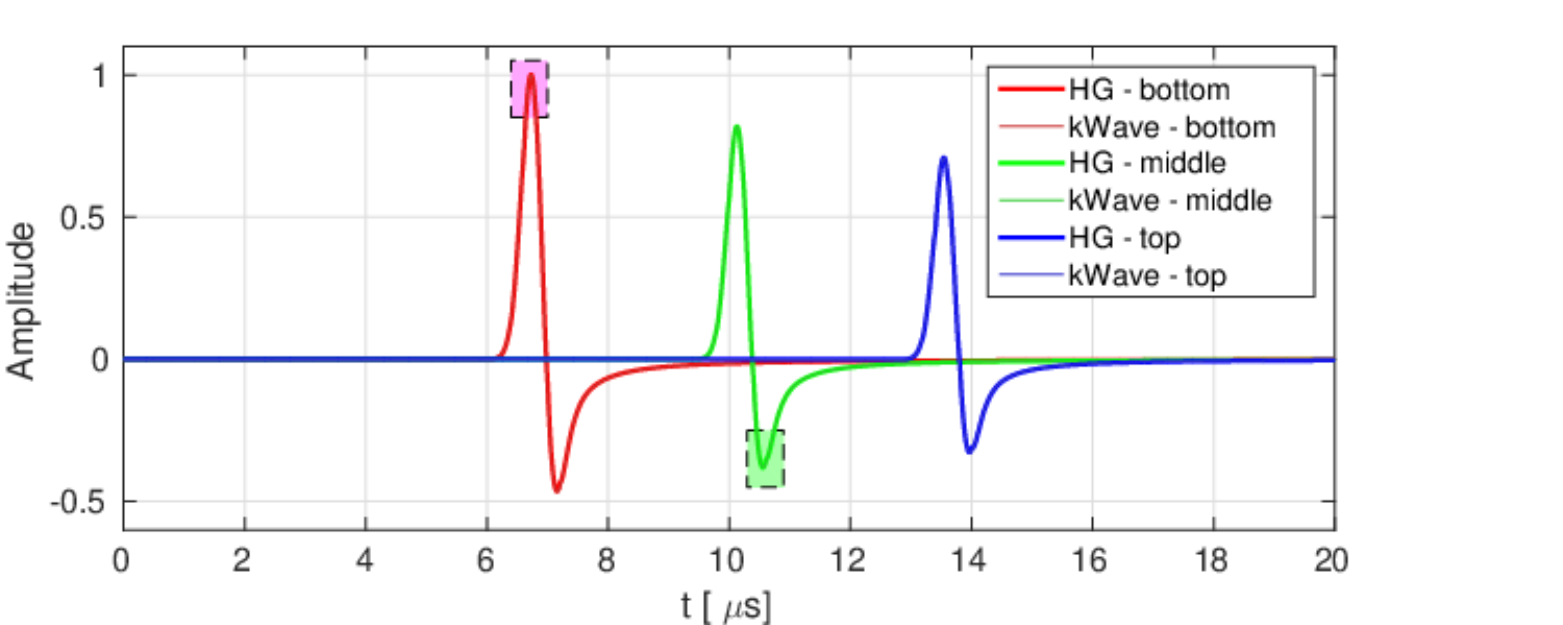} \label{sim:fig:HGvsKWAVE_homo}} \hspace*{0.2cm}
\subfloat[][{$\Delta t = 8$ [ns].}]{\includegraphics[width = 6.8cm]{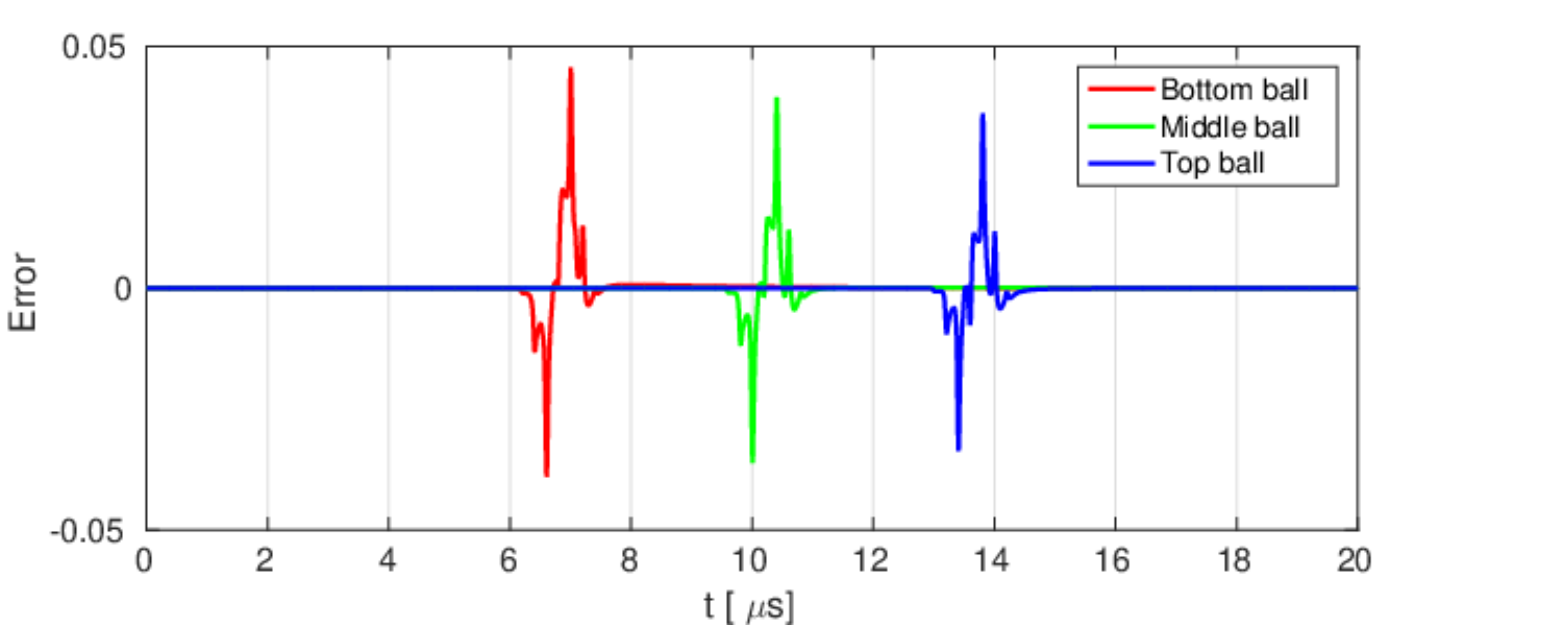} \label{sim:fig:HGvsKWAVE_homo_8}}\\[10pt]
\subfloat[][{$\Delta t = 4$ [ns].}]{\includegraphics[width = 6.8cm]{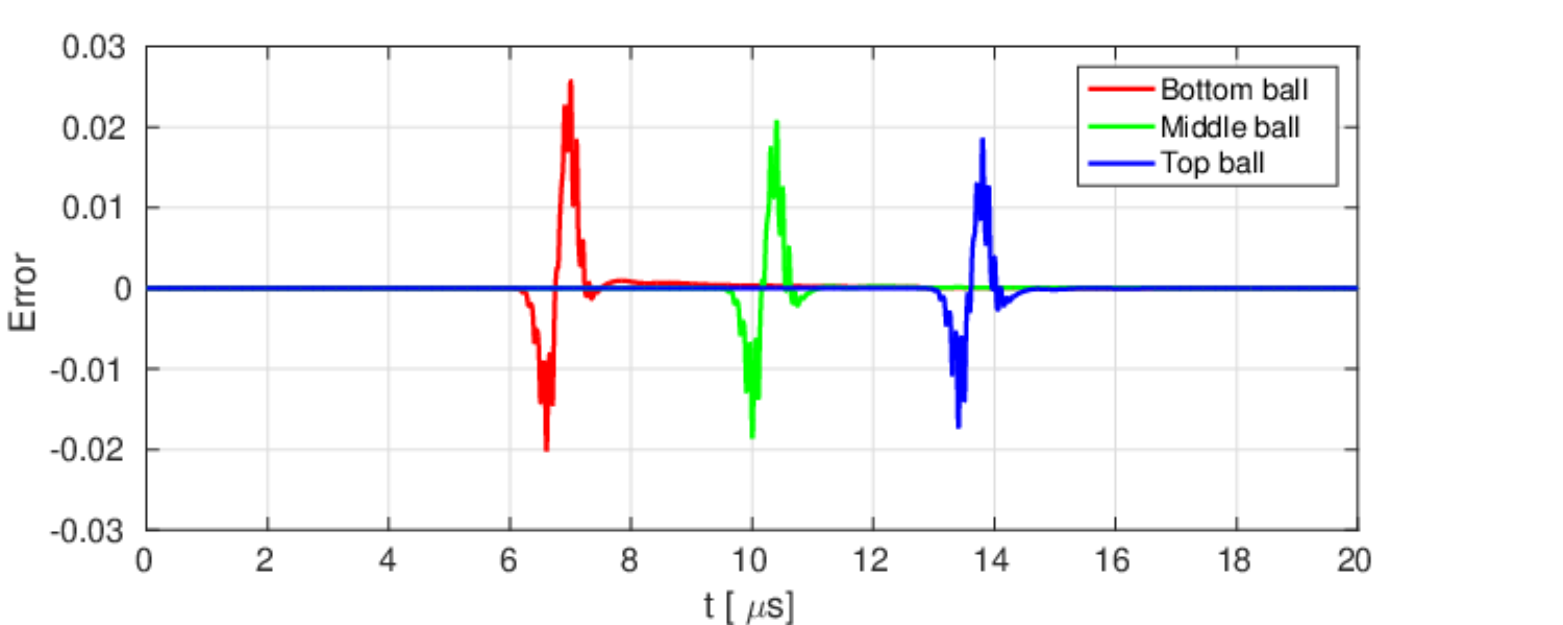} \label{sim:fig:HGvsKWAVE_homo_4}}\hspace*{0.2cm}
\subfloat[][{$\Delta t = 2$ [ns].}]{\includegraphics[width = 6.8cm]{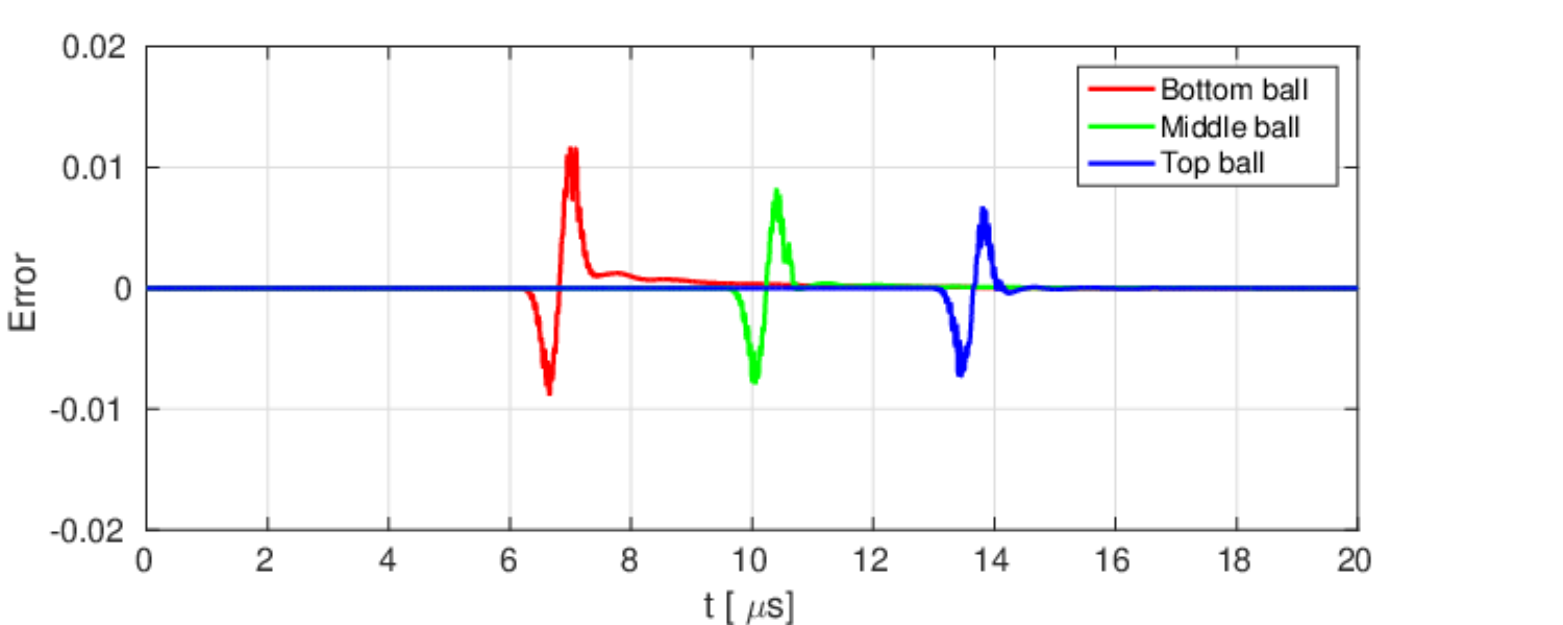} \label{sim:fig:HGvsKWAVE_homo_2}}\\[10pt]
\hspace*{-0.2cm}\subfloat[][{HG vs. k-Wave for $\Delta t = 8 $ [ns];  \textcolor{red}{bottom} ball zoomed in.}]{\includegraphics[width = 6.8cm]{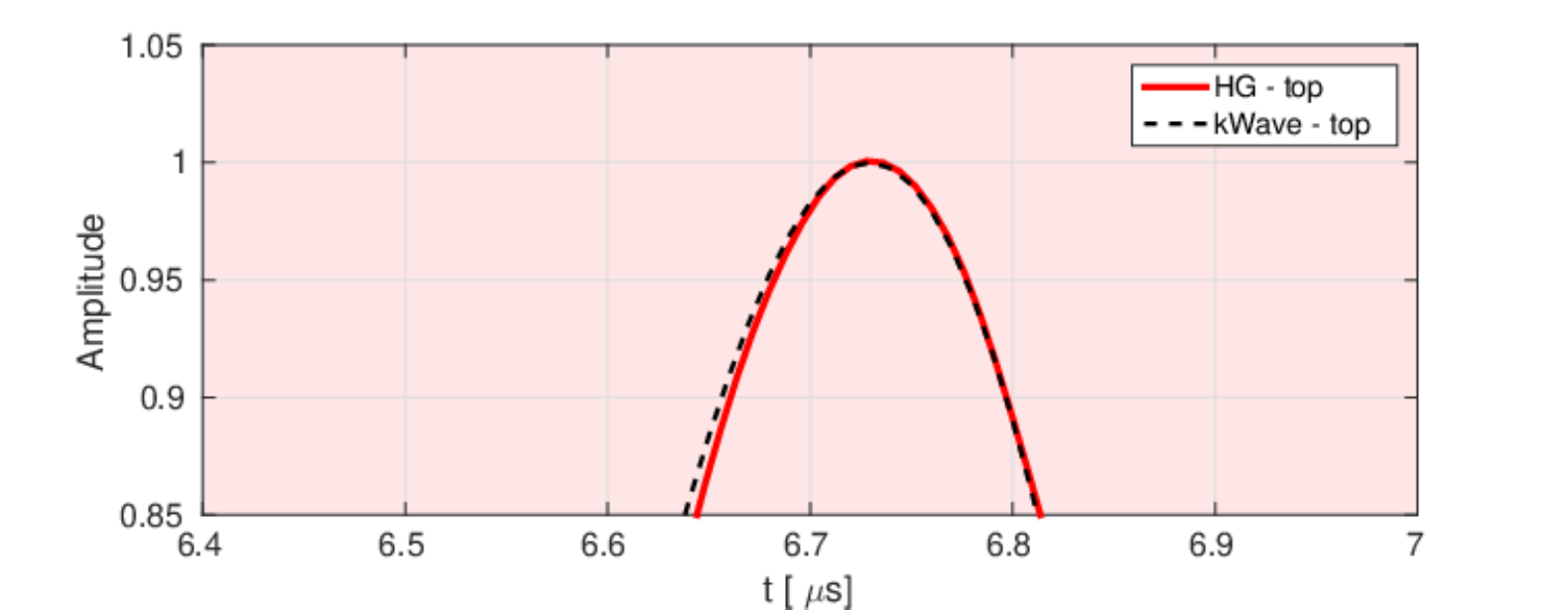} } \hspace*{0.2cm}
\subfloat[][{HG vs. k-Wave for $\Delta t = 8 $ [ns]; \textcolor{green}{middle} ball zoomed in.}]{\includegraphics[width = 6.8cm]{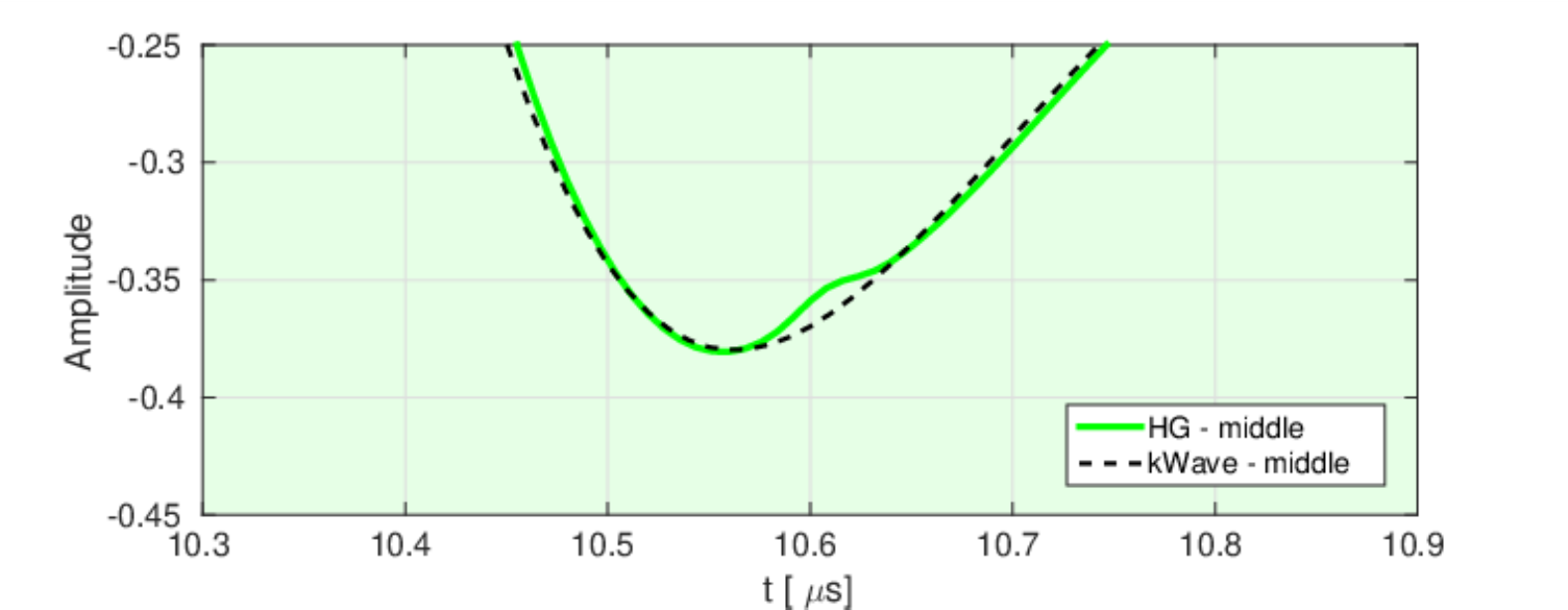}  }
\caption{Homogeneous sound speed. Comparison between HG and \texttt{k-Wave} forward solutions at $x_0$ for decreasing time step length $\Delta t \in \{8, 4, 2\}$ [ns].}
\label{sim:fig:main_HGvsKWAVE_homo}
\end{center}
\end{figure}

\begin{figure}
\begin{center}
\subfloat[][REE.]{\includegraphics[width = 6cm]{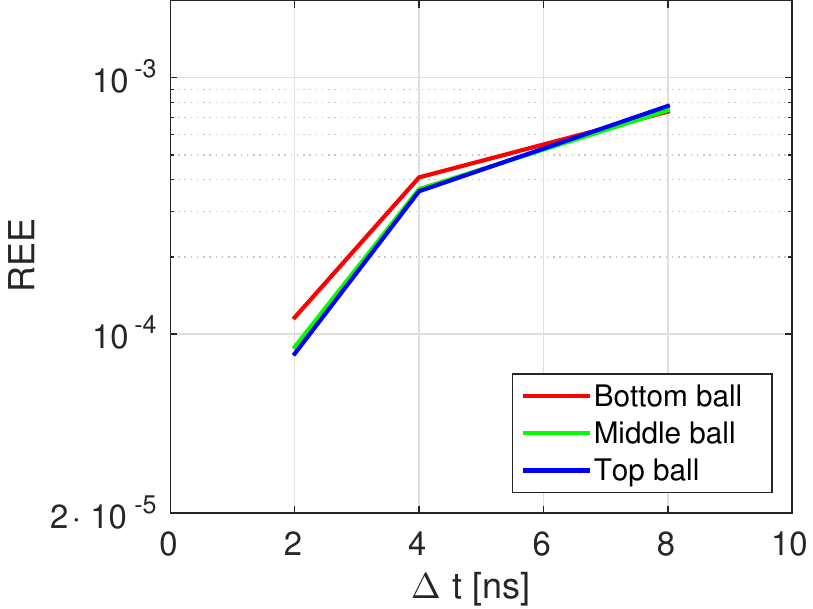}\label{sim:fig:error_HGvsKWAVE_homo}}
\hspace*{0.5cm}
\subfloat[][$\ell_\infty$ error.]{\includegraphics[width = 6cm]{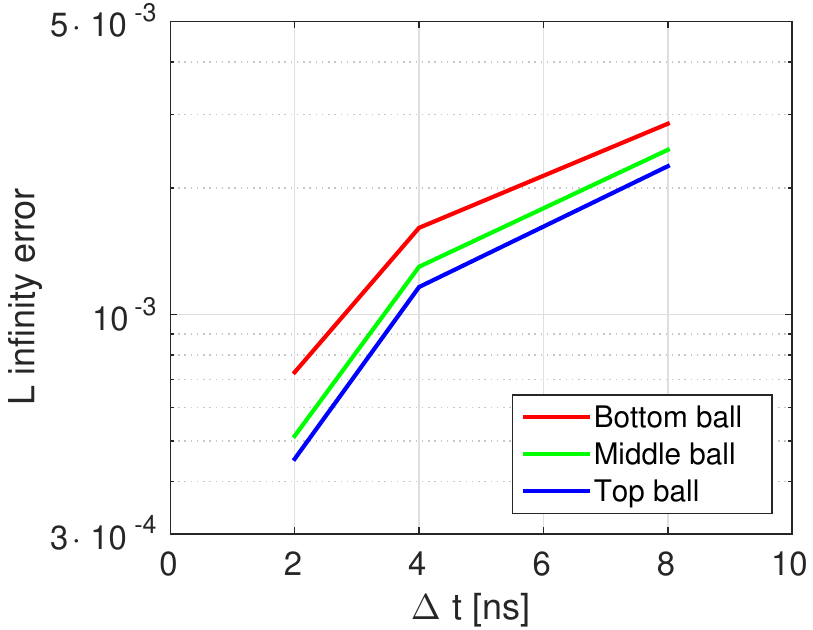}\label{sim:fig:error_HGvsKWAVE_homo_max}}
\caption{Homogeneous sound speed. Convergence of the forward HG solver with decreasing time step $\Delta t \in \{8, 4, 2\}$ [ns] in terms of (a) REE and (b) $\ell_\infty$ error.}
\end{center}
\end{figure}

%The sound speed in the domain is homogeneous and equal to $c_0 = 1,500$ [m/s]. In figure \ref{sim:fig:InitialPressure} we plot the initial pressure. 
We place a single sensor at $x_0 = (6.35, 0)$ [mm] and shoot 4,000 rays in directions equisampling the angle range $[0,\pi]$. %$[\phi_{\min} = 0$,\phi_{\max} = \pi]$.
A subset of the trajectories obtained is shown in figure \ref{sim:fig:Rays_homo}.
These trajectories are straight lines, as expected in the homogeneous case. 
We denote with $g^{kW}_{i} = \Pat_{kW} (u_i)$ and $g^{HG}_{i} = \Pat_{HG} (u_i)$ the \texttt{k-Wave} and HG solutions of the forward problem for the initial pressure $u_i$. We evaluated both methods for time step lengths $\Delta t \in \{2, 4, 8\}$ [ns].  
%We compute the pressure time series $g_{HG, i}, g_{kW, i}$ at sensor $S$ using time steps $\Delta t \in \{2, 4, 8\}$ [ns] both with the HG solver and \texttt{k-Wave} and compare the results between the two. 
In figure \ref{sim:fig:HGvsKWAVE_homo} we plot the solutions at $x_0$ 
for each of the initial conditions individually. The HG solution (thick line) is superposed over \texttt{k-Wave} solution (thin line), the time step length is $\Delta t = 8$ [ns].
We colour coded the initial conditions: bottom (the ball closest to the sensor) in \textcolor{red}{red}, middle in \textcolor{green}{green} and top in \textcolor{blue}{blue}.
Even for this coarse step, both solutions differ by at most 5\% (see figure \ref{sim:fig:main_HGvsKWAVE_homo}(b)).
All pressure time series were normalised by $\ell_\infty$ norm of the \texttt{k-Wave} solution for the bottom ball, $\|g^{kW}_{1}\|_{\infty} = 0.062$.

Figures \ref{sim:fig:main_HGvsKWAVE_homo}(b-d) show the difference between the \texttt{k-Wave} and HG solutions for decreasing time steps $\Delta t \in \{8, 4, 2\}$ [ns].
%We observe a decreasing trend as the step is reduced, with error values in the order of 1\% for $\Delta t = 2$ [ns]. 
We observe that halving the step essentially halves the $\ell_\infty$ error. 
%These errors are consistent across the three initial pressures, which are propagated independently. 
This is in agreement with HG solver for homogeneous sound speed being \textit{exact} up to the numerical errors due to interpolation between the grids, numerical integration and the approximation of the derivative of the Green's function. 
%The use of a small grid spacing (with respect to the sound speed) allows us to show the convergence of our solution in this numerical example.
Figure \ref{sim:fig:error_HGvsKWAVE_homo} shows the convergence of the  \textit{relative energy error} 
\begin{equation}\label{eq:REE}
REE^{\textnormal{forward}}(g^{HG}) = \frac{\sum_{t_i \in (0, T)}(g^{kW}(t_i) - g^{HG}(t_i))^2}{\sum_{t_i \in (0, T)}(g^{kW}(t_i))^2},
\end{equation}
for each initial pressure with the decreasing time step length.
%where $g^{kW}$ is the forward pressure time series computed with \texttt{k-Wave} and $g^{HG}$ the one obtained with the HG solver.
%The sum is calculated over all time steps $t_i$ and normalised to make the errors for different initial pressures comparable.
The corresponding non-normalised $\ell_\infty$ errors,
\begin{equation}\label{eq:sim:linf}
\|g^{kW} - g^{HG}\|_\infty,
\end{equation}
are shown in figure \ref{sim:fig:error_HGvsKWAVE_homo_max}. Both metrics behave similarly, thus henceforth we use REE with its convenient normalisation. 
%We observe a similar behaviour in both error plots, REE and $\ell_\infty$.
%However, since the $\ell_\infty$ error is not normalised by the energy of the initial pressure, it is sensitive to its location within the domain.
%For the remainder of this chapter we use the REE as our reference metric, as it allows us to compare errors for initial pressures in different locations.
The linear convergence is due to convergence of the grid interpolation error as the RK scheme is of higher order and the Green's function derivative $G'$ approximation is performed on a finer grid. Beyond time step length of $\Delta t$ = 2 [ns], the error due to spatial interpolation of $u_i$ with fixed grid spacing $\Delta x$ and accuracy of $G'$ approximation become significant, and we cannot arbitrarily reduce the error solely by further time step refinement without also refining the spacial grid for $u_0$ and the $G'$ approximation.

%In these plots we observe the reduction in the interpolation error between the grids, because the magnitude of all other sources of error is fixed for different time steps.
%For time step $\Delta t$ = 2 [ns], the magnitude of the other sources of error (Green's function approximation, RK solution) is already comparable, and hence we cannot arbitrarily reduce the error by taking finer time steps.

%====================================================================================================
% ERROR ANALYSIS FOR ADJOINT PROBLEM
%====================================================================================================

\subsection{Error Analysis of the Adjoint Problem}\label{sec:sim:adj}

\begin{figure}
\begin{center}
\subfloat[][]{\includegraphics[height = 7.4cm]{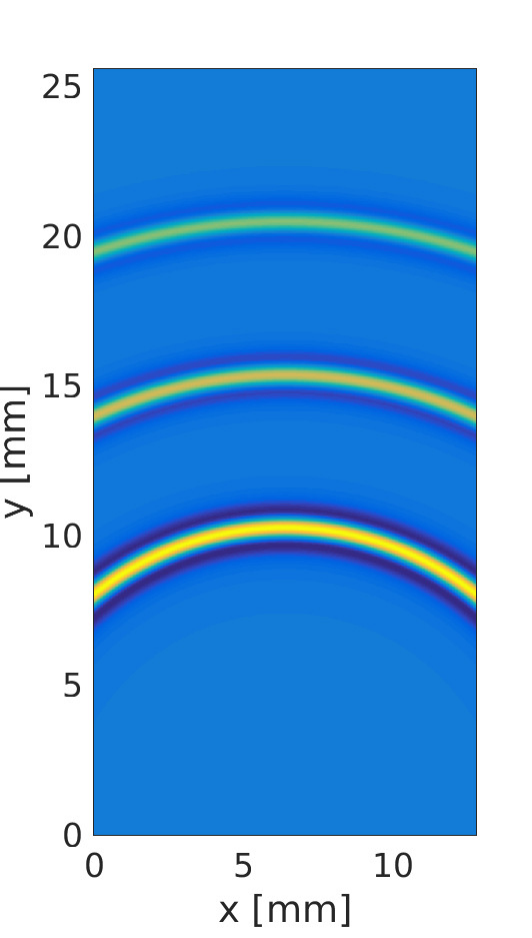}}
\subfloat[][]{\includegraphics[height = 7.5cm]{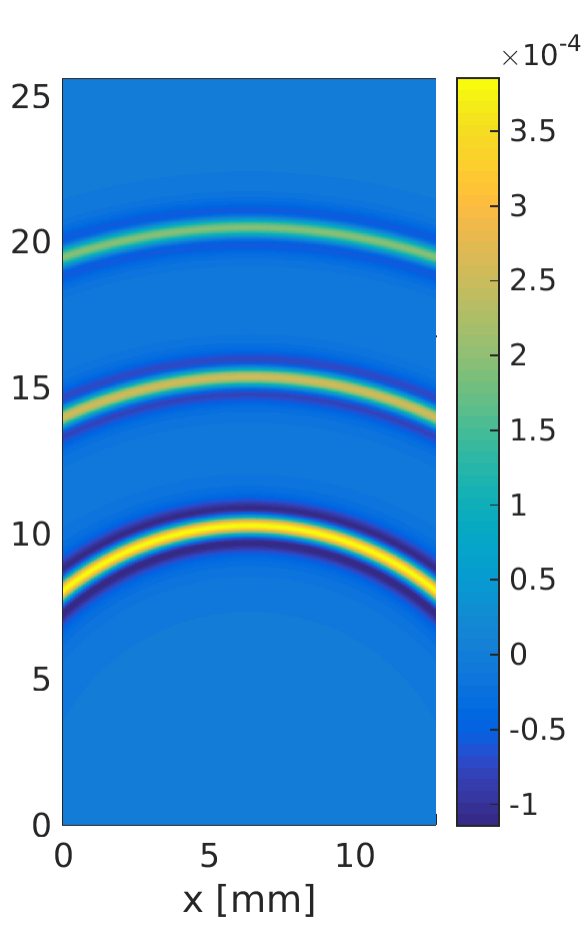}}
\subfloat[][]{\includegraphics[height = 7.5cm]{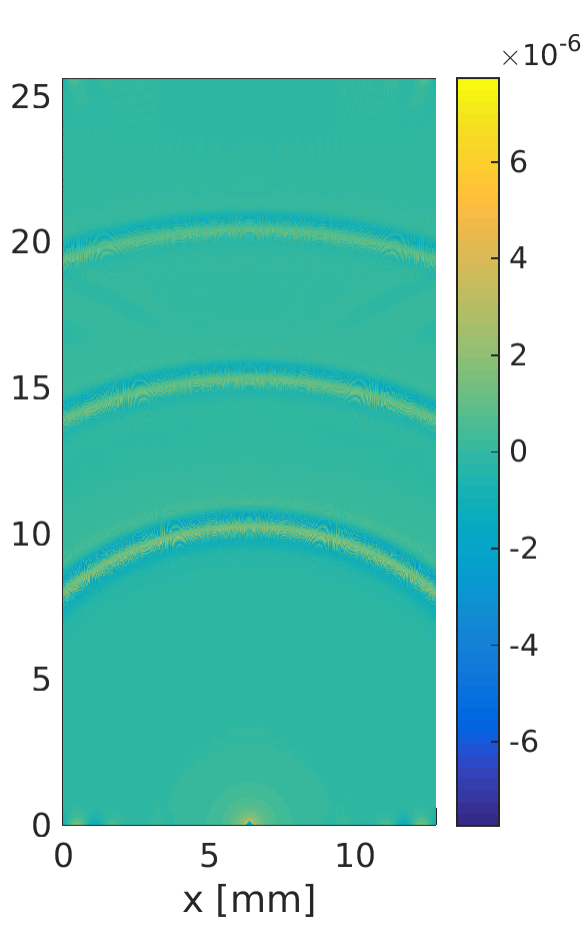}\label{sim:fig:adjoint_2D_homo_error}}
\caption{Homogeneous sound speed. Time step $\Delta t = 2$ [ns]. (a) Sum of adjoint \texttt{k-Wave} solutions $\sum_j \Pat_{kW}^\star(g^{kW}_j)$. (b) Sum of adjoint HG solutions $\sum_j \Pat_{HG}^\star(g^{kW}_j)$. (c) Difference between \texttt{k-Wave} and HG. }
\label{sim:fig:adjoint_2D_homo}
\end{center}
\end{figure}

%In this subsection, we analyse the discrepancies between HG and \texttt{k-Wave} for the adjoint problem. To this aim, we use the domain characterised in the previous subsection.

%==================== HOMOGENEOUS DOMAIN
We solve the adjoint PAT problem for the simulated data $g^{kW}_j,\; j\in \{1, 2, 3\}$ obtained with \texttt{k-Wave} in the previous section \ref{sec:sim:fwd}.  %this is, we propagate back the time series from the three solved forward problems.
The \texttt{k-Wave} solution is assumed to have no errors, which allows us to explore the error in our adjoint solver independently from the forward.
Because of the reversibility of rays, we can reuse the ray trajectories obtained for the forward problem.
In figure \ref{sim:fig:adjoint_2D_homo} we plot the sum of the three adjoint solutions for $g^{kW}_j$ obtained with \texttt{k-Wave} (a),  HG (b) and their difference (c) for time step $\Delta t = 2$ [ns]. The amplitude profile along the central cross-section $x = 6.35$ [mm] is shown in Figure \ref{sim:fig:adjoint_2D_homo_error_line}. 
%Figure \ref{sim:fig:adjoint_2D_homo_error} contains the difference between the \texttt{k-Wave} and HG adjoints for $\Delta t = 2$ [ns].
Similarly, as for the forward problem, the $\ell_\infty$ error of the adjoint solution is in the order of 1\%.
%Similarly as for the forward problem, the error in the amplitude for this particular time step are in the order of 1\%.

The convergence of the \textit{relative energy error} for the adjoint problem for $g = g^{kW}_j, \; j=\{1,2,3\}$
\begin{equation}
REE^{\textnormal{adjoint}}(g) = \frac{\sum_{x_i \in \Omega}((\Pat^{*}_{kW}\,g)(x_i) - (\Pat^{*}_{HG}\,g)(x_i))^2}{\sum_{x_i \in \Omega}((\Pat^{*}_{kW}\,g)(x_i))^2}
\end{equation}
%where $(\Pat^{*}_{kW}\,g)(x_i)$ and $(\Pat^{*}_{kW}\,g)(x_i)$ are the solution to the problem for \texttt{k-Wave} and HG, respectively, the input $g$ is set to $g^{kW}_j$ and the sum is computed over all pixels $x_i$ in the domain $\Omega$.
with the decreasing time step is depicted in Figure \ref{sim:fig:adjoint_2D_homo_error_3balls}
evidencing that HG adjoint solver is also \textit{exact} for homogeneous sound speeds up to the same numerical errors as discussed for the forward problem.
%We observe that the convergence is consistent for the bottom, middle and top balls, with very similar magnitudes for the three time steps.
%Further reductions in the step size may not lead to an arbitrary decrease in the discrepancy between the solvers, as we reach the nearest-neighbour interpolation error and the error in the approximation of the Green's function.

\begin{figure}
\begin{center}
\subfloat[][]{\includegraphics[width = 6.8cm]{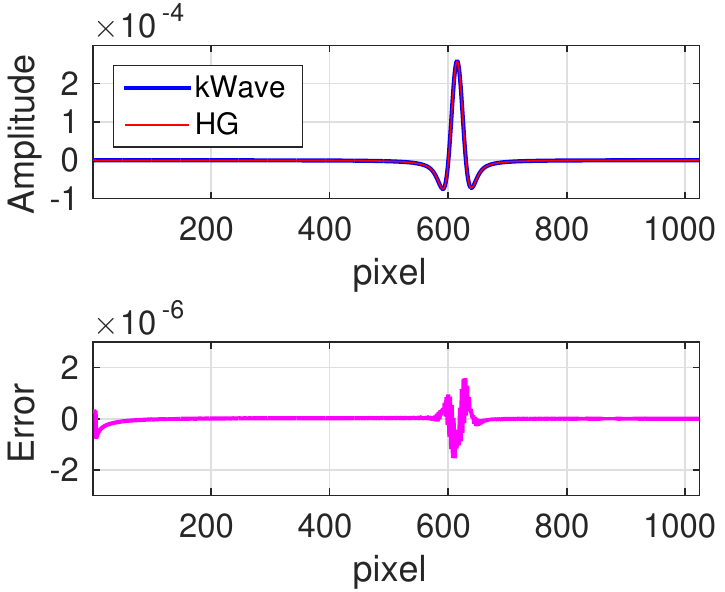}\label{sim:fig:adjoint_2D_homo_error_line}}
\subfloat[][]{\includegraphics[width = 6.8cm]{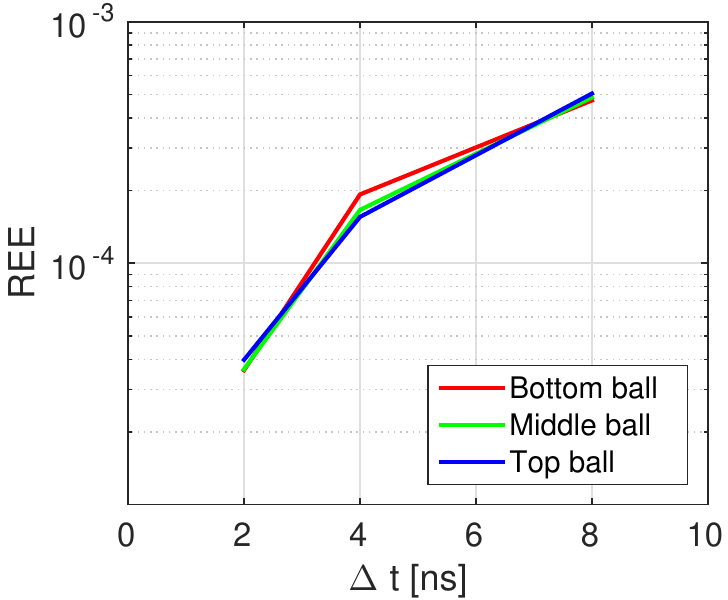}\label{sim:fig:adjoint_2D_homo_error_3balls}}\hspace*{0.2cm}
\caption{Homogeneous sound speed. (a) Adjoint solution profile through $x = 6.35$ [mm] for \texttt{k-Wave} and HG with $\Delta t = 2$ [ns]. (b) Convergence of adjoint HG solver with decreasing time step $\Delta t = \{8,4,2\}$ [ns].}
\label{sim:fig:adjoint_2D_homo_error}
\end{center}
\end{figure}
}

%====================================================================================================
% ACOUSTICALLY HETEROGENEOUS DOMAIN
%====================================================================================================
\subsection{Simulated 2D Vessel Phantom with Heterogeneous Sound Speed}

We complete our simulations with an experiment using a more realistic vessel phantom visualised in \ref{sim:fig:acoustically_homo_domain_nonsmooth}.
The domain $\Omega$ was discretised with  $128\times 256$ equispaced grid points with $\Delta x = 100$ [$\mu$m].
The sound speed $c(x)$ shown in figure \ref{sim:fig:soundSpeed_2D} was generated with the \texttt{peaks()} function in \texttt{MATLAB} and $c(x) \in 1,500 \pm 10\%$ [m/s].
The time step was set to $\Delta t = 10$ [ns] for both \texttt{k-Wave} and HG simulations. To avoid Gibbs phenomena in \texttt{k-Wave}, the initial pressure is a priori smoothed, see figure \ref{sim:fig:acoustically_homo_domain_smooth}.

%The initial pressure is shown in figure \ref{sim:fig:acoustically_homo_domain_nonsmooth}.This pressure tries to emulate the structure of vessels, which is the typical object what we may encounter in a real PAT scenario.
%Due to the limitations of \texttt{k-Wave}, it is not possible to use the initial pressure as it is since it contains high spatial frequencies. 
%To propagate them in time, we would need to use a very fine time step. Given a time step, the \texttt{k-Wave} toolbox provides a smoothing function that applies a smooth filter to cut off all the frequencies higher than those allowed by the specified step. We choose $\Delta t = 10$ [ns] and apply the filter to obtain the initial pressure in \ref{sim:fig:acoustically_homo_domain_smooth}.

% Initial Pressure
\begin{figure}
\begin{center}
\subfloat[][]{\includegraphics[height = 7.5cm]{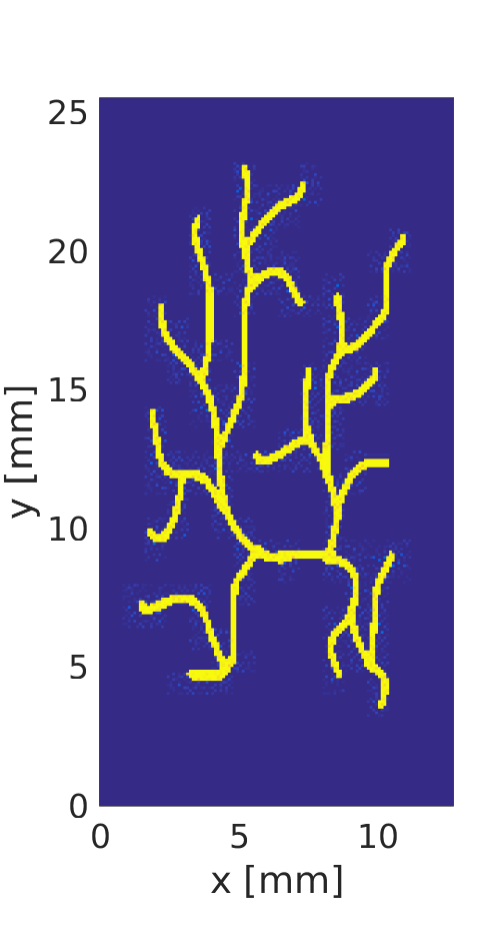}\label{sim:fig:acoustically_homo_domain_nonsmooth}}
\subfloat[][]{\includegraphics[height = 7.45cm]{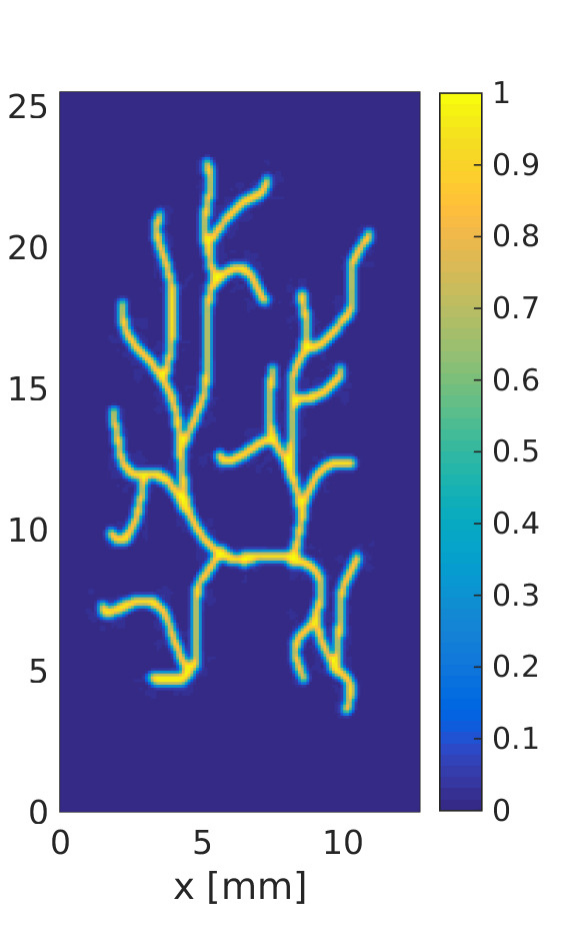}\label{sim:fig:acoustically_homo_domain_smooth}}
\subfloat[][]{\raisebox{6pt}{\includegraphics[height = 7.1cm]{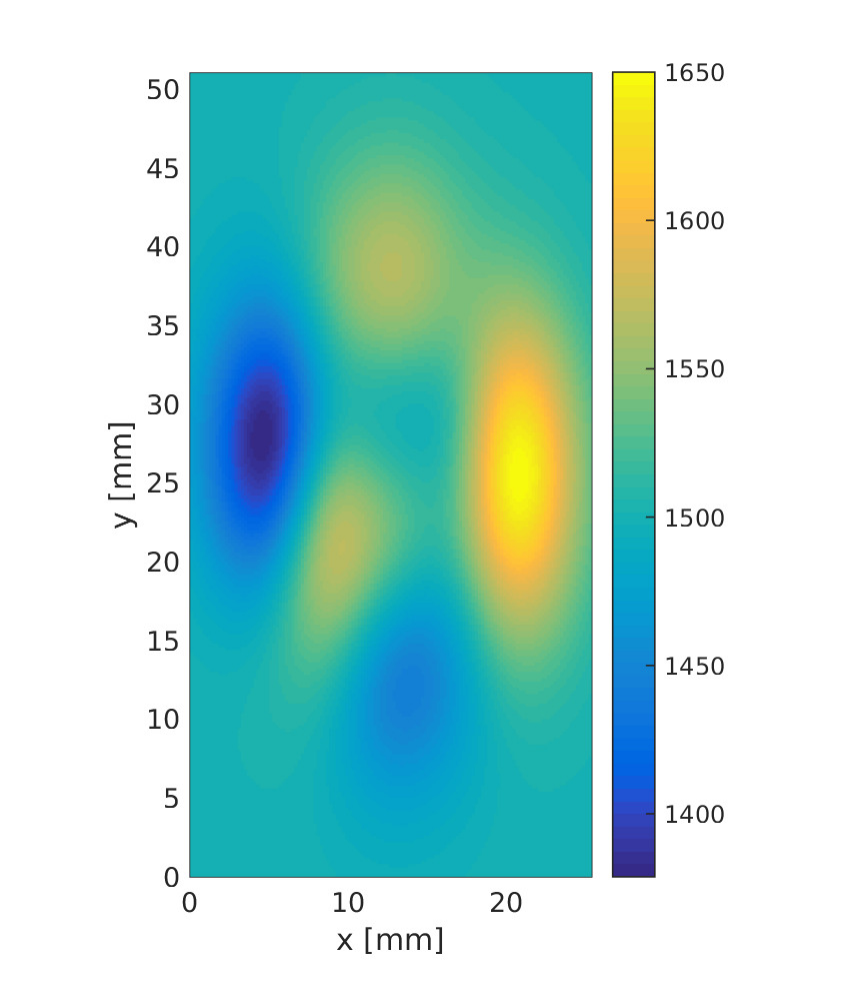}\label{sim:fig:soundSpeed_2D}}}
\caption{(a) Initial pressure (without smoothing). (b) Initial pressure $u_0(x)$ (\texttt{k-Wave} smoothing). (c) Sound speed $c(x)$.}
\label{sim:fig:acoustically_homo_domain}
\end{center}
\end{figure}

% Trajectories 
\begin{figure}
\begin{center}
\vspace*{1cm}
\subfloat[][]{\includegraphics[height = 6.5cm]{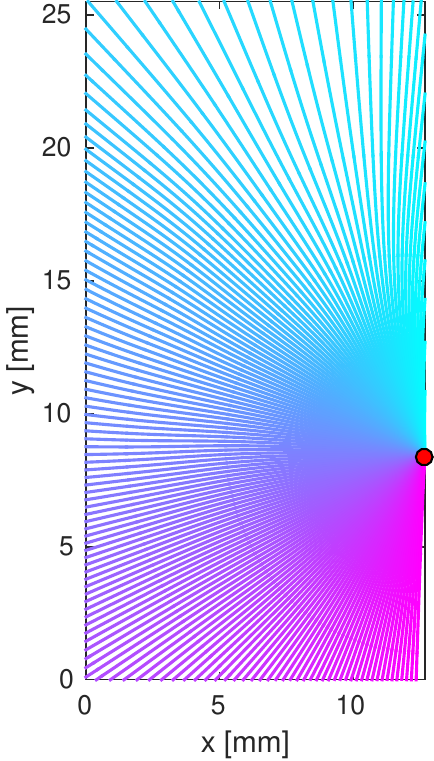}\label{sim:fig:acoustic_het_2D_trajec_1}}\hspace*{0.2cm}
\subfloat[][]{\includegraphics[height = 6.5cm]{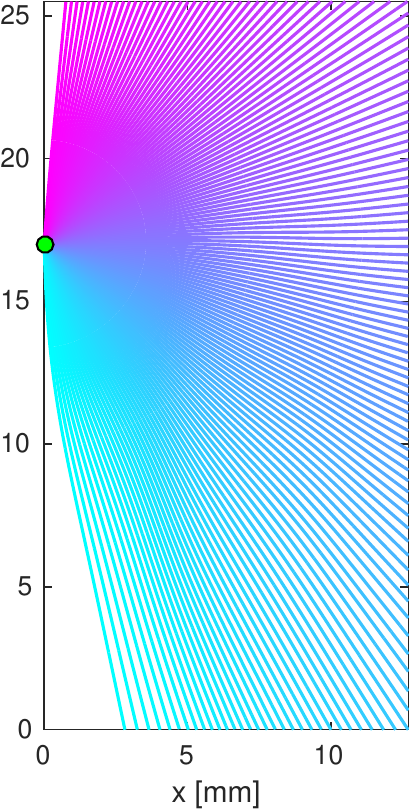}\label{sim:fig:acoustic_het_2D_trajec_2}}\hspace*{0.2cm}
\subfloat[][]{\includegraphics[height = 6.5cm]{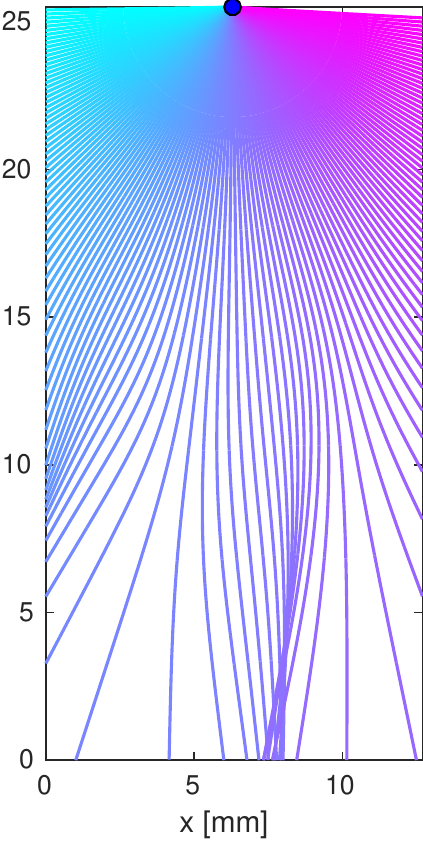}\label{sim:fig:acoustic_het_2D_trajec_3}}
\caption{Ray trajectories from 3 randomly chosen sensors.}
\label{sim:fig:acoustic_het_2D_trajec}
\end{center}
\end{figure}

\paragraph{Forward problem} 
We place a sensors at each grid point on the boundary of $\Omega$ amounting to 764 sensor in total.
From each sensor location we shoot 4,000 rays with direction angles uniformly distributed between 0 and 2$\pi$ to ensure that we cover the entire domain.
In figure \ref{sim:fig:acoustic_het_2D_trajec} we exhibit a subset of trajectories for three sensors chosen at the random. %in the homogeneous sound speed example.
Due to heterogeneity of $c$, these trajectories are no longer straight lines and present the two problematic scenarios: caustics (figure \ref{sim:fig:acoustic_het_2D_trajec_3}) and 
shadow regions (figure \ref{sim:fig:acoustic_het_2D_trajec}(b-c)). 
\rB{We note that the pseudo-spectral finite difference method in \texttt{k-Wave} is a full wave solver and as such not affected neither by caustics nor shadow regions.}
In figure \ref{sim:fig:sinogram_het_2D} we plot the sinograms for \texttt{k-Wave}, HG and their difference.
The sinograms were normalised to 1 using the $\ell_\infty$ norm of the \texttt{k-Wave} solution, $\|g^{kW}\|_\infty= 0.9329$.
The $\ell_\infty$ error of the HG sinogram is significantly larger than in the homogeneous case (albeit $\ell_\infty$ appears a pessimistic metric in this case) while its REE is 1.99\%.

% Sinogram
\begin{figure}
\begin{center}
\subfloat[][]{\includegraphics[height = 5cm]{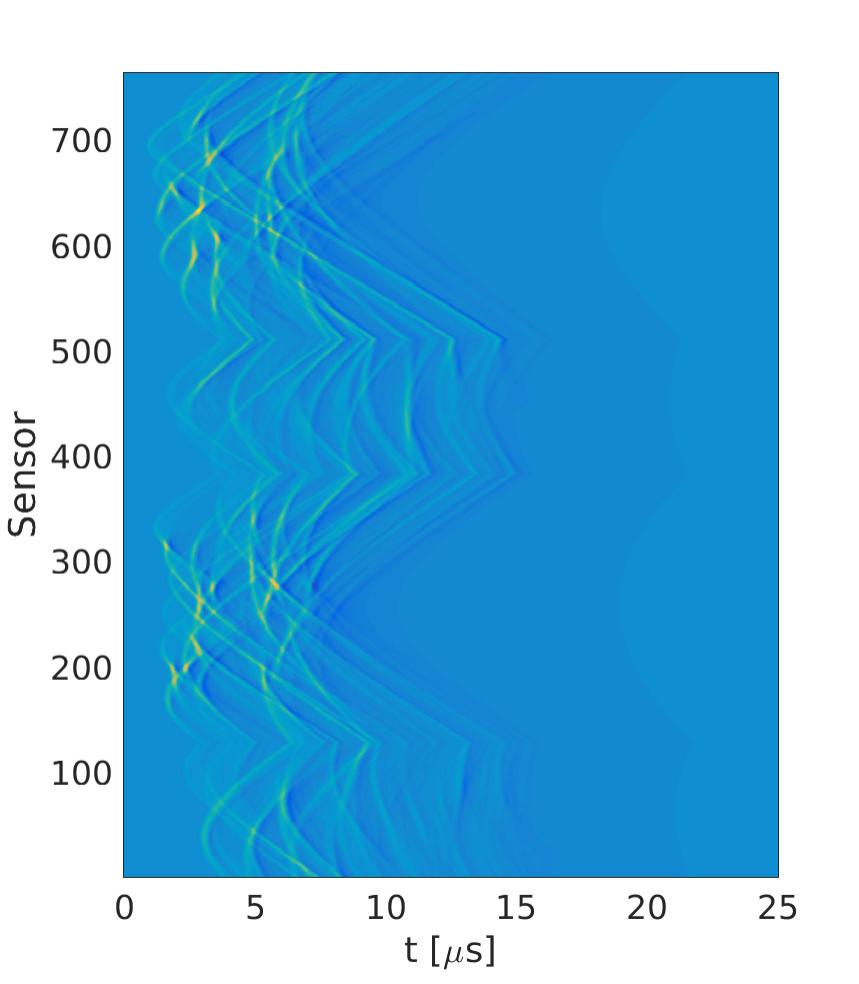}}
\subfloat[][]{\includegraphics[height = 5cm]{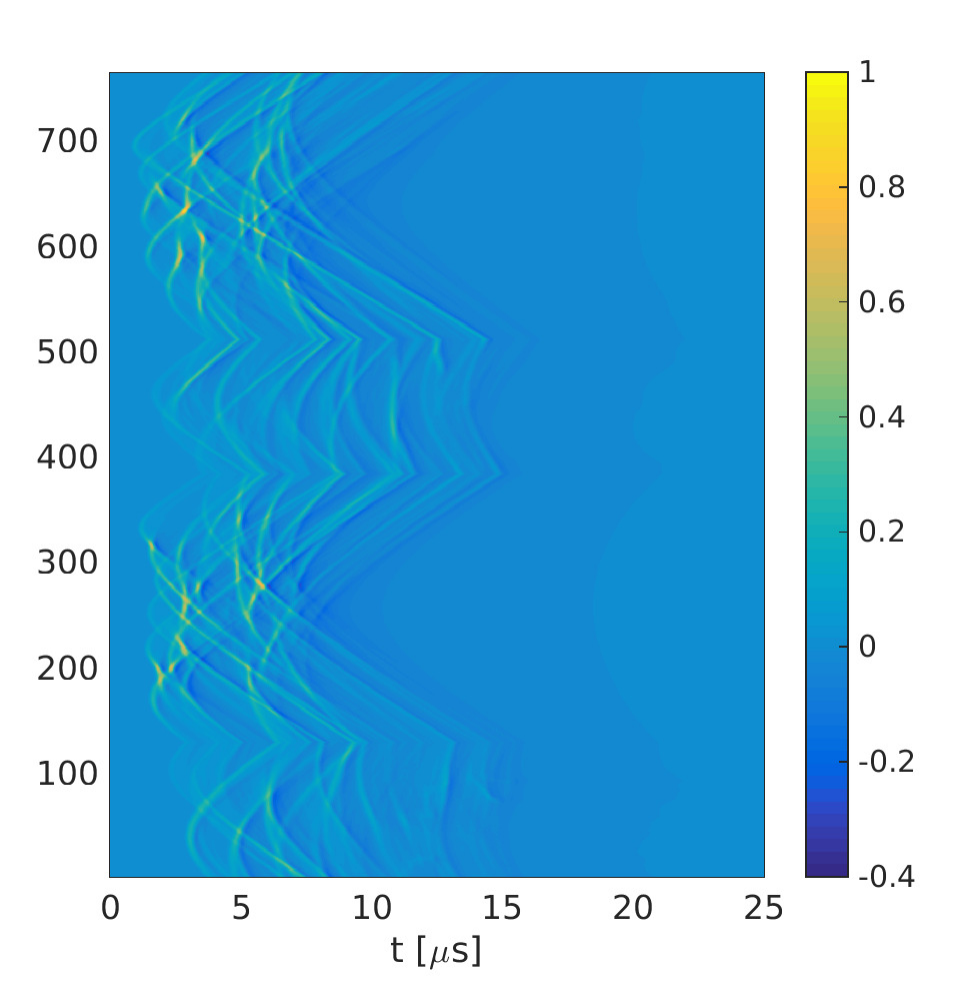}}
\subfloat[][]{\includegraphics[height = 5cm]{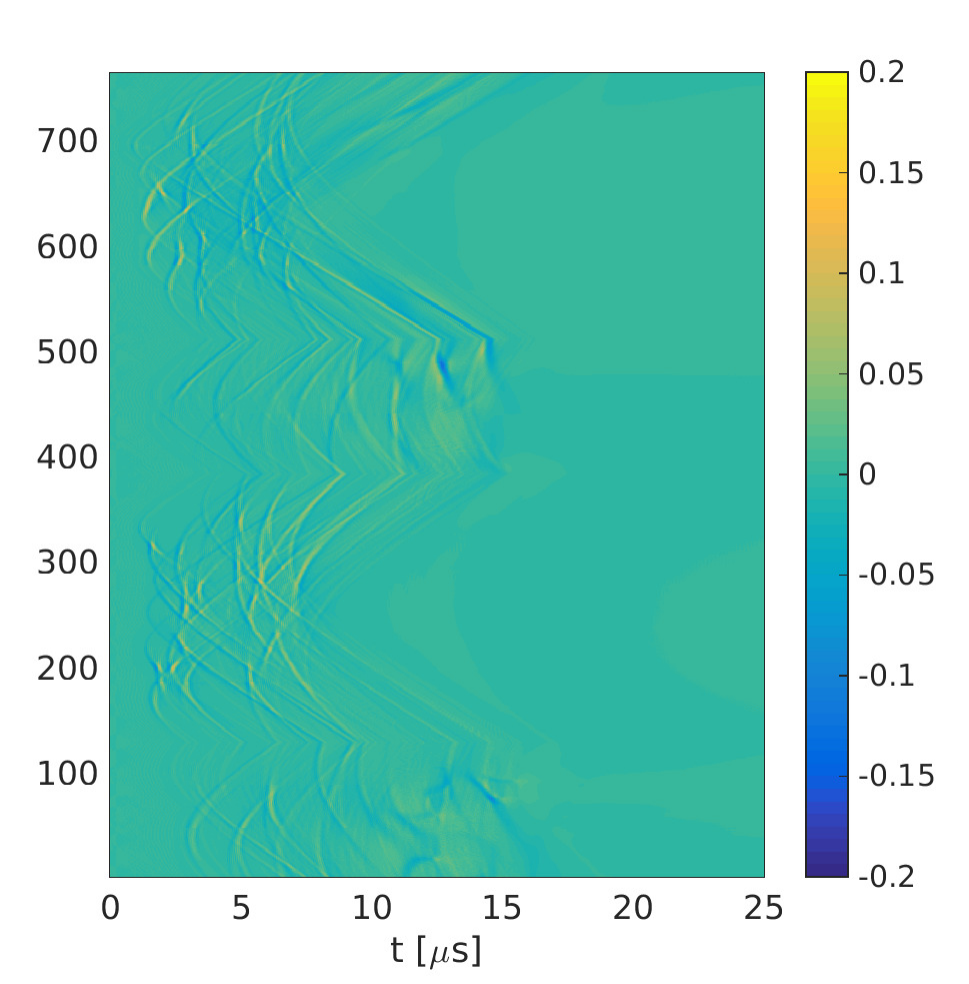}}
\caption{(a) \texttt{k-Wave} sinogram. (b) HG sinogram. (c) Difference between \texttt{k-Wave} and HG sinograms.}
\label{sim:fig:sinogram_het_2D}
\end{center}
\end{figure}

\paragraph{Consistent forward-adjoint}
We compute the \texttt{k-Wave} adjoint with the \texttt{k-Wave} forward data and the HG adjoint with its corresponding HG forward data.
The results for these simulations are shown in figures \ref{sim:fig:adjoint_het_2D}(a-b).
The overall structure of the adjoints is the same, containing similar artefacts for both methods.
However, we observe more significant differences in the bottom left part of the image (c.f.~figure \ref{sim:fig:adjoint_het_2D_error}).
The $\ell_\infty$ norm of the HG adjoint is on the order of 20\%, although the large errors are concentrated in small region of the image. Here, the REE is equal to 2.23\%.

% Adjoint
\begin{figure}
\begin{center}
\subfloat[][]{\includegraphics[height = 7.5cm]{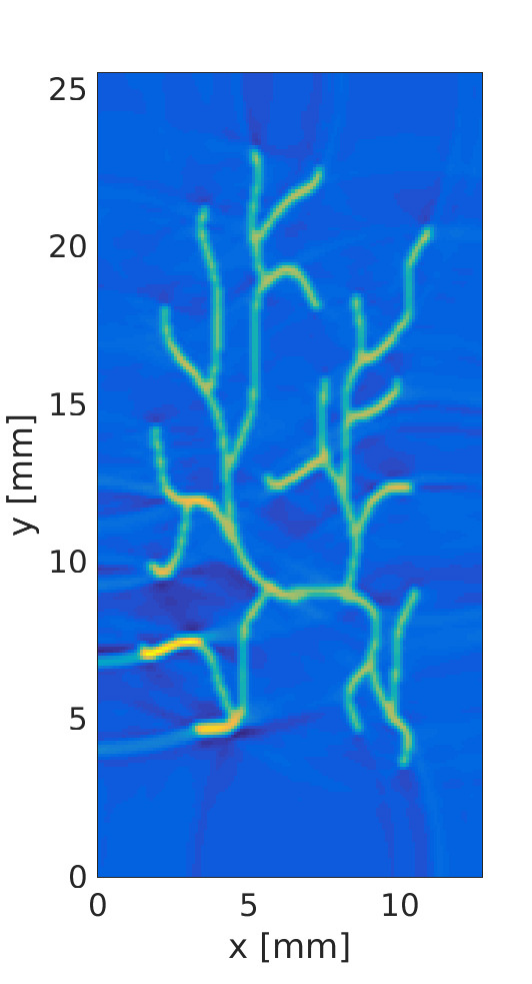}}
\subfloat[][]{\includegraphics[height = 7.5cm]{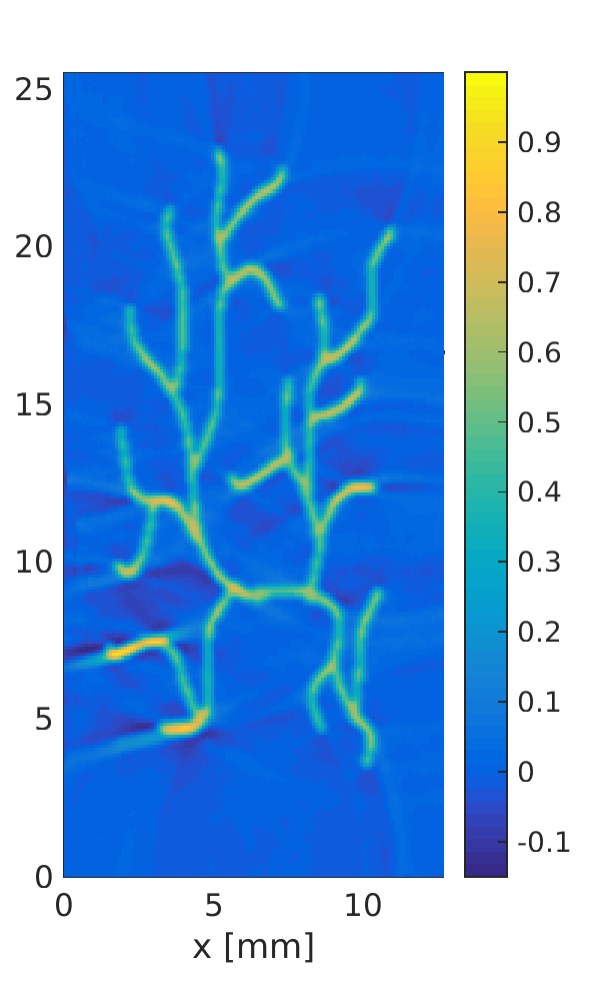}}
\subfloat[][]{\includegraphics[height = 7.5cm]{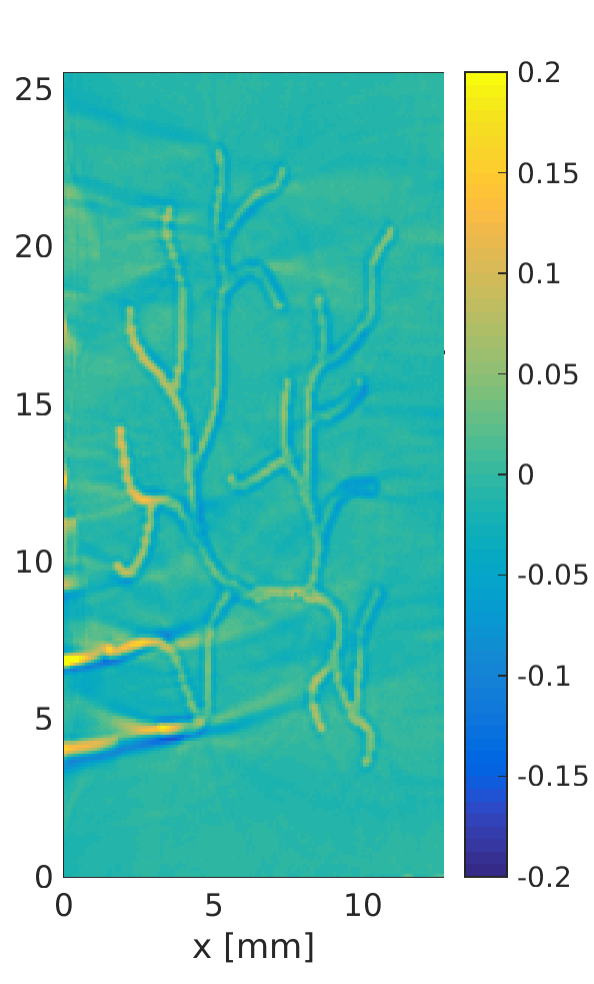}\label{sim:fig:adjoint_het_2D_error}}\\
\subfloat[][]{\includegraphics[height = 7.5cm]{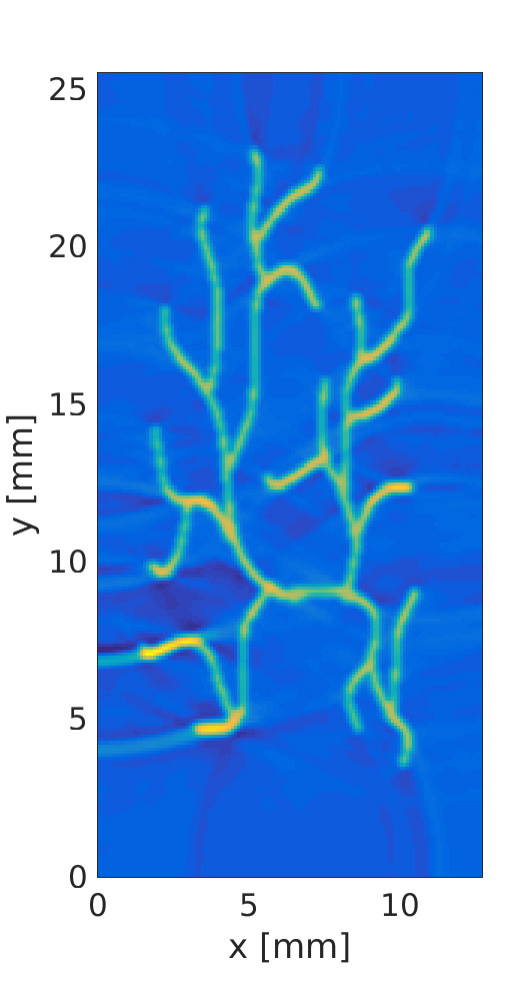}\label{sim:fig:adjoint_het_2D_mixKW}}
\subfloat[][]{\includegraphics[height = 7.5cm]{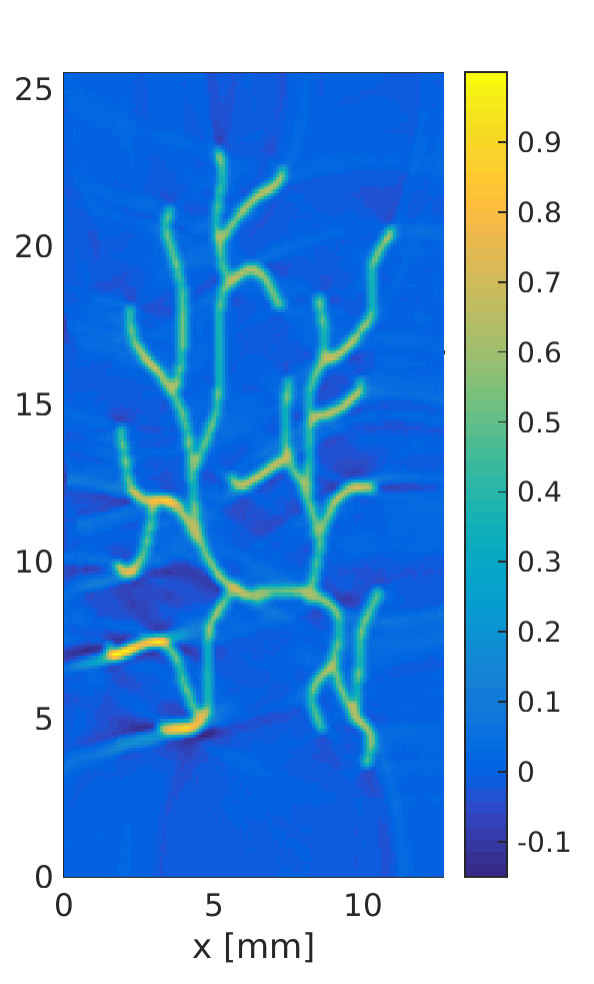}\label{sim:fig:adjoint_het_2D_mixHG}}
\subfloat[][]{\includegraphics[height = 7.5cm]{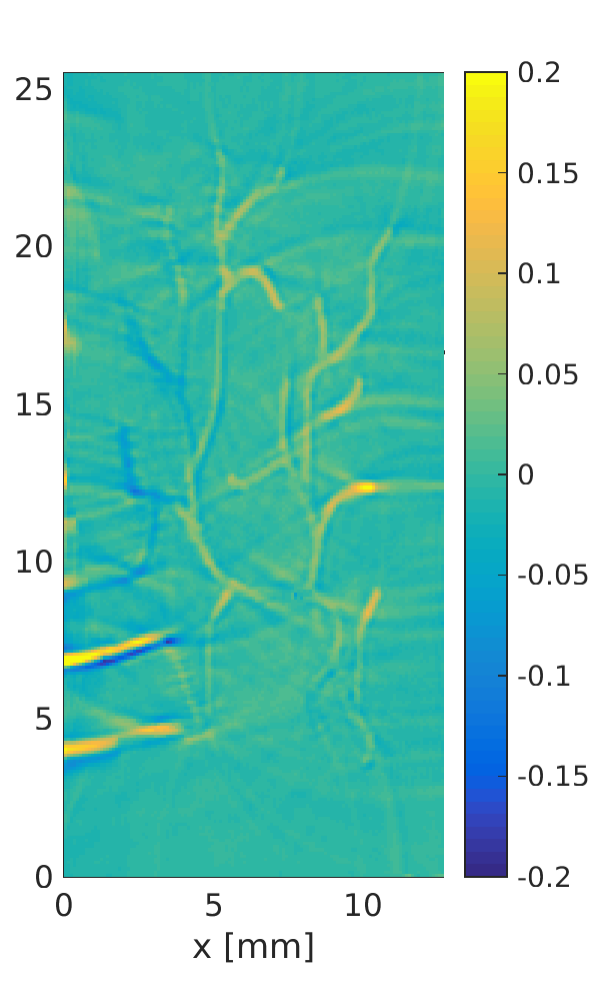}\label{sim:fig:adjoint_het_2D_mixError}}
\caption{(a) \texttt{k-Wave} adjoint (gold standard). (b) HG adjoint. (c) Difference between \texttt{k-Wave} and HG adjoints (same as error of HG adjoint).
         (d) \texttt{k-Wave} adjoint - HG data. (e) HG adjoint - \texttt{k-Wave} data. (f) Difference between mixed adjoints.}
\label{sim:fig:adjoint_het_2D}
\end{center}
\end{figure}

\paragraph{Mixed forward-adjoint}
We solve two additional adjoint problems combining HG forward data with \texttt{k-Wave} adjoint solver and vice-versa to isolate forward/adjoint errors.
Figure \ref{sim:fig:adjoint_het_2D_mixKW} shows the \texttt{k-Wave} solution of adjoint problem with HG forward data and
figure \ref{sim:fig:adjoint_het_2D_mixHG} the HG adjoint with the \texttt{k-Wave} forward data.
Both solutions have similar structure and artefacts also similar to the consistent case.
The difference between the two mixed adjoints is shown in figure \ref{sim:fig:adjoint_het_2D_mixError}, with $\ell_1$ norm of around 20\%.
The REE for the combination HG forward data / \texttt{k-Wave} adjoint is 0.74\% and for \texttt{k-Wave} forward data / HG adjoint is 1.79\%, 
both lower than the REE for the consistent forward-adjoint, as one would expect.
}
%\input{52_Simulations_heterogeneous}

%==============================
% 5. Conclusions
%==============================
%==============================
% Conclusions and Future Work
%==============================
\section{Conclusions and outlook}\label{sec:conclusions}
In this paper we presented a novel acoustic Hamilton-Green solver for the forward and adjoint problems in PAT. The name is derived from the two main ingredients of the solver: the Hamiltonian system and the Green's integral formula for the solution of the wave equation with initial data and time varying sources. The solver can be viewed as a ray tracing approximation of the heterogeneous Green's function (more precisely its time derivative) in the Green's integral formula evaluated on the ray induced grid. 
The main benefit of this approach is that the solution is computed for each sensor independently, opening to PAT the realm of ray based solvers which are well established for other tomographic modalities e.g.~X-ray CT or PET. 
Furthermore, ray based solvers allow the flexibility to efficiently deal with
sparse sensors and regions of interest thus potentially enabling online/real-time or dynamic reconstructions.

We analysed the error the proposed Hamilton-Green method, discussed its convergence in the homogeneous case and evaluated its accuracy in heterogeneous case on a 2D numerical vessel phantom against an established finite difference pseudo-spectral method implemented in \texttt{k-Wave} toolbox \cite{treeby2010k}. 

There are many open questions that need addressing to advance understanding and performance of the Hamilton-Green solver.
A theoretical challenge is to find solution that holds valid in the presence of caustics without compromising the computational efficiency. In this context we are currently investigating two promising approaches: the dynamical ray tracing and Gaussian beams which both can be integrated with the Hamilton-Green solver. 
We are also working on an efficient implementation on GPU hardware which is imperative to apply the solver to large scale 3D real data problems. A particularly challenging implementation aspect is adaptive inserting and removing rays into and from the cone to control the impact of the changing ray density while maintaining the efficiency offered by one shot methods. 
Furthermore, we are working on coupling the Hamilton-Green solver with iterative methods for variational reconstruction which performance hinges on the ability to independently evaluate partial operators such as for instance the stochastic primal-dual hybrid gradient method recently proposed in \cite{ehrhardt2017faster}.

\appendix

\section{Implementation details}\label{app}

%\section{Discretization Aspects}
%===============================================
% Sound Speed on the rays
%===============================================

\subsection{Interpolation: from Cartesian to ray based grid}\label{app:intC2R}
When interpolating a function specified on Cartesian grid at a point in a ray based grid, $x(\ell, \theta; x_0)$, we use the Cartesian grid point closest to $x(\ell, \theta; x_0)$ in the $\ell_2$ sense. 

Quantities in the code interpolated from Cartesian to ray based grid are: slowness $\eta$ and its gradient $\nabla \eta$, initial pressure $u_0$. 

\subsection{Interpolation: from ray based to Cartesian grid}\label{app:intR2C}
When interpolating a function specified on ray based grid at a point in Cartesian grid, $x$, we use the the ray based grid point closest to $x$ in the $\ell_2$ sense.

Quantities in the code interpolated from ray based to Cartesian grid are: the solution of the adjoint problem \eqref{eq:GHAdj}, the reversed phase (propagation time), the reversed amplitude attenuation ($\mu_R$) used in the computation of the adjoint. To be precise, the latter can be obtained according to equations \eqref{eq:phaseRev} and \eqref{eq:muR} from the (non-reversed) phase and amplitude attenuation (or ray density), respectively.  

\subsection{Approximating the homogeneous Green's function at a point in the domain} \label{app:G0}
Equation \eqref{eq:GHAdj} requires knowledge of time derivative of the homogeneous Green's function $G_0$ for a sound speed value $c_0 = c(x(\ell,\theta; x_0))$. Obtaining $G_0$ for an exact value of $c_0 = c(x(\ell,\theta; x_0))$ (in particular if $c$ was interpolated on the rays using an interpolation formula of an order higher than the nearest neighbour) would undue increase the computational cost. Instead, we precompute the Green's function $G_0$ for a range of sound speeds $c_{\textrm{min}} \leq c(x) \leq c_{\textrm{max}}, \forall x\in\Omega$ at equal increments, $\Delta c$ fine enough and use the tabulated $G_0$ with the closest sounds speed value to approximate $G_0$ with $c_0 = c(x(\ell,\theta; x_0))$. The time derivative of $G_0$ is precomputed analogously.

\begin{figure}
\begin{center}
\vspace*{.1cm}
\hspace*{-0.3cm}\subfloat[]{\includegraphics[height = 5.8cm]{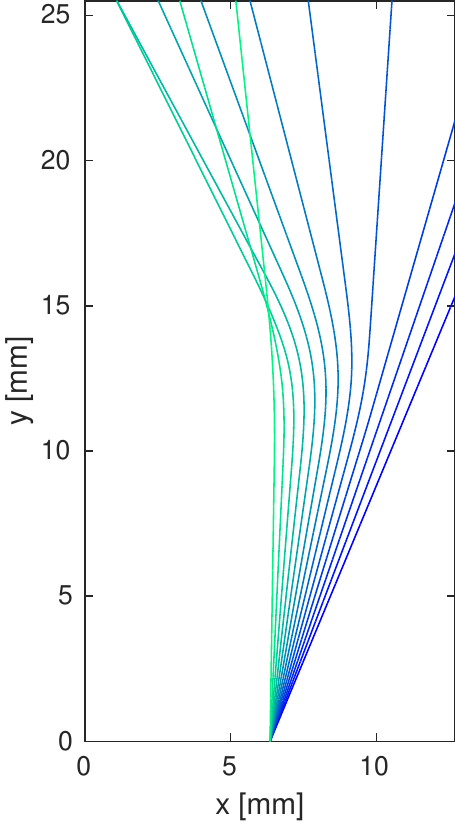}\label{sim:fig:caustic_subrays}}\hspace*{0.3cm}
\subfloat[]{\includegraphics[height = 6cm]{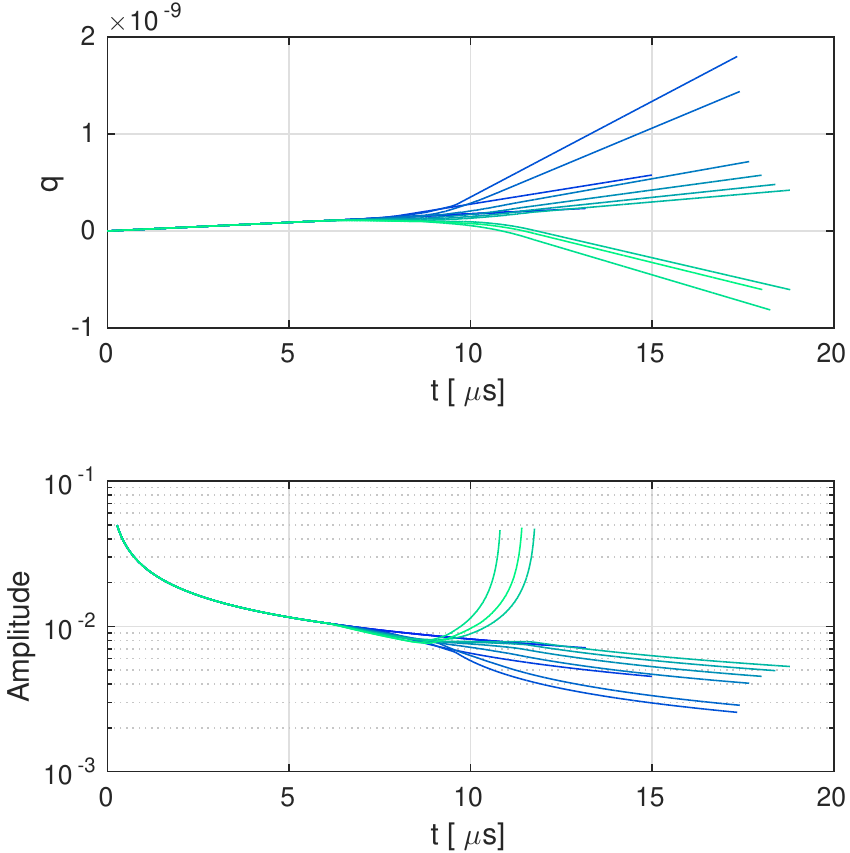}\label{sim:fig:caustic_QA}}
\caption{(a) Subset of ray trajectories shot from $x_0$ illustrating the transition between the rays which form the caustic (\textcolor{myGreen}{light green}) and those which bend around it (\textcolor{myBlue}{dark blue}). 
         (b) $q(t; x_0)$ (top) and $A(x(t; x_0))$ (bottom) corresponding to the subset of ray trajectories in (a).}
\label{}
\end{center}
\end{figure}

%=============== CAUSTICS - LENS =========
\rA{\section{Effect of caustics on the HG forward solution in an acoustic lens}\label{app:lens}}
{\color{myGreen}

\begin{figure}
\begin{center}
\captionsetup[subfigure]{justification=centering}
\subfloat[][Lens sound speed $c(x)$.]{\includegraphics[height = 7.5cm]{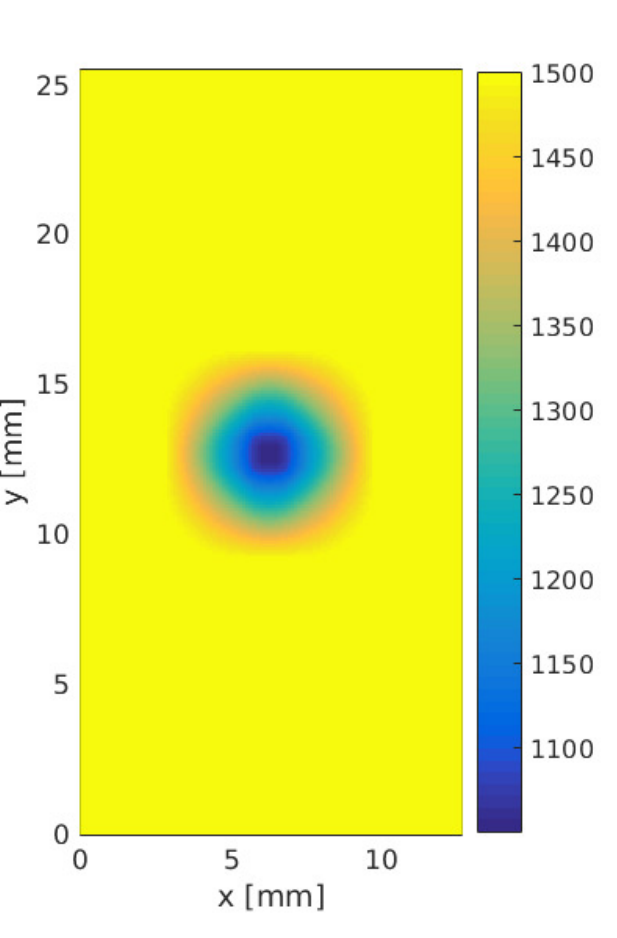}\label{sim:fig:soundSpeed_caustic}}\hspace*{0.2cm}
\subfloat[][Ray trajectories\\ from sensor $S$.]{\includegraphics[height = 7.5cm]{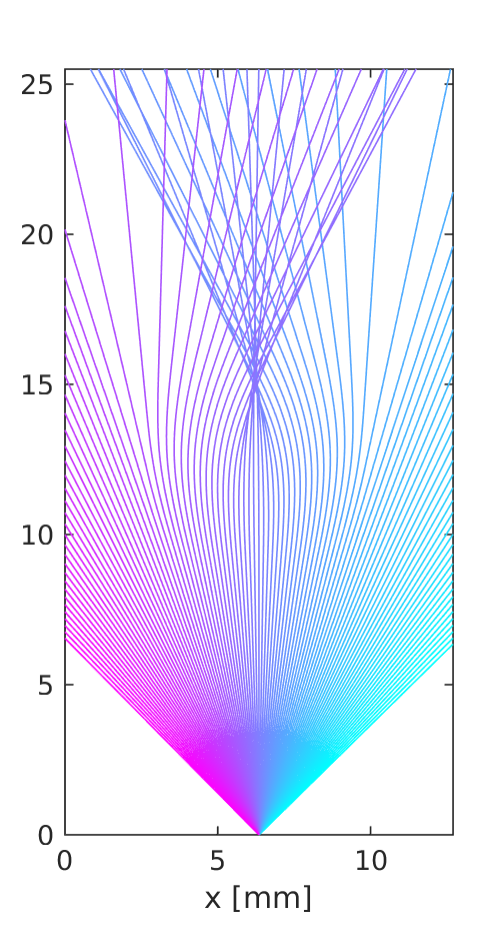}\label{sim:fig:soundSpeed_caustic_rays}}\hspace*{0.2cm}
\subfloat[][Superposition of initial \\ pressures $\sum_i u_{i}$.]{\includegraphics[height = 7.5cm]{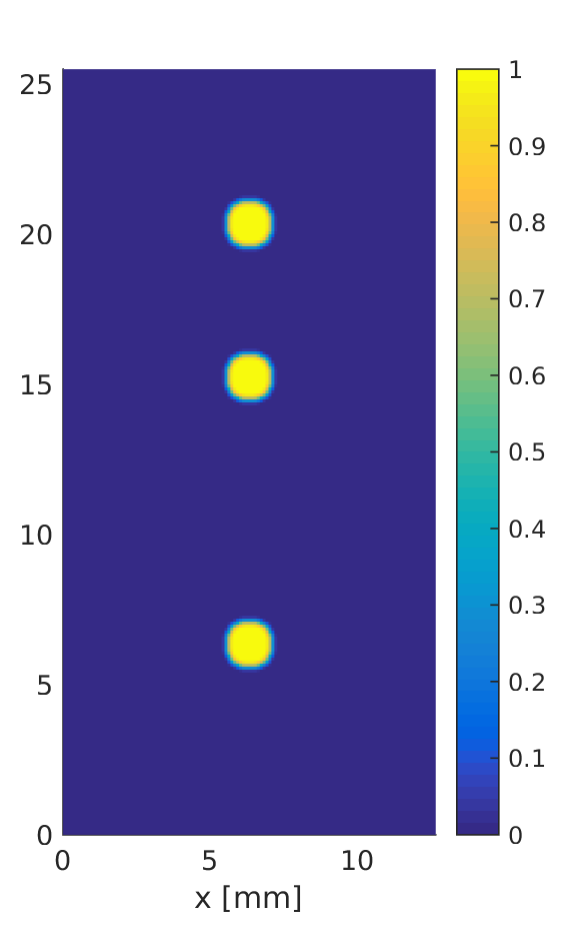} \label{sim:fig:soundSpeed_p0}}
\caption{(a) Acoustic lens sound speed $c(x)$.
         (b) Ray trajectories with a caustic in vicinity of $(6.35, 15)$ [mm]. 
         (c) Superposition of $\sum_i u_i$, each $u_i$ is an independent initial condition to the forward PAT problem.}
\label{}
\end{center}
\end{figure}
%========================== EXAMPLE
To illustrate the effect of caustics we compute the forward propagation through an acoustic lens. We would like to stress, that the construction relies on a large variation of the speed of sound (here some 33\%), hence such scenario is unrealistic to transpire in soft tissue where variations of sound speed are on the order of 5-10\%.  

Let $\Omega$ be a rectangular domain which we discretise with $N_x \times N_y = 128 \times 256$ equispaced grid points, $\Delta x$ = 100 [$\mu$m] with mostly homogeneous sound speed $c = 1,500$ [m/s] which dips in the centre to  $c_{\min} =1,050$ [m/s] forming an \textit{acoustic lens} depicted in figure \ref{sim:fig:soundSpeed_caustic}. 
%The sound speed is homogeneous over $\Omega$ with magnitude 1,500 [m/s] except for an \textit{acoustic lens} in the centre which is a circle with decreasing sound speed towards its centre, with minimum magnitude 1,050 [m/s], as shown in figure \ref{sim:fig:soundSpeed_caustic}. 
This lens has the effect of focusing the rays. %coming from any point outside it. 
%This example shows the behaviour of the solver in an extreme situation, with sound speed changes in the order of 50\%, where the effects of caustics are pronounced. Such a large variation in sound speed is not realistic in soft tissue, where sound speed changes are of the order $\pm 5$\%.

We place a single sensor $x_0$ at the point $(6.35, 0)$ [mm] on the boundary and compute the PAT data for this sensor for three initial pressures $u_i(x)$, $i\in\{b, m, t\}$, each $u_i(x)$, up to smoothing, assumes value 1 on a shifted ball of radius 8$\Delta x$. Figure \ref{sim:fig:soundSpeed_p0}, for compactness, shows their superposition $\sum_i u_i$. 

To evaluate the impact of the caustic on HG solver, we compare the respective solutions at the sensor to the  
%We evaluate the Hamilton-Green (HG) solutions against 
full wave solutions obtained with \texttt{k-Wave}. Both solvers use the same time step $\Delta t = 10$ [ns]. %(which we consider a gold standard acoustic solver for PAT, as explained in the manuscript). 

HG solver shoots 4,000 rays from the sensor $x_0$ evenly sampling the angular range $[\pi/4, 3\pi/4]$ to cover the relevant part of the domain. 
The ray cone is illustrated in figure \ref{sim:fig:soundSpeed_caustic_rays}. The inner trajectories focus in the vicinity of $(6.35, 15) [mm]$. These rays form the caustic while the outer trajectories bend around that point and hence are unaffected by the caustic.
The balls, $u_i$, were placed so to illustrate the differences in accuracy of the solution at the sensor in dependence of their location relative to the caustic:
%carried by rays before, at and after the caustic: 
the bottom ball $u_b(x)$ is seen before any focusing, the middle ball $u_m(x)$ overlaps with the focus point and the top ball $u_t(x)$ is seen after a subset of rays passed the caustic.

The solutions at the sensor for each initial condition $u_i(x)$ computed with HG and \texttt{k-Wave} are plotted in \ref{sim:fig:caustic_signals}, and their difference in 
\ref{sim:fig:caustic_error}.
The results in \ref{sim:fig:caustic_error} demonstrate that the solution for the initial pressure before the caustic is not affected (\textcolor{red}{bottom ball, $u_b$}). The error is largest for the initial pressure at the caustic (\textcolor{green}{middle ball, $u_b$}). This is mainly due to setting to 0 the ray amplitudes $A(x(t,x_0))$ for all times for which $|q(t; x_0)| < \epsilon$ with a small parameter $\epsilon = q(0; x_0)$ %enough $|q|< \epsilon$ 
(note that the pulse shape is essentially preserved).  
The initial pressure behind the caustic is also affected but to a lesser degree (\textcolor{blue}{top ball, $u_t$}). The amplitude contributions carried by the subset of rays that bends around the caustic are not affected at all, while the amplitude contribution of rays passing though the caustic undergo a phase shift (note the stretching of the pulse)
according to \eqref{eq:AMaslov}.
%Maslov showed that for the amplitude after the caustic it holds (see e.g.~\cite{runborg2007mathematical}, section 2.2)  
%\begin{equation} \label{req:amplitudeMaslov}
%A(x(t; x_0)) = A(x_0)\frac{\eta(x_0)}{\eta(x(t; x_0))}\sqrt{\frac{|q(0; x_0)|}{|q(t; %x_0)|}} e^{-i m \frac{\pi}{2}},
%\end{equation}
%where $m\in\matbb{N}$ is the Maslov's index. 
HG solution after the caustic does not account for the complex factor $\exp(-i m(t) \pi/2)$
but it compensates for the caustic induced reversing of the associated ray tube coordinate system by taking the absolute value of $q(t; x_0)$.

Figure \ref{sim:fig:caustic_subrays} shows the trajectories of a small subset of rays in the transition region between the inner rays which form the caustic in $\Omega$ (\textcolor{myGreen}{light green}) and the outer rays that bend around the caustic (\textcolor{myBlue}{dark blue}). 
In the top of figure \ref{sim:fig:caustic_QA} we plot the corresponding determinants as functions of time $q(t;x_0)$. $q(t; x_0)$ are linear increasing until the rays reach the lens, for $t \approx 8$ [$\mu$s].
Thereafter, the determinants $q(t; x_0)$ for different rays diverge, with $q(t; x_0)$ for the three \textcolor{myGreen}{light green} trajectories bending downwards and eventually becoming negative. As $q(t; x_0)$ approach 0, the individual amplitudes blow up and become unbounded for $q(t; x_0)=0$, as illustrated in the bottom of figure \ref{sim:fig:caustic_QA}. Numerically, these amplitudes will be set to 0 for all $t:\; |q(t; x_0)| < \epsilon$ thus this interval will not contribute to the solution at the sensor. Once $|q(t; x_0)| \geq \epsilon$, these amplitudes will again contribute to the solution.  
%The values where $q<0$ the amplitude is not well defined.

\begin{figure}
\begin{center}
\subfloat[]{\includegraphics[width = 6.8cm]{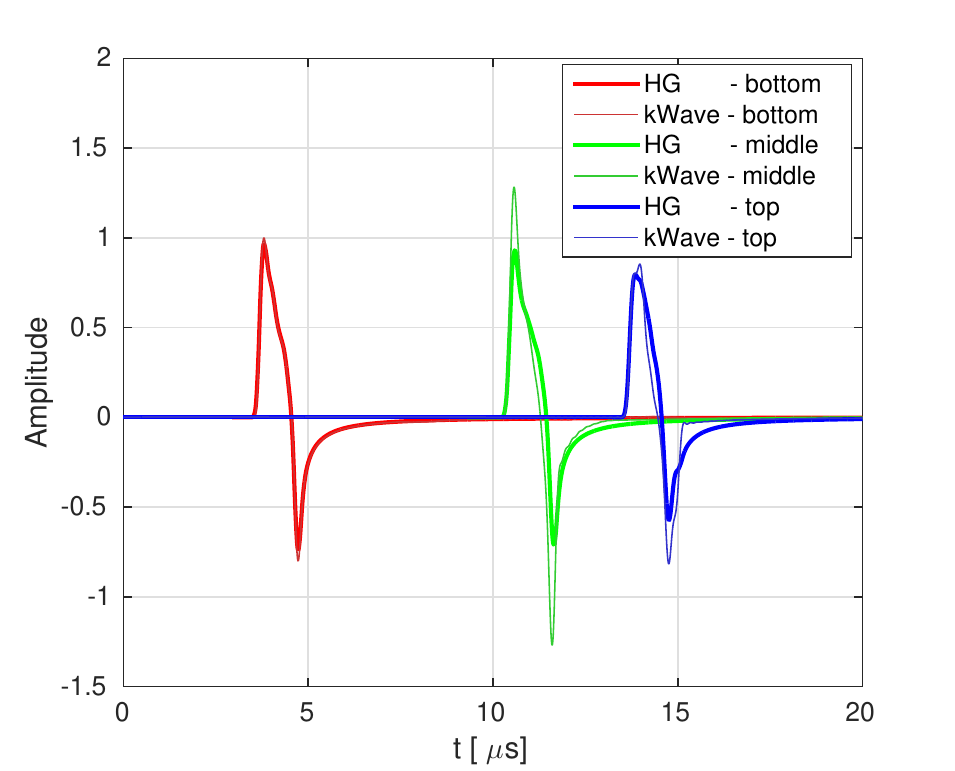}\label{sim:fig:caustic_signals}}
\subfloat[]{\includegraphics[width = 6.8cm]{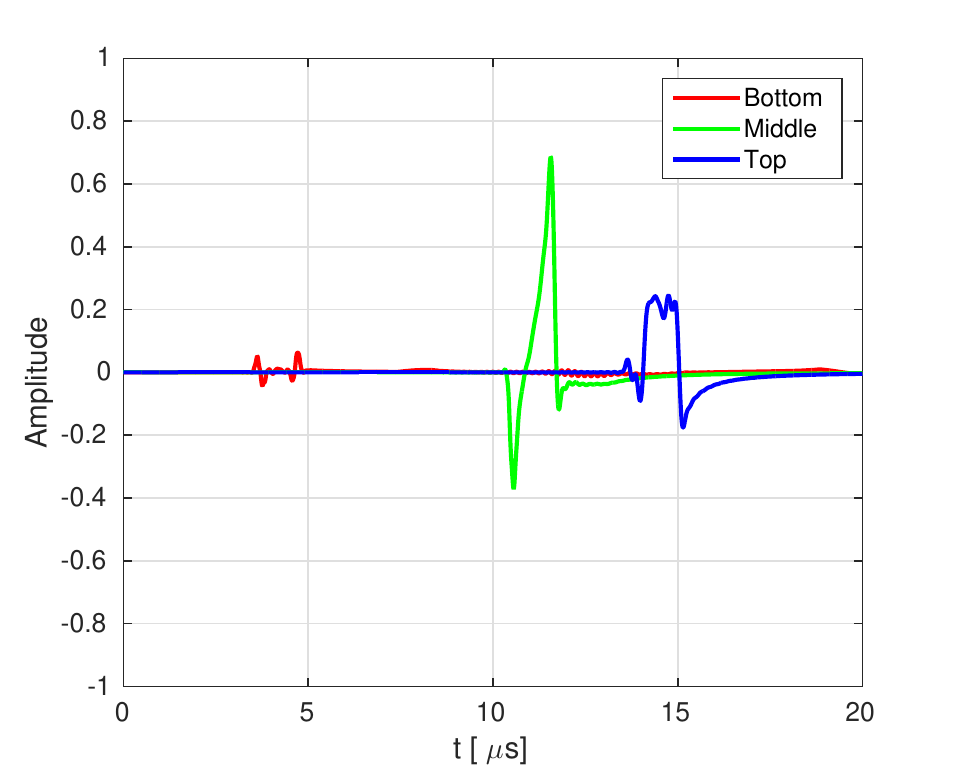}\label{sim:fig:caustic_error}}
\caption{(a) Forward solution $u(t, x_0)$ at $x_0$ for initial conditions $u_i, \; i\in \{b,m,t\}$ computed with HG and \texttt{k-Wave}. (b) Difference between HG and \texttt{k-Wave} solutions.}
%The error for the bottom ball, before the caustic, is in the order of 7\%, while the errors on and after the caustic are much larger.}
\label{}
\end{center}
\end{figure}}

\section*{Acknowledgements} 
We would like to acknowledge helpful discussions and encouraging feedback on the solver development from Ben Cox and Brad Treeby, 
on numerical solvers for time domain wave equation from Bill Symes, as well as on theory of ray approximations from Alden Waters. We thank the anonymous referees for their constructive criticism and helpful pointers which lead to the present version of the manuscript.

%==============================
% Bibliography
%==============================
\bibliographystyle{model1-num-names}
\bibliography{Bibliography}

\end{document}